\RequirePackage{fix-cm}
\documentclass[twocolumn]{svjour3}          % twocolumn
\smartqed  % flush right qed marks, e.g. at end of proof
\usepackage{hyperref}
\hypersetup{colorlinks,allcolors=blue}
\usepackage{graphicx}
\usepackage{amssymb}
\usepackage{amsmath}
\usepackage{float}
\usepackage{bm}
\usepackage{enumitem}
\usepackage[title]{appendix}
\usepackage[misc,geometry]{ifsym} 
\usepackage[caption=false]{subfig}

 \journalname{Nonlinear Dynamics}
\begin{document}

\title{Design of active network filters as hysteretic sensors %\thanks{Grants or other notes
}


\author{Yu Mao \and Harry Dankowicz %etc.
}


\institute{Y. Mao (\Letter) \and H. Dankowicz \at
              Department of Mechanical Science and Engineering, University of Illinois at Urbana-Champaign, Urbana, IL 61801, USA \\
            %   Tel.: +\\
            %   Fax: +\\
              \email{yumao4@illinois.edu}           %  \\
}

\date{Received: date / Accepted: date}

\maketitle

\begin{abstract}
This work aims to propose and design a class of networks of coupled linear and nonlinear oscillators, in which short bursts of exogenous excitation result in sustained endogenous network activity that returns to a quiescent state only after a characteristic time and along a different path than when originally excited. The desired hysteretic behavior is obtained through the coupling of self-excited oscillations with purposely designed rate laws for slowly-varying nodal parameters, governed only by local interactions in the network. The proposed architecture and the sought dynamics take inspiration from complex biological systems that combine endogenous energy sources with a paradigm for distributed sensing and information processing. In this paper, the network design problem considers arbitrary topologies and investigates the dependence of the desired response on model parameters, as well as on the placement of a single nonlinear node in an otherwise linear network. Perturbation analysis in various asymptotic parameter limits is used to define the proposed internal dynamics. Parameter continuation techniques validate the asymptotic results numerically and demonstrate their robustness over finite ranges of parameter values. Both approaches suggest a nontrivial dependence of the optimal distribution of nonlinearity on the network topology.

\end{abstract}

\section{Introduction}\label{sec1}

Using theory and computation, this paper considers the purposeful integration of self-excitability, bistability, and hysteresis in a network of coupled oscillators in order to obtain an \emph{active hysteretic sensor} that i) responds to excitation with sustained activity also after the excitation is removed and ii) resets autonomously after a refractory period through slow variations of internal parameters. This proposed design paradigm ensures that 
\begin{itemize}
    \item excitation at the appropriate frequency, even if applied only locally within the network, results in a significant global response;
    \item a memory of the response persists in the subsequent time history; and
    \item no further excitation or trigger is required to subdue this response after some characteristic time. 
\end{itemize}
A design problem of a similar qualitative nature is that envisioned in~\cite{MaMaFu2020} of a network of piezoelectric ultrasonic active sensors embedded throughout a structure and powered by harvested vibrational energy from the structure. In~\cite{UrYaYaMa2021}, self-excitability and parameter hysteresis are leveraged to develop a vibrational viscometer with improved performance relative to traditional sensor designs relying on frequency-response analysis. Coupling between a physical micro-cantilever and a virtual cantilever simulated by an analog circuit is used in~\cite{KaYaYaMa2020} to achieve high-sensitivity mass measurements of tens of nanograms.

In our system design, we take inspiration from ideas about self-organization in complex systems (e.g., task differentiation among social insects~\cite{duarte2011evolutionary}), the use of internal state variables to drive or explain critical transitions (e.g., in the dynamics of continental glaciers during the Pleistocene~\cite{engler2017dynamical}), and the role of active processes in enhancing sensor sensitivity (e.g., of human cochlea to input frequency~\cite{wang2017explosive}). We amalgamate these concepts in an original framework with rich mathematical structure and potential. We leave for future work a discussion of the utility of the proposed framework for characterizing naturally occurring dynamics or engineering design.

The relationship between local interactions among agents of a complex system and their collective behavior is of critical importance both for understanding natural systems and for designing engineered systems. In a study of bumble bees~\cite{CrGrMoKoOpPiCo2018}, the authors showed that spatial occupancy distributions play ``a functional role in patterning information flow within insect colonies.'' In~\cite{saghafi2017emergent}, local interactions in a static network designed to drive locally optimum-seeking behaviors were found associated with a potential destabilization of the global optimum and system collapse. Optimal protocols for randomized gossip-based algorithms on sparse peer-to-peer networks were studied in~\cite{ChPa2012} to minimize the time and number of interactions until all nodes had computed global properties of the system state. In our network design problem, we assume a static topology and restrict attention to deterministic information flows only between adjacent agents.

The components of our sensor networks are damped, linear mechanical oscillators and suitably modified van der Pol oscillators, linearly coupled by their relative displacements according to an arbitrary network topology. Through such coupling and under variations in the linear damping parameters, we obtain the possibility of either i) a globally stable trivial equilibrium, ii) a globally attracting limit cycle oscillation coexisting with an unstable trivial equilibrium, or iii) a locally attracting limit cycle oscillation coexisting with a locally stable trivial equilibrium. The phenomenology is richer than that available to a solitary van der Pol oscillator under one-parameter variations, even one designed to exhibit bistability~\cite{defontaines1990chain}, and depends in a non-trivial way on the network topology. Indeed, we find bifurcations of the trivial equilibrium associated with each of the linear modes of the network Laplacian with the possibility of multiple co-existing limit cycles and secondary limit cycle bifurcations.

In our previous work~\cite{mao2021topology}, we explored the response to steady harmonic excitation of networks of damped, linear mechanical oscillators and (regular) van der Pol oscillators for fixed values of the model parameters. Here, in contrast, we are principally concerned with the transient response of such networks under slowly-varying linear damping coefficients, allowing for exogenous excitation at most as an initial trigger of the subsequent system dynamics. By suitable design of the internal dynamics governing the damping coefficients, we leverage the full bifurcation structure available to these networks, and allow hysteretic transitions between self-excited oscillations and the trivial equilibrium. This achieves the sought sensor behavior.

It remains, of course, to define the internal dynamics governing the damping coefficients. In~\cite{JoScTh1997}, slow monotonic dynamics of the concentration of bromomalonic acid were shown to govern observed transitions between different oscillatory behaviors and quiescence in a model of the Belousov-Zhabotinsky reaction, even in a system with unlimited reactants. Transitions from stationary cutting to chatter in milling operations due to slowly-varying system parameters were considered in~\cite{DoMuKuSt2018}. A modeling study in~\cite{roongthumskul2011multiple} of spontaneous oscillations of hair bundles of the bullfrog sacculus interspersed with quiescent intervals were attributed to slow cyclic variations in the stiffness of the hair cells across a supercritical Hopf bifurcation. In this paper, we rely on bistability to design the internal dynamics of the damping coefficients to result in a delayed transition back to quiescence after an initial trigger, and with an autonomous return of these coefficients to an equilibrium value.

Coupled oscillators as platforms for sensor design have been widely studied~\cite{CaLiMaAr2017,HuHo2020,LuAdYaLeUr2016,MaPhCr2019}. A common target of measurement is added mass deposited on the mechanism~\cite{spletzer2006ultrasensitive} resulting in a shift in the resonance frequency~\cite{MaAlPlViDoLe2017} or amplitude ratio at resonance~\cite{WeZhPuBoSaKr2016}. Similar phenomenology is used in~\cite{ZhNi2020} to detect the presence of specific biomolecules through their effect on the dielectric properties of the surrounding medium. In~\cite{PyBaRhWeQu2019}, a network of Colpitts oscillators is proposed as a potential threshold color sensor. A network of coupled self-oscillating Belousov-Zhabotinsky gels and piezoelectric MEMS units was proposed in~\cite{YaLeBa2018} for the design of materials able to perform ``computational tasks such as pattern recognition.'' Agnostic to a particular source of stimuli, the analysis in the current paper is concerned with a particular form of transient system response, triggered by an excitation exceeding some critical threshold.

The remainder of the paper is organized as follows. In Sect.~\ref{sec:a model problem}, we study the dynamics of a simple 4-node network as an accessible precursor to a general discussion of arbitrary networks. The analysis aims to build intuition and guide the construction of autonomous rate laws for the nodal damping coefficients that accomplish the desired hysteretic response to a triggering excitation. We use the multiple-scale perturbation method to derive rigorous results in the asymptotic limits of very small or very large damping and use these results to parameterize the internal dynamics. Importantly, in these limits, we show that local interactions within the network provide sufficient input to the nodal rate laws. In Sect.~\ref{sec_networkfilters}, we generalize these results to networks of arbitrary topology but with only a single nonlinear node. An extension of the asymptotic theory to the case of finite damping is proposed in Sect.~\ref{sec:examples} and validated using numerical simulations on a generic 15-node network example. These results show that the desired hysteretic behavior is possible for a finite range of values of a common scaling factor. Parameter continuation is used here and in Sect.~\ref{sec_arameterrobustness} to characterize this range as a function of the network topology and location of the nonlinear node. To complete the construction, Sect.~\ref{sec_hysteresisresonance} briefly considers the relationship between the magnitude and duration of a burst of harmonic excitation required to trigger the desired hysteretic behavior. Finally, a concluding discussion is presented in Sect.~\ref{sec_conclusion}. 

%%%%%%%%%%%%%%%%%%%%%%%%%%%%%%%%%%%%%%%%%%%%%%%%%%%%%%%%%%%%%%%%%%%%%%%

\section{A 4-node network model design}\label{sec:a model problem}

As explained in the introduction, we use networks of coupled mechanical oscillators as a template upon which additional dynamics may be imposed. We aim to design an \emph{active hysteretic sensor}, i.e., a dynamical system that i) responds to excitation by transitioning from a quiescent state to one with sustained activity, ii) continues to exhibit such activity also after the excitation is removed, and iii) resets autonomously to the state of quiescence after a refractory period and along a different path than when excited.

\begin{figure}[htbp]
    \centering
    \includegraphics[width=.4\textwidth]{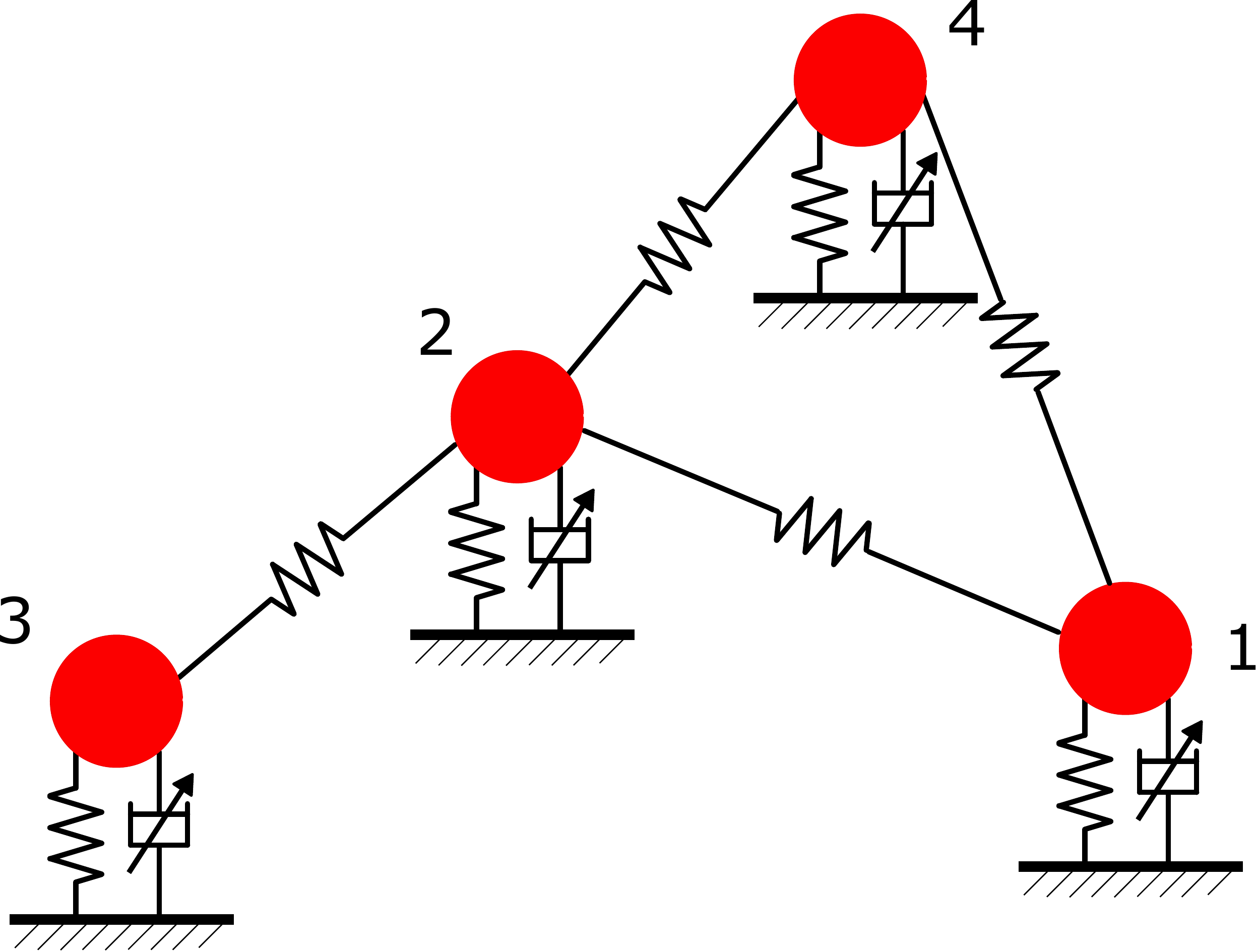}
    \caption{A schematic representation of a four-node network of linearly coupled, linear and weakly nonlinear oscillators with all vertical displacements.}
    \label{4node_network_example}
\end{figure}

In this section, we illustrate the main characteristics of this design problem by considering the dynamics of a linearly coupled, mechanical oscillator network with four-node topology shown in Fig.~\ref{4node_network_example} and with displacement vector $u$ governed by the differential equation 
\begin{equation}\label{eq_4nodegenform}
    \ddot{u}+C(u)\dot{u}+Ku=F
\end{equation}
in terms of the diagonal damping matrix $C(u)$, stiffness matrix
\begin{equation}\label{eq_4nodestiffness}
    K=\begin{pmatrix}3 & -1 & 0 & -1\\-1 & 4 & -1 & -1\\0 & -1 & 2 & 0\\-1 & -1 & 0 & 3\end{pmatrix},
\end{equation}
and time-dependent exogenous excitation vector $F(t)$. As shown here and in subsequent sections, this four-node network captures all the essential ingredients of our network design problem and allows us to develop a paradigm for implementing the desired hysteretic behavior that generalizes to networks of arbitrary size.

We restrict attention to the four possible designations of a single node as a nonlinear oscillator with damping coefficient
\begin{equation}
    \epsilon\left(-1-10u_Q^2+10u_Q^4\right),\,Q\in\{1,2,3,4\}
\end{equation}
with the remaining nodes assumed to be linear with damping coefficients $\epsilon\zeta_i$, $i\ne Q$. For example, the choice \begin{equation}\label{eq_4nodedamping}
    C(u)=\epsilon\begin{pmatrix}-1-10u_1^2+10u_1^4 & 0 & 0 & 0\\0 & \zeta_2 & 0 & 0\\0 & 0 & \zeta_3 & 0\\0 & 0 & 0 & \zeta_4\end{pmatrix},\,\epsilon,\zeta_i>0
\end{equation}
implies that nodes 2 through 4 are linear oscillators with positive damping, while node 1 has a nonlinear damping coefficient that is negative for small nodal displacements, increasingly negative for medium-sized displacements, and increasingly positive for large nodal displacements. To allow us to deploy analytical techniques in the analysis, we concern ourselves here with small values of the parameter $\epsilon$. We model the imposition of additional dynamics through possibly time-varying $\zeta_i(t)$ obtained, for example, from internal dynamics $\dot{\zeta_i}=g_i(u,\dot{u},\zeta_i)$ for some set of functions $g_i$.

Our design problem is then reduced to a selection, as appropriate, of the network topology, including the value of $Q$, as well as a determination of a choice of the functions $g_i$ that is mathematically and physically realizable. Importantly, we aim for $g_i$ to depend only on information about $u$ and $\dot{u}$ that is available through nearest neighbor interactions in the network.

\subsection{Small linear damping} \label{sec:small linear damping}

Before proceeding to impose time-dependence on the damping parameters, we use perturbation analysis to characterize the autonomous dynamics for constant $\zeta_i$ and $F=0$.

\subsubsection{Periodic response}

The stability of the trivial equilibrium $u=0$ in the absence of exogenous excitation depends on a competition between the stabilizing effect of the positive damping on the linear oscillators and the destabilizing effect of the negative damping for small-to-medium-sized displacements of the nonlinear oscillator. For the choice of $C(u)$ in Eq.~\eqref{eq_4nodedamping} and constant $\zeta_2,\zeta_3,\zeta_4\sim\mathcal{O}(1)$, analysis of the first-order form of the governing equations in $(u_1,\dot{u}_1,\ldots,\allowbreak u_4,\dot{u}_4)$ finds the local dynamics governed by eight complex exponential rates
\begin{equation}\label{eq_4nodeeigvalue}
\begin{split}
    \lambda_{1,2}&=\pm\mathrm{j}+\frac{1-\zeta_2-\zeta_3-\zeta_4}{8}\epsilon+\mathcal{O}(\epsilon^2),\\
    \lambda_{3,4}&=\pm\sqrt{2}\mathrm{j}+\frac{1-4\zeta_3-\zeta_4}{12}\epsilon+\mathcal{O}(\epsilon^2),\\
    \lambda_{5,6}&=\pm2\mathrm{j}+\frac{1-\zeta_4}{4}\epsilon+\mathcal{O}(\epsilon^2),\\
    \lambda_{7,8}&=\pm\sqrt{5}\mathrm{j}+\frac{1-9\zeta_2-\zeta_3-\zeta_4}{24}\epsilon+\mathcal{O}(\epsilon^2),
\end{split}
\end{equation}
where $\mathrm{j}^2=-1$. In this approximation, it follows that the trivial equilibrium is asymptotically stable for $\zeta_4>1$ and unstable for $\zeta_4<1$. The transition at $\zeta_4=1$ corresponds to a Hopf bifurcation associated with a branch of periodic orbits with limiting frequency $2$ and period $\pi$ corresponding to the crossing of the imaginary axis by the eigenvalue pair $\lambda_{5,6}$. Since the corresponding eigenvectors are parallel to $v_{5,6}=\begin{pmatrix}\pm\mathrm{j} & -2 & 0 & 0 & 0 & 0 \mp\mathrm{j} & 2\end{pmatrix}^\top$, it follows that the periodic orbits are approximated in the small-amplitude limit by the normal-mode oscillations
$u_1(t)=-u_4(t)=\left(A/\sqrt{2}\right)\cos(2t+\phi)$, $u_2(t)=u_3(t)=0$ for constant amplitude $A$ and phase $\phi$.

\begin{figure}[htbp]
\centering
\includegraphics[width=.48\textwidth]{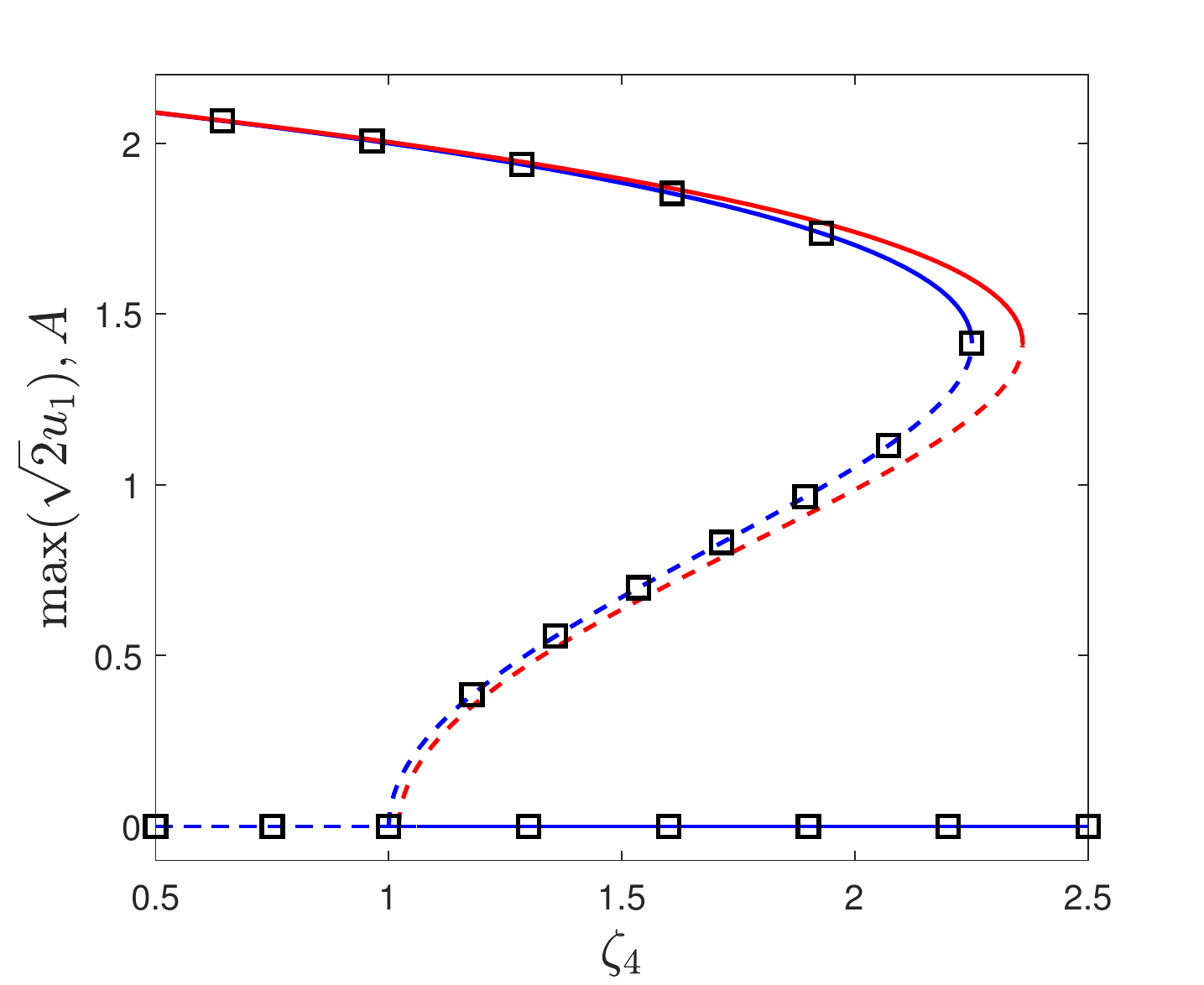}
\caption{Variations of $\max(\sqrt{2}u_1)$ in the 4-node network for $\epsilon=0.01$ (blue lines) and $\epsilon=0.1$ (red lines) obtained using numerical continuation along the branch of periodic orbits emanating from the Hopf bifurcation at $\zeta_4=1$. Solid lines denote stable periodic responses and dashed lines denote unstable periodic responses. Black squares represent the slow-flow approximation for the amplitude $A$ given by Eq.~\eqref{eq_4nodeslowflowamp}.}\label{4_node_epvar_zeta_vs_A}
\end{figure}

As detailed in Appendix~\ref{appendix_4node_multiplescale}, we may use the method of multiple scales \cite{nayfeh1995nonlinear} to determine the fate of the periodic orbits as their amplitude grows with deviations of $\zeta_4$ from $1$. To this end, substitution of
\begin{equation}\label{eq_4nodemultscale}
\begin{split}
    &u_1(t)=\frac{1}{\sqrt{2}}A(\epsilon t)\cos\left( 2t+\phi(\epsilon t )\right)+\epsilon v_1(t),\\
    &u_2(t)=\epsilon v_2(t),\,u_3(t)=\epsilon v_3(t),\\
    &u_4(t)=-\frac{1}{\sqrt{2}}A(\epsilon t)\cos\left( 2t+\phi(\epsilon t )\right)+\epsilon v_4(t)
\end{split}
\end{equation}
into the fully nonlinear governing equations and elimination of secular terms $\cos(2t+\phi)$ and $\sin(2t+\phi)$ in the dynamics of the difference $\left(v_1-v_4\right)/\sqrt{2}$ (with natural frequency $2$) yield the conditions
\begin{equation}\label{eq_4nodeslowflow1}
    4 A\phi' = 0
\end{equation}
and
\begin{equation}\label{eq_4nodeslowflow2}
    4 A'-(1-\zeta_4)A-\frac{5}{4}A^3+\frac{5}{16}A^5 = 0
\end{equation}
(here and below, $'$ denotes differentiation with respect to the slow time scale $\epsilon t$). These conditions are satisfied at nontrivial equilibrium values (in the slow time scale) of $A$ and $\phi$ provided that
\begin{equation}\label{eq_4nodeslowflowamp}
    A^2=2\pm 2\sqrt{\frac{9-4\zeta_4}{5}},
\end{equation}
i.e., for $\zeta_4\in[0,9/4]$ with two co-existing solutions on the interval $\zeta_4\in[1,9/4)$. In this asymptotic limit, we conclude that the Hopf bifurcation at $\zeta_4=1$ is subcritical and that the branch of periodic orbits has a geometric fold at a saddle-node bifurcation at $\zeta_4=9/4$, as shown in Fig.~\ref{4_node_epvar_zeta_vs_A} (black squares). We validate the predictions of the perturbation analysis using numerical continuation along the branch of periodic solutions emanating from the Hopf bifurcation at $\zeta_4=1$ for different values of $\epsilon$. The results (solid and dashed curves) in Fig.~\ref{4_node_epvar_zeta_vs_A} show close agreement for small $\epsilon$.

\begin{figure}[htbp]
\centering
\includegraphics[width=.48\textwidth]{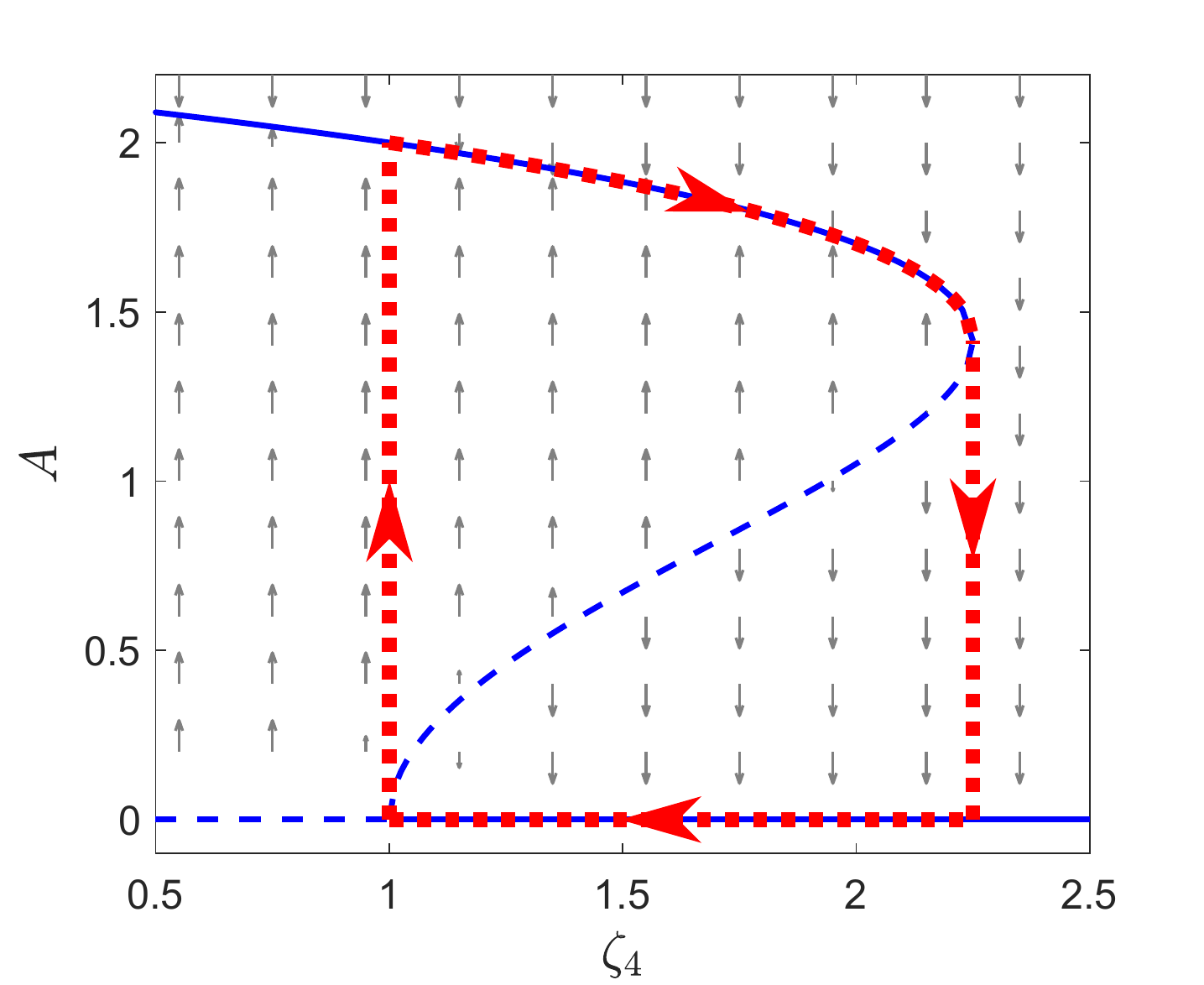}
\caption{Stable (solid) and unstable (dashed) steady state value of $A(\epsilon t)$ in the 4-node network in the slow-flow approximation given by Eq.~\eqref{eq_4nodeslowflow2} for different values of the damping at node 4. Grey arrows denote the vector field of the $A$ dynamics. Red dotted curves denote a possible hysteresis loop. }\label{4_node_zeta_vs_A}
\end{figure}

\subsubsection{Hysteretic response}

We reproduce the predictions of the asymptotic analysis in Fig.~\ref{4_node_zeta_vs_A} and include arrows (grey) describing the direction of flow of $A(t_1)$ per \eqref{eq_4nodeslowflow2} away from equilibrium solutions. Indeed, for $\zeta_4<1$, it follows from Eq.~\eqref{eq_4nodeslowflow2} that arbitrarily small perturbations $A(0)$ from the trivial equilibrium result in transient dynamics that asymptotically converge to a limit-cycle oscillation in the displacement vector along the upper branch of Fig.~\ref{4_node_zeta_vs_A}. For $\zeta_4\in[1,9/4)$ such perturbations have to exceed a critical threshold in order to trigger a transition to this branch of self-excited oscillatory dynamics. Clearly, the trivial equilibrium is the globally asymptotically-stable steady-state behavior for $\zeta_4>9/4$ in the absence of excitation. The existence of a region of bistability for $\zeta_4\in[1,9/4)$ suggests the possibility of a hysteretic system response under a prescribed, slow, and cyclic variation in the damping coefficient $\zeta_4$ over an interval containing both the Hopf and saddle-node bifurcation points, as suggested by the red dotted curves in Fig.~\ref{4_node_zeta_vs_A}. As desired, this response shows sustained activity even after the triggering excitation has been removed and resets to the trivial equilibrium after a refractory period and along a different path than when excited.

\begin{figure}[htbp]
\centering
\includegraphics[width=.48\textwidth]{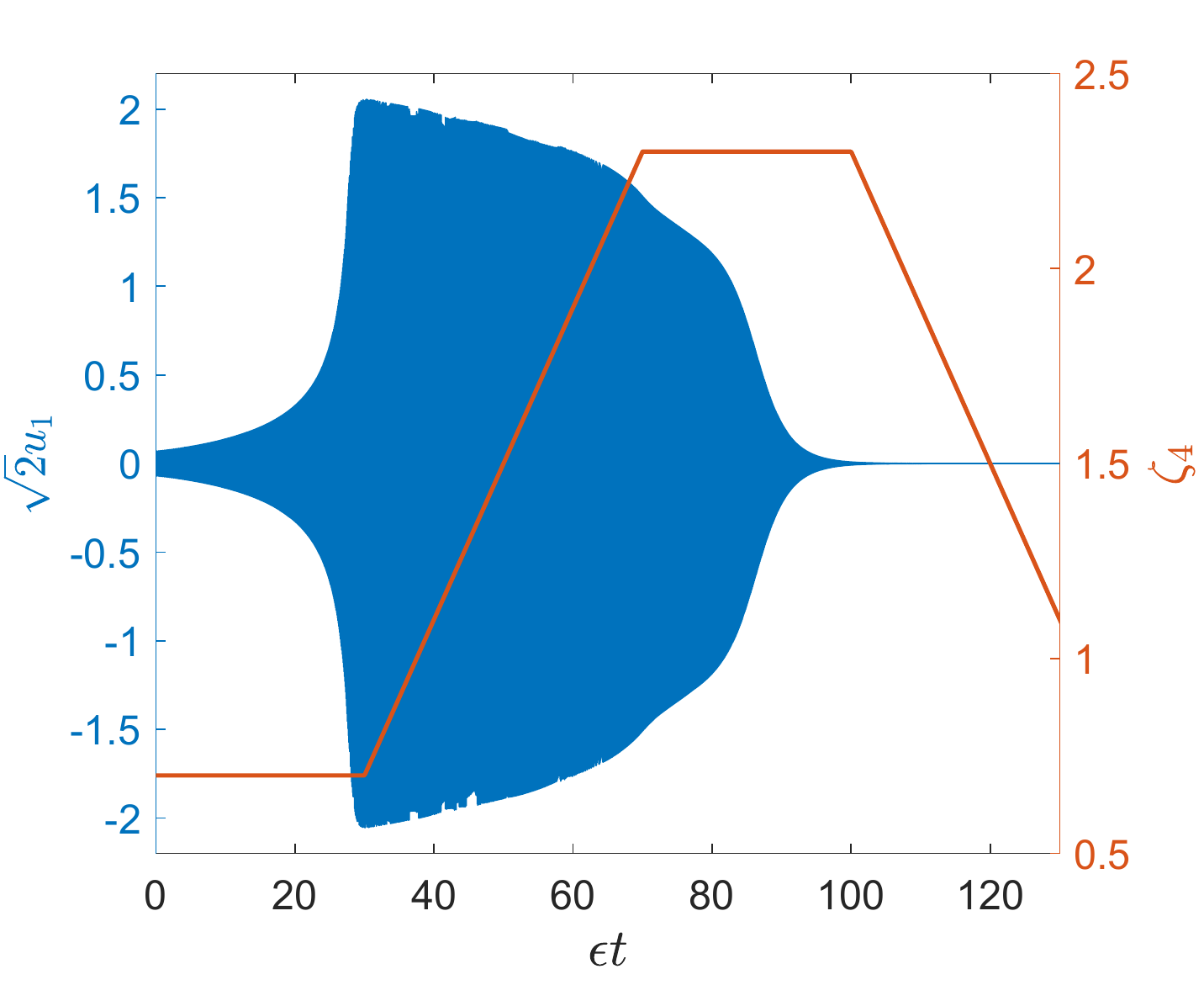}
\caption{Time histories of the displacement of node 1 (blue line) in the 4-node network with $\epsilon=0.01$ and varying $\zeta_4$ (orange line). Small random noise is added to $u(0)$ to trigger the hysteretic response.}\label{4_node_example_Q1_amp_ep0_1}
\end{figure}

As shown in Fig.~\ref{4_node_example_Q1_amp_ep0_1}, these predictions carry over to the full dynamics of the 4-node network system in the presence of small amount of noise. Here, a slow, cyclic, piecewise-linear variation of $\zeta_4$ (red) triggers a rapid transition to large oscillations in the displacement of node 1 (blue) that persist over an extended period of time (relative to the period of oscillation) and reset to the quiescent state for a different value of $\zeta_4$.

\subsubsection{Damping dynamics}

\begin{figure}[htbp]
\centering
\includegraphics[width=.48\textwidth]{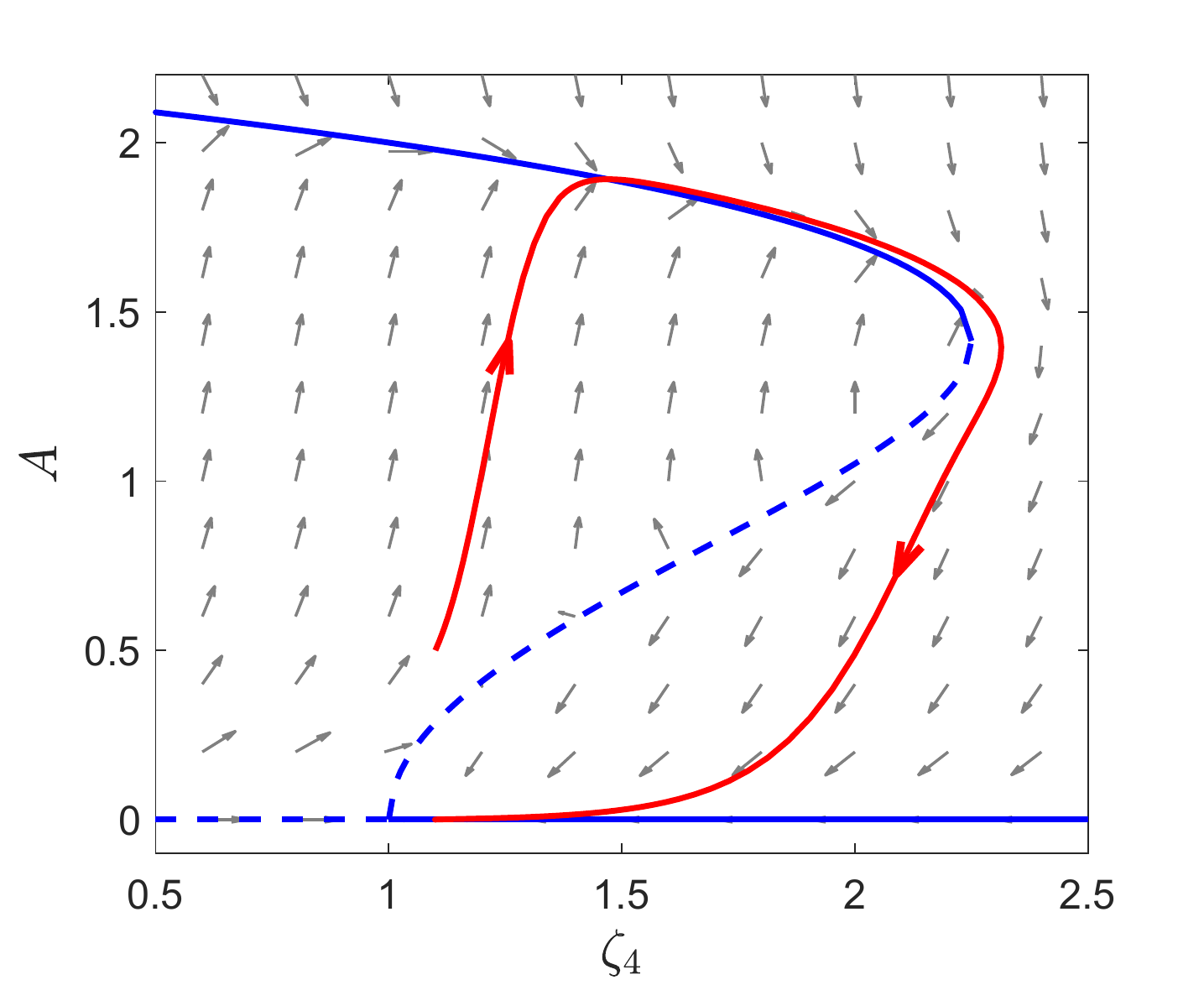}
\caption{A hysteresis loop (red line) governed by Eqs.~\eqref{eq_4nodeslowflow2} and \eqref{eq_4node_damplaw1} with $\delta=0.1$, $\tau=20$, $\zeta_4(0)=1.1$, and $A(0)=0.5$. Grey arrows denote the vector field of the coupled $A$ and $\zeta_4$ dynamics.}\label{4_node_network_hysteretic_example}
\end{figure}

\begin{figure}[htbp]
\centering
\includegraphics[width=.48\textwidth]{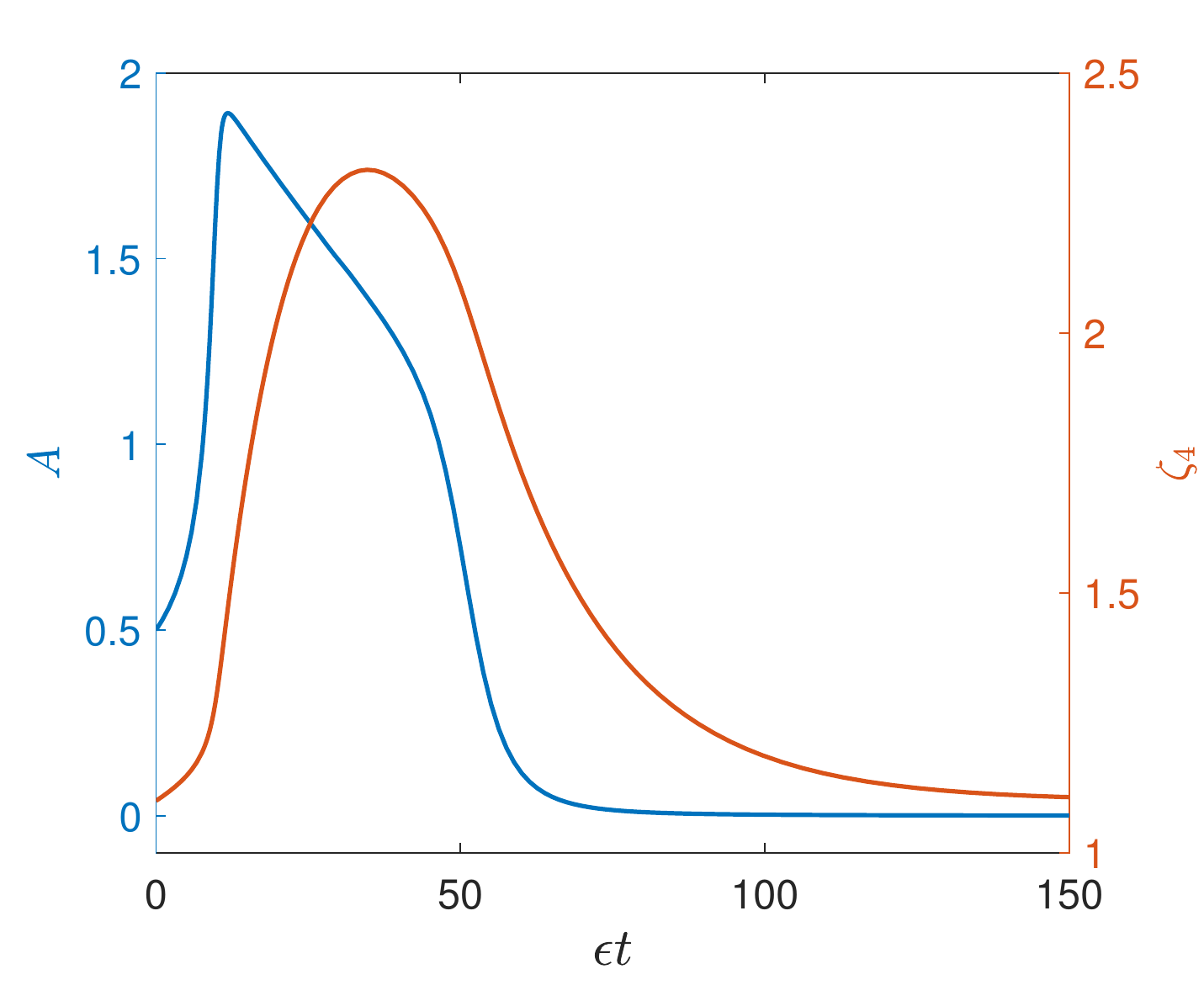}
\caption{Time histories of $A$ (blue line) and $\zeta_4$ (orange line) along the hysteresis loop in Fig.~\ref{4_node_network_hysteretic_example} versus the slow time variable $\epsilon t$. Comparison with Fig.~\ref{4_node_example_Q1_amp_ep0_1} shows the desired characteristics. }\label{4_node_network_hysteretic_example_time}
\end{figure}

As an alternative to a prescribed variation in the damping coefficient $\zeta_4$, we may imagine that such variations would be autonomously triggered by a transition of the network system from the locally-stable trivial equilibrium to the self-excited oscillation, induced by a short burst of exogenous excitation applied to the network when $\zeta_4$ is close to, but to the right of, the Hopf bifurcation. Consider, for example, a model coupling the vanishing of Eq.~\eqref{eq_4nodeslowflow2} with the damping dynamics
\begin{equation}\label{eq_4node_damplaw1}
    \zeta'_4=\tau^{-1}\left(-\zeta_4+\delta+1+\frac{5}{8}A^2\right)
\end{equation}
for some appropriately large time scale $\tau>0$ and some appropriately small offset $\delta>0$. From this equation, we expect $\zeta_4$ to be a slowly increasing function of time along the stable branch of limit-cycle oscillations, ultimately resulting in values of $\zeta_4$ that exceed the value at the saddle-node bifurcation at $9/4$ (where $A=\sqrt{2}$). Once the dynamics again converge to the trivial equilibrium, $\zeta_4$ decreases slowly toward $1+\delta$ allowing enough time for the system to reset. An example of such a hysteretic response is shown in Figs.~\ref{4_node_network_hysteretic_example} and \ref{4_node_network_hysteretic_example_time} for the case that $\delta=0.1$, $\tau=20$, $\zeta_4(0)=1.1$, and $A(0)=0.5$. In contrast to the schematic in Fig.~\ref{4_node_zeta_vs_A} and due to the finite value of $\tau$, the trajectory overshoots the upper branch of the $A$ nullcline and drops back toward $A=0$ only after $\zeta_4$ moves well past the saddle-node point.

\begin{figure}[htbp]
\centering
\includegraphics[width=.48\textwidth]{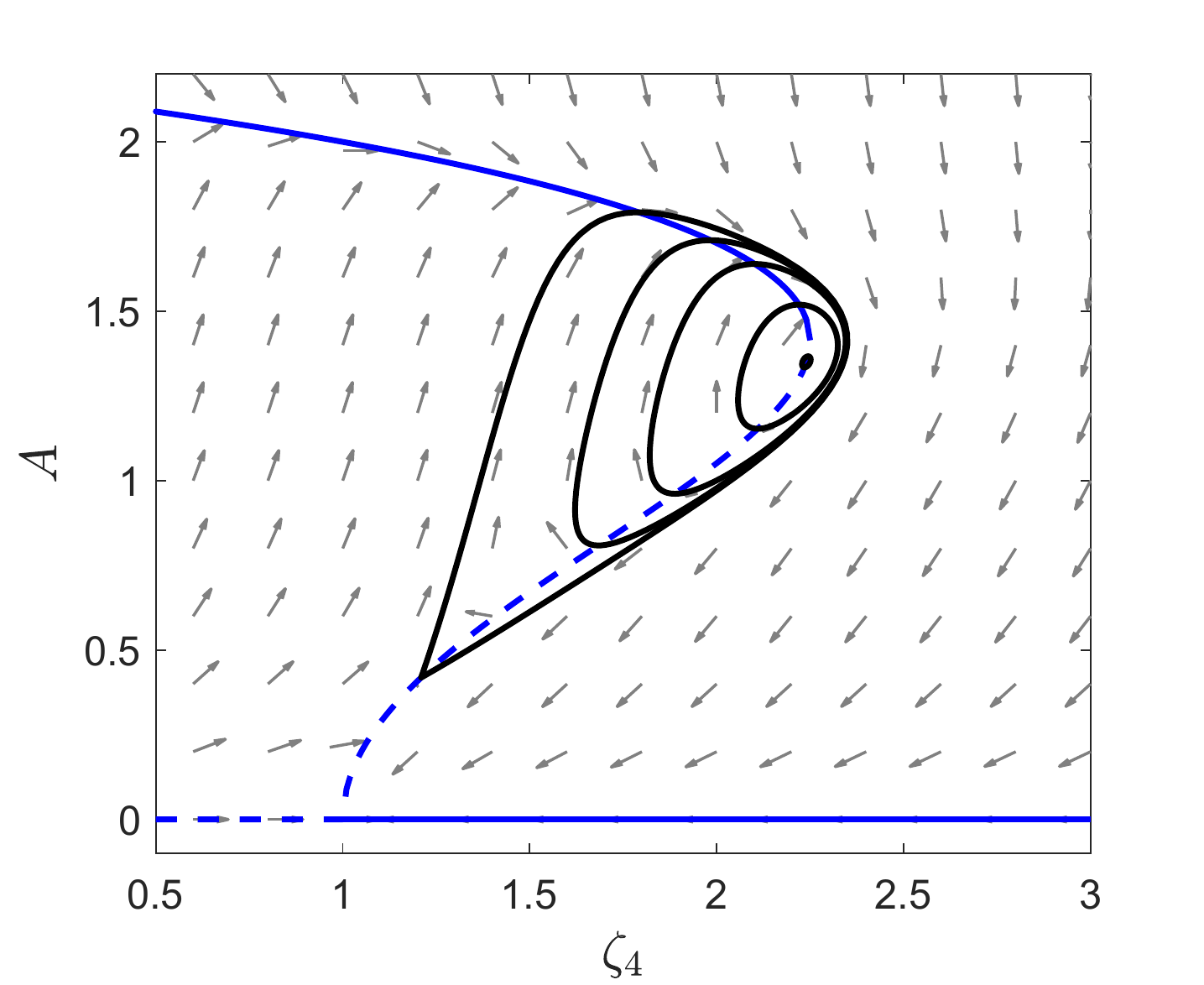}
\caption{The family of periodic orbits (black lines) born at a Hopf bifurcation with different value of $\tau$ increasing from $1/\delta=10$ and terminating on a homoclinic orbit. Grey arrows denote the vector field of the coupled $A$ and $\zeta_4$ dynamics with $\tau=10$.}\label{HB_PO_small_damping}
\end{figure}

While Eq.~\eqref{eq_4node_damplaw1} appears to achieve the desired behavior, it behooves us to explore the dependence on the model parameters $\tau$ and $\delta$. Indeed, for arbitrary $\tau$ and $\delta$, the coupled dynamics has a trivial equilibrium at $(\zeta_4,A)=(1+\delta,0)$. Additional equilibria are found along the middle branch of the $A$ nullcline, where
\begin{equation}
    \zeta_4=\delta+\frac{13\pm\sqrt{25-80\delta}}{8}
\end{equation}
provided that $\delta\le 5/16$. These equilibria emerge from the end points of the middle branch of the $A$ nullcline as $\delta$ increases from $0$ and disappear at a saddle-node bifurcation when $\delta=5/16$, $\zeta_4=31/16$, and $A=1$. Indeed, while the trivial equilibrium is asymptotically stable for any $\tau$ and $\delta$, the nontrivial equilibrium with $A<1$ is always a saddle and the nontrivial equilibrium with $A>1$ is stable for $\tau<1/\delta$ and unstable for $\tau>1/\delta$. The bifurcation at $\tau=1/\delta$ is always a Hopf bifurcation, out of which emanates a ``vertical'' canard family of limit cycles, with representatives shown in Fig.~\ref{HB_PO_small_damping} in the case that $\delta=0.1$. These in turn limit on a homoclinic orbit at a homoclinic bifurcation of the nontrivial saddle equilibrium with $A<1$ and $\tau	> 1/\delta$.

\begin{figure}[htbp]
\centering
\includegraphics[width=.48\textwidth]{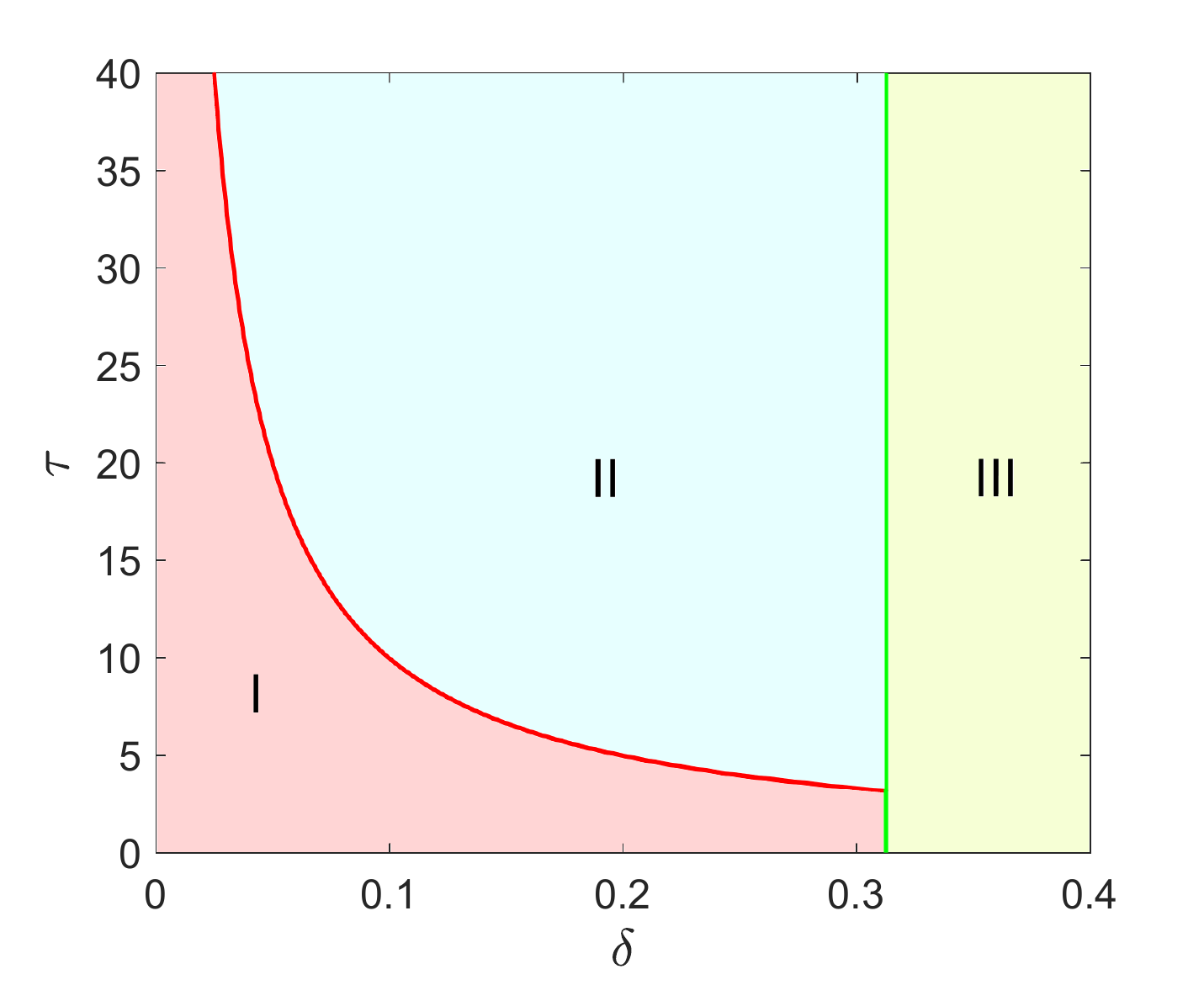}
\caption{Hopf (red curve) and saddle-node (green curve) bifurcations for nontrivial equilibria of the dynamics governed by Eqs.~\eqref{eq_4nodeslowflow2} and \eqref{eq_4node_damplaw1} in the $(\delta,\tau)$ space. }\label{HB_and_SN_small_damping}
\end{figure}

\begin{figure}[htbp]
\centering
\includegraphics[width=.48\textwidth]{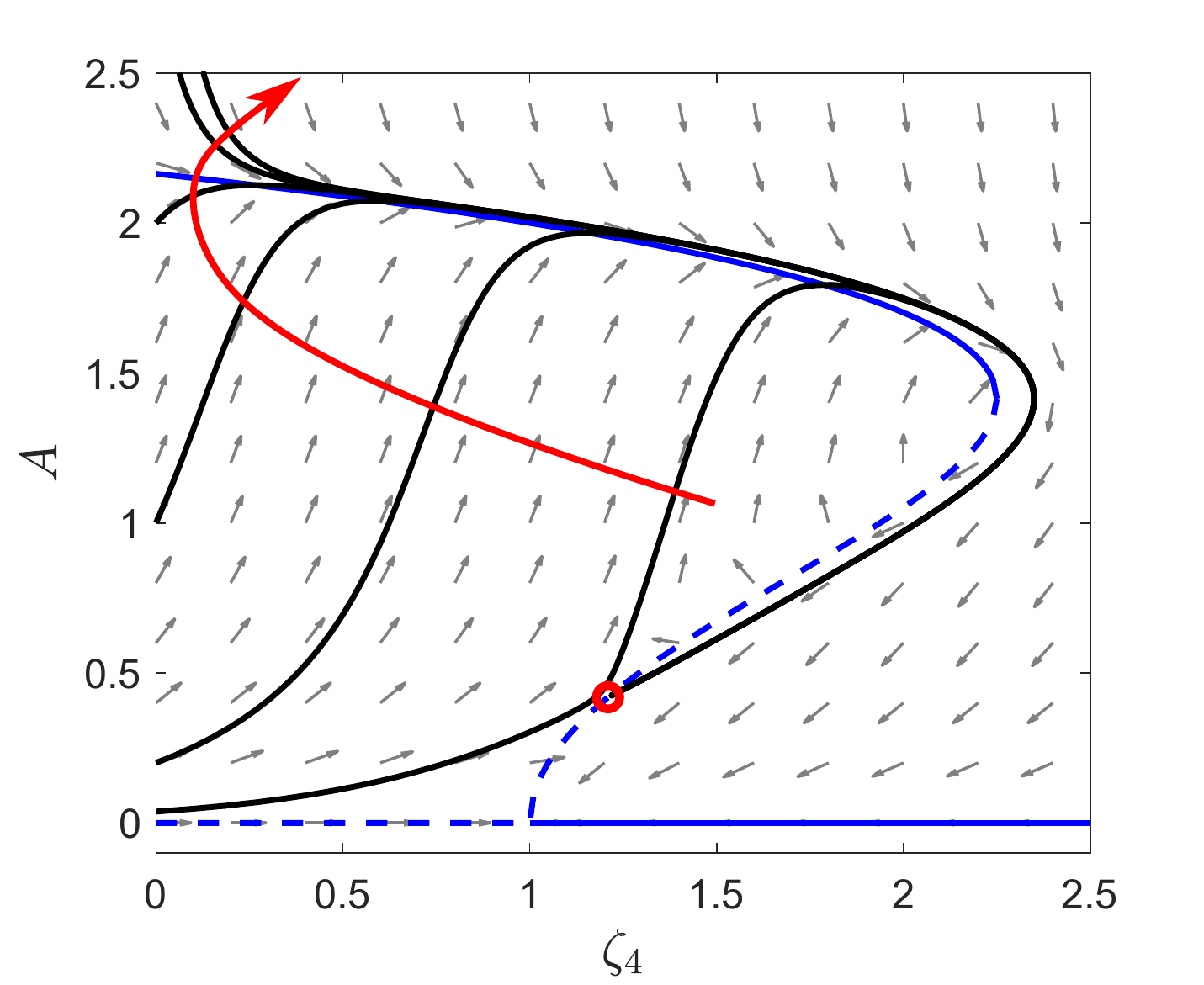}
\caption{One branch of the stable manifold (black) of the saddle equilibrium (red circle) in the 4-node network example with different values of $\tau$ decreasing from $1/\delta=10$. Grey arrows represent the vector field of the coupled $A$ and $\zeta_4$ dynamics with $\tau=10$, and the red arrow indicates increasing values of $\tau$.}\label{4_node_stable_manifold}
\end{figure}

We summarize these observations regarding the nontrivial equilibrium with $A>1$ in the bifurcation diagram in Fig.~\ref{HB_and_SN_small_damping}. This provides some insight into the selection of $\tau$ to ensure the desired hysteretic behavior that eventually returns to the trivial equilibrium following a triggering event. In particular, for $(\delta,\tau)$ in region I the equilibrium with $A>1$ is stable, while it is unstable for $(\delta,\tau)$ in region II and absent for $(\delta,\tau)$ in region III. Either of the latter regions may provide for a hysteretic behavior. For $(\delta,\tau)$ in region I, the burst of exogenous excitation may result in a subsequent autonomous trajectory in the basin of attraction of the stable nontrivial equilibrium, rather than the trivial equilibrium. The corresponding basins are separated by one branch of the stable manifold of the nontrivial equilibrium with $A<1$. At the homoclinic bifucation, this coincides with one branch of the unstable manifold. Careful numerical analysis illustrated in Fig.~\ref{4_node_stable_manifold} shows this branch sweeping upward rapidly as $\tau$ decreases from $1/\delta$, quickly allowing the basin of attraction to the nontrivial equilibrium to encompass the entire region under the upper branch of the $A$ nullcline.

A similar analysis applies to other configurations of the network. For example, as shown in Appendix~\ref{appendix_4node_smalldamping}, if the nonlinear oscillator is located at node 2 with $\zeta_1,\zeta_3,\zeta_4\sim\mathcal{O}(\epsilon)$, the multiple-scales analysis yields the differential equation
\begin{equation}
    A'=\frac{1}{24}(9-\zeta)A+\frac{45}{64}A^3-\frac{135}{512}A^5,
\end{equation}
where $\zeta:=\zeta_1+\zeta_3+\zeta_4$. This governs the slow dynamics of the amplitude of periodic oscillations born at a Hopf bifurcation at $\zeta=9$. Nonzero equilibria of $A$ are found along solutions to
\begin{equation}
    A^2=\frac{4}{3}\pm\frac{4}{9}\sqrt{\frac{81-4\zeta}{5}},
\end{equation}
i.e., for $\zeta\in[0,81/4]$ with two co-existing solutions on the interval $\zeta\in[9,81/4)$. In particular, at the non-trivial equilibrium with $\zeta=81/4$, we find $A=2/\sqrt{3}$ and
\begin{equation}
    -\zeta+9+\frac{135}{16}A^2=0.
\end{equation}
Consequently, we obtain the desired hysteretic response by letting 
\begin{equation}
    \zeta'=\tau^{-1}\left(-\zeta+\delta+9+\frac{135}{16}A^2\right)
\end{equation}
for some large time scale $\tau>0$ and some small offset $\delta>0$. From the definition of $\zeta$, we see that this holds if
\begin{equation}
\label{eq:damp_node2}
\zeta'_{1,3,4}=\tau^{-1}\left(-\zeta_{1,3,4}+\frac{1}{3}\delta+3+\frac{45}{16}A^2\right).
\end{equation}

Although, in all cases, we may obtain the desired hysteretic response through the appropriate coupling of some equivalent damping coefficient $\zeta$ to the amplitude $A$, we deem this scheme \emph{unrealizable} for the full network, since $A$ does not directly correspond to a measurable state of the system. A further complication arises in the case of large linear damping, to which we now turn.

\subsection{Large linear damping}\label{sec:large linear damping}

Recall, again, the choice of $C(u)$ in Eq.~\eqref{eq_4nodedamping} but this time assuming large linear damping with $\zeta_2,\zeta_3,\zeta_4\sim\mathcal{O}(1/\epsilon^2)$. In this limit, the local dynamics near the $u=0$ equilibrium are governed by the exponential rates
\begin{equation}
\begin{split}
    &\lambda_{1,2,3}\approx -\epsilon\zeta_{2,3,4}+\mathcal{O}(1),\\
    &\lambda_{4,5,6}\lesssim 0, \\
    &\lambda_{7,8}=\pm\mathrm{j}\sqrt{3}+\frac{1}{6}\left(3-\tilde{\zeta}_2-\tilde{\zeta}_4\right)\epsilon+\mathcal{O}(\epsilon^2),
\end{split}
\end{equation}
where $\tilde{\zeta}_2:=1/\epsilon^2\zeta_2$ and $\tilde{\zeta}_4:=1/\epsilon^2\zeta_4$. With $\tilde{\zeta}:=\tilde{\zeta}_2+\tilde{\zeta}_4$, it follows that the trivial equilibrium is asymptotically stable for $\tilde{\zeta}>3$ and unstable for $\tilde{\zeta}<3$, with $\tilde{\zeta}=3$ corresponding to a Hopf bifurcation out of which emanates a branch of periodic orbits with limiting frequency $\sqrt{3}$ and period $2\pi/\sqrt{3}$. Since the corresponding eigenvectors are parallel to $v_{7,8}=(1,\pm\mathrm{j}\sqrt{3},0,\allowbreak 0,0,0,0,0)$, it follows that the periodic orbits are approximated in the small-amplitude limit by the normal-mode oscillations $u_1=A\cos(\sqrt{3}t+\phi)$, $u_2(t)=u_3(t)=u_4(t)=0$ for constant amplitude $A$ and phase $\phi$. A consistent multiple-scale ansatz yields the slow-flow amplitude equation
\begin{equation}
    A'=\frac{1}{6}\left(3-\tilde{\zeta}\right)A+\frac{5}{4}A^3-\frac{5}{8}A^5
\end{equation}
with nontrivial equilibria at
\begin{equation}
    A^2=1\pm\sqrt{\frac{27-4\tilde{\zeta}}{15}},
\end{equation}
i.e., for $\tilde{\zeta}\in(0,27/4]$ with two co-existing solutions on the interval  $\tilde{\zeta}\in[3,27/4)$. The Hopf bifurcation at $\tilde{\zeta}=3$ is again subcritical with a geometric fold along the family of periodic orbits at $\tilde{\zeta}=27/4$ where $A=1$. 

Following the principles in the previous section, we may again obtain a hysteretic system response, alternating between the trivial equilibrium and the family of stable self-excited oscillations, by assuming a triggering perturbation of short duration and allowing autonomous variations in the scaled damping coefficients according to the differential equations
\begin{align}
\label{eq_4node_damplaw_large}
    \tilde{\zeta}_{2,4}'=\tau^{-1}\left(-\tilde{\zeta}_{2,4}+\frac{1}{2}\delta+\frac{3}{2}+\frac{15}{8}A^2\right),
\end{align}
since these imply that
\begin{equation}
    \tilde{\zeta}'=\tau^{-1}\left(-\tilde{\zeta}+\delta+3+\frac{15}{4}A^2\right).
\end{equation}
As in the previous section, we deem this scheme unrealizable in practice, since $A$ does not correspond to a measurable state of the system. In contrast to the analysis in the case of small linear damping, here, $A$ is not even accessible as an approximation of the amplitude of oscillation of either of nodes 2 or 4. Instead, $A$ may be understood as an approximation of the amplitude of oscillation of node 1, to which nodes 2 and 4 are connected through the network topology.

Similar observations apply to different configurations of the network with $Q=2$, $3$, or $4$, as shown in Appendix~\ref{appendix_4node_largedamping}. They prompt the analysis in the next section, which seeks to replace the explicit dependence on $A$ in the dynamics of $\zeta_i$ with some quantity $A_i$ that is computable from local measurements within the network.

\subsection{Local interactions}

We return to Eq.~\eqref{eq_4node_damplaw1} with the intent of replacing the occurrence of $A^2$ in the last term with a function of $u$ and its time history. In particular, we seek a dependence only on the components of $u$ visible to node 4 per the network topology. Finally, we seek to propose a formulation that generalizes to all the possible placements of the nonlinear oscillator in the limits of $\mathcal{O}(\epsilon)$ and $\mathcal{O}(1/\epsilon)$ damping for the linear oscillators.

By Eq.~\eqref{eq_4nodemultscale}, it follows that $A^2$ may be obtained from the $\mathcal{O}(1)$ squared amplitude of the projection of $\sqrt{2}u_4(t)$ onto $e^{\mathrm{j}2t}$. Notably, the analogous projections of $u_2(t)$ and $u_3(t)$ both vanish. Similarly, from Eq.~\eqref{eq_4nodemultscale_2} in Appendix~\ref{appendix_4node_smalldamping}, the term $A^2$ in each of the components of Eq.~\eqref{eq:damp_node2} may be obtained from the $\mathcal{O}(1)$ squared amplitude of the projections of $2\sqrt{3}u_1(t)$, $2\sqrt{3}u_3(t)$, and $2\sqrt{3}u_4(t)$, respectively, onto $e^{\mathrm{j}\sqrt{5}t}$ when the nonlinear oscillator is at node 2. Finally, from Eq.~\eqref{eq_4nodemultscale_3} in Appendix~\ref{appendix_4node_smalldamping}, when the nonlinear oscillator is at node 3, the term $A^2$ in each of the components of Eq.~\eqref{eq:damp_node3} may be obtained from the $\mathcal{O}(1)$ squared amplitude of the projections of $\sqrt{6}u_1(t)$ and $\sqrt{6}u_4(t)$, respectively, onto $e^{\mathrm{j}\sqrt{2}t}$. In this case, the analogous projection of $u_2(t)$ vanishes. 
While promising, however, these observations do not generalize to the case of high linear damping, since only the nonlinear oscillator has an $\mathcal{O}(1)$ amplitude proportional to $A$ in that limit.

As an alternative to extracting $A^2$ in each case from the displacement history of the corresponding linear node, consider instead the $\mathcal{O}(1)$ projections on the appropriate oscillatory mode of the \emph{net interaction forces} $u_2(t)+u_4(t)-2u_1(t)$, $u_1(t)+u_3(t)+u_4(t)-3u_2(t)$, $u_2(t)-u_3(t)$, and $u_1(t)+u_2(t)-2u_4(t)$ experienced by nodes 1 through 4. In the case that the nonlinear oscillator is located at node 1, Eq.~\eqref{eq_4nodemultscale} implies that these projections for nodes 2, 3, and 4 have amplitude $0$, $0$, and $3A/\sqrt{2}$ respectively. Similarly, in the case that the nonlinear oscillator is located at node 2, these projections for nodes 1, 3, and 4 have amplitude $2A/\sqrt{3}$, $2A/\sqrt{3}$, and $2A/\sqrt{3}$, respectively. Finally, in the case that the nonlinear oscillator is located at node 3, these projections for nodes 1, 2, and 4 have amplitude $A/\sqrt{6}$, $0$, and $A/\sqrt{6}$, respectively. This approach generalizes immediately to the case of high linear damping, since then the $\mathcal{O}(1)$ net interaction force equals the displacement of the nonlinear oscillator for nodes in its neighborhood in the network topology and $0$ for nodes outside of this neighborhood.

We are now in a position to design our four-node network with identical linear building blocks with damping dynamics given by
\begin{equation}
    \tau\zeta_i'=-\zeta_i+\frac{\delta+\zeta_{\text{HB}}+(\zeta_{\text{SN}}-\zeta_{\text{HB}})A_i^2/A_{\text{SN}}^2}{M}
\end{equation}
in the case of weak linear damping and 
\begin{equation}
    \tau\tilde{\zeta}_i'=-\tilde{\zeta}_i+\frac{\delta+\tilde{\zeta}_{\text{HB}}+(\tilde{\zeta}_{\text{SN}}-\tilde{\zeta}_{\text{HB}})A_i^2/A_{\text{SN}}^2}{M}
\end{equation}
in the case of strong linear damping with $A_i$ obtained from the $\mathcal{O}(1)$ projection of the net interaction force experienced by the $i$-th node onto the appropriate harmonic exponential $e^{\mathrm{j}\omega t}$. In these expressions, the $_{\text{HB}}$ and $_{\text{SN}}$ subscripts denote values at the Hopf and saddle-node bifurcations, while $M$ is some integer that ensures that the appropriate linear combination of damping coefficients generates a desired hysteretic response.

It remains to operationalize this scheme so as to eliminate explicit reference to $\mathcal{O}(1)$ terms and, consequently, the scaling with respect to $\epsilon$. Before doing so, we turn next to a generalization of the theory in this and the preceding subsections to an entire class of network topologies.

%%%%%%%%%%%%%%%%%%%%%%%%%%%%%%%%%%%%%

\section{General active network filters} \label{sec_networkfilters}

We proceed to generalize the observations in Sect.~\ref{sec:a model problem} for the particular four-node network to arbitrary network topologies. We remain concerned with the case of a nonlinear oscillator occupying a single node in an otherwise linear oscillator network.

To this end, consider an undirected network topology of size $N$, represented by the Laplacian $L$, with the nonlinear oscillator located at node $Q$. Let the equations of motion be given in the general form
\begin{equation}
\label{eq:gennetworkdyneq}
    \ddot{u}+\epsilon C(u)\dot{u}+Ku=F,\,u\in\mathbb{R}^N,
\end{equation}
in terms of the diagonal damping matrix $\epsilon C(u)$ with \begin{equation}
    C_{Q,Q}(u)=-\nu+\eta u_Q^2-\eta u_Q^4    
\end{equation}
and $C_{i,i}(u)=\zeta_i$ for $i\ne Q$, stiffness matrix $K=I_N+L$, and exogenous excitation vector $F(t)$. Since the stiffness matrix is positive definite, there exists an orthogonal matrix $P$ such that the modal stiffness matrix $\tilde{K}:=P^\text{T} K P=\mathrm{diag}(\omega_1^2,\ldots,\omega_N^2)$ in terms of the natural frequencies $\omega_1=1\le\ldots\le\omega_N$. We assume below that these frequencies are all different. 

\subsection{Small linear damping}

Let $u=Px$ in terms of the modal coordinate vector $x$. Equation~\eqref{eq:gennetworkdyneq} then becomes
\begin{equation}
    \ddot{x}+\epsilon\tilde{C}(x)\dot{x}+\tilde{K}x=\tilde{F},
\end{equation}
where $\tilde{C}(x):=P^\text{T} C(x) P$ is the modal damping matrix and $\tilde{F}(t):=P^\text{T} F(t)$. The stability of the trivial equilibrium at $x=0$ in the absence of excitation is now determined by the roots of the determinant 
\begin{equation}
    \left|\lambda^2 I_N+\epsilon \lambda\tilde{C}(0)+\tilde{K}\right|    
\end{equation}
Indeed, provided that $\epsilon\ll 1$ and $\zeta_i,\nu,\eta\sim\mathcal{O}(1)$, these roots are of the form
\begin{equation}
    \lambda_i=\pm\mathrm{j}\omega_i-\frac{\epsilon}{2}\tilde{C}_{i,i}(0)+\mathcal{O}(\epsilon^2),\,k=i,\ldots,N.
\end{equation}
Asymptotic stability thus follows as long as all diagonal elements of $\tilde{C}(0)$ are positive, while instability results as soon as at least one such element is negative. Generically, one-parameter variations crossing the threshold of instability are associated with Hopf bifurcations from which emanate one-parameter families of periodic orbits corresponding to self-excited sustained oscillations. By definition,
\begin{equation}
    \tilde{C}_{i,i}(0)=-P^2_{Q,i}\nu+\sum_{k\ne Q}P^2_{k,i}\zeta_k.
\end{equation}
In terms of the notation $I:=\mathrm{argmax}_iP^2_{Q,i}$ and $\zeta:=\sum_{k\ne Q}P^2_{k,I}\zeta_k$, it follows that the trivial equilibrium is asymptotically stable when $\zeta>P^2_{Q,I}\nu$ and unstable when $\zeta<P^2_{Q,I}\nu$. By the earlier assumption, $\omega_I\ne 1$.

We proceed to again use the method of multiple scales to determine the fate of the periodic orbits born at the Hopf bifurcation under variations in $\zeta$. To this end, substitution of
\begin{equation}\label{eq_multscale}
    u(t)=P_{\cdot,I}A(\epsilon t)\cos\left(\omega_I t+\phi(\epsilon t)\right)+\epsilon v(t)
\end{equation}
($P_{\cdot,I}$ denotes the $I$-th column of matrix $P$) into the fully nonlinear governing equations yields the coefficients
\begin{equation}
    -2\omega_I A\phi'
\end{equation}
and
\begin{equation}
\label{eq:Adyn}
    -2\omega_I A'+\omega_I\left(P^2_{Q,I}\nu-\zeta\right)A+\frac{P^4_{Q,I}\eta\omega_I}{4}A^3-\frac{P^6_{Q,I}\eta\omega_I}{8}A^5
\end{equation}
in front of the secular terms $\cos(\omega_I t+\phi)$ and $\sin(\omega_I t+\phi)$, respectively, in the dynamics of $\sum_k P_{k,I}v_k$ (with natural frequency $\omega_I$). These expressions vanish at nontrivial equilibrium values (in the slow time scale) of $A$ and $\phi$ provided that
\begin{equation}
\label{eq:Anull}
    A^2=\frac{1}{P^2_{Q,I}}\pm\sqrt{\frac{P^2_{Q,I}(8\nu+\eta)-8\zeta}{P^6_{Q,I}\eta}}
\end{equation}
for $\zeta\in[0,P^2_{Q,I}(8\nu+\eta)/8]$ with two co-existing solutions on the interval $\zeta\in[P^2_{Q,I}\nu,P^2_{Q,I}(8\nu+\eta)/8)$. We conclude that the Hopf bifurcation at $\zeta=P^2_{Q,I}\nu$ is subcritical and that the branch of periodic orbits has a geometric fold at a saddle-node bifurcation at $\zeta=P^2_{Q,I}(8\nu+\eta)/8$, where $A=1/|P_{Q,I}|$.

For any $k\ne Q$, the net interaction force experienced by the $k$-th node is given to $\mathcal{O}(1)$ by the matrix product
\begin{align}
    -L_{k,\cdot}u(t)&=-L_{k,\cdot}P_{\cdot,I}A(\epsilon t)\cos\left(\omega_I t+\phi(\epsilon t)\right)\nonumber\\
    &=-P_{k,I}(\omega_I^2-1)A(\epsilon t)\cos\left(\omega_I t+\phi(\epsilon t)\right)
\end{align}
($L_{k,\cdot}$ denotes the $k$-th row of Laplacian $L$) which vanishes at all displacement nodes of the $I$-th mode shape (i.e., nodes $k$ with $P_{k,I}=0$). These are also the terms that do not contribute to $\zeta$. The equations
\begin{equation}\label{eq_smalldampingdynamics_nodek}
    \zeta_k'=\tau^{-1}\left(-\zeta_k+\frac{\delta+P^2_{Q,I}\nu+P^4_{Q,I}\eta A_k^2/8}{1-P^2_{Q,I}}\right)
\end{equation}
with $A_k$ equal to the amplitude of the $\mathcal{O}(1)$ projection onto $e^{\mathrm{j}\omega_I t}$ of $L_{k,\cdot}u/P_{k,I}(\omega_I^2-1)$ for $P_{k,I}\ne 0$ and $0$ otherwise
now imply that
\begin{equation}\label{eq_smalldampingdynamics}
    \zeta'=\tau^{-1}\left(-\zeta+\delta+P^2_{Q,I}\nu+P^4_{Q,I}\eta A^2/8\right),
\end{equation}
which produces the desired hysteretic response following a triggering excitation of sufficient magnitude and duration.

For arbitrary $\tau$ and $\delta$, the coupled dynamics has a trivial equilibrium at $(\zeta,A)=(P_{Q,I}^2\nu+\delta,0)$. Additional equilibria are found along the middle branch of the $A$ nullcline, where
\begin{equation}
    \zeta= P_{Q,I}^2\nu + \delta+\frac{P_{Q,I}^2\eta \pm P_{Q,I}\sqrt{P_{Q,I}^2\eta^2-32\eta\delta} }{16}
\end{equation}
provided that $\delta\le P_{Q,I}^2\eta/32$. These equilibria emerge from the end points of the middle branch of the $A$ nullcline as $\delta$ increases from $0$ and disappear at a saddle-node bifurcation when $\delta=P_{Q,I}^2\eta/32$, $\zeta=P_{Q,I}^2(32\nu+3\eta)/32$, and $A=1/\sqrt{2}|P_{Q,I}|$. As in the example in Sect.~\ref{sec:a model problem}, the trivial equilibrium is asymptotically stable for any $\tau$ and $\delta$, the nontrivial equilibrium with $A<1/\sqrt{2}|P_{Q,I}|$ is always a saddle, the nontrivial equilibrium with $A>1/\sqrt{2}|P_{Q,I}|$ is stable for $\tau<1/2\delta$ and unstable for $\tau>1/2\delta$ (the extra factor of $2$ is a result of a different definition of $\zeta$ in Sect.~\ref{sec:a model problem}), and the bifurcation at $\tau=1/2\delta$ is always a Hopf bifurcation.

\subsection{Large linear damping}

We repeat this analysis in the case of large linear damping such that $\zeta_i\sim\mathcal{O}(1/\epsilon^2)$ while $\nu,\eta\sim\mathcal{O}(1)$. In this case, the trivial equilibrium loses stability at a Hopf bifurcation from which emanates a one-parameter family of periodic responses with limiting natural frequency $\sqrt{K_{Q,Q}}$. Substitution of the multiple-scale ansatz
\begin{equation}
    u(t)=A(\epsilon t)\cos\left(\sqrt{K_{Q,Q}}t+\phi(\epsilon t)\right)\mathbf{e}_Q+\epsilon v(t)
\end{equation}
with $\mathbf{e}_{Q,i}=\delta_{Qi}$ then yields the amplitude equation
\begin{equation}
    A'=\frac{1}{2}\left(\nu-\tilde{\zeta}\right)A+\frac{\eta}{8}A^3-\frac{\eta}{16}A^5,
    \label{eq:ampeqhighdamp1}
\end{equation}
where
\begin{equation}
    \tilde{\zeta}:=\sum_{k\ne Q}\frac{K_{Q,k}K_{k,Q}}{K_{Q,Q}}\tilde{\zeta}_k, \quad \tilde{\zeta}_k = \frac{1}{\epsilon^2\zeta_k}.
\end{equation}
Nontrivial equilibria are located at
\begin{equation}
    A^2=1\pm\sqrt{\frac{\eta+8(\nu-\tilde{\zeta})}{\eta}}
    \label{eq:ampeqhighdamp2}
\end{equation}
i.e., for $\tilde{\zeta}\in(0,(8\nu+\eta)/8]$ with two co-existing solutions on the interval $\tilde{\zeta}\in[\nu,(8\nu+\eta)/8)$. We conclude that the Hopf bifurcation occurs when $\tilde{\zeta}=\nu$ and that the saddle-node bifurcation occurs when $\tilde{\zeta}=\nu+\eta/8$.

For $k\ne Q$,
\begin{equation}
    K_{Q,k}K_{k,Q}=L_{Q,k}L_{k,Q},
\end{equation}
which is nonzero only if the $k$-th node is in the neighborhood of the $Q$-th node. It follows that $\tilde{\zeta}$ only depends on the damping coefficients associated with nodes in the neighborhood of $Q$ and that we can extract $A^2$ for each such node from the $\mathcal{O}(1)$ net interaction force projected onto $e^{\mathrm{j}\sqrt{K_{Q,Q}}t}$, whereas this projection vanishes for all other nodes. Since for our network design
\begin{equation}
    \sum_{k\ne Q}K_{Q,k}K_{k,Q}=K_{Q,Q}-1,
\end{equation}
the equations
\begin{equation}
    \tilde{\zeta}_k'=\tau^{-1}\left(-\tilde{\zeta}_k+\frac{K_{Q,Q}(\delta+\nu+\eta A_k^2/8)}{K_{Q,Q}-1}\right)
\end{equation}
with $A_k$ equal to the amplitude of the $\mathcal{O}(1)$
projection onto $e^{\mathrm{j}\sqrt{K_{Q,Q}}t}$ of $u_Q$ for $L_{k,Q}\ne 0$ and $0$ otherwise now imply that
\begin{equation}
    \tilde{\zeta}'=\tau^{-1}\left(-\tilde{\zeta}+\delta+\nu+\eta A^2/8\right),
    \label{eq:ampeqhighdamp3}
\end{equation}
which again produces the desired hysteretic response following a triggering excitation of sufficient magnitude and duration.

For arbitrary $\tau$ and $\delta$, the coupled dynamics has a trivial equilibrium at $(\tilde{\zeta},A)=(\nu+\delta,0)$. Additional equilibria are found along the middle branch of the $A$ nullcline, where now
\begin{equation}
    \tilde{\zeta}=\nu + \delta+\frac{\eta \pm \sqrt{\eta^2-32\eta\delta} }{16}
\end{equation}
provided that $\delta\le \eta/32$. These equilibria emerge from the end points of the middle branch of the $A$ nullcline as $\delta$ increases from $0$ and disappear at a saddle-node bifurcation when $\delta=\eta/32$, $\tilde{\zeta}=\nu+3\eta/32$, and $A=1/\sqrt{2}$. As before, the trivial equilibrium is asymptotically stable for any $\tau$ and $\delta$, the nontrivial equilibrium with $A<1/\sqrt{2}$ is always a saddle, the nontrivial equilibrium with $A>1/\sqrt{2}$ is stable for $\tau<1/2\delta$ and unstable for $\tau>1/2\delta$, and the bifurcation at $\tau=1/2\delta$ is always a Hopf bifurcation.

%%%%%%%%%%
% Examples
%%%%%%%%%%

\section{Implementation and numerical results}\label{sec:examples}

In this section, we operationalize the recipe in Sect.~\ref{sec_networkfilters} without restricting attention to the $\epsilon\ll 1$ limit. Only by doing so, do we achieve the physical realizability sought at the outset. In this case, we omit any mention of $\mathcal{O}(1)$ in seeking to estimate $A_k$ from the net interaction force acting on the $k$-th oscillator. Instead, in the case of initially small linear damping, we let
\begin{align}\label{eq_projectionsmalldamping}
   A_k(t) &= \left| \int_{t-2\pi/\omega_I}^t \frac{\omega_I L_{k,\cdot}u(s)e^{\mathrm{j}\omega_I s}}{\pi P_{k,I}(\omega_I^2-1)}\, \mathrm{d} s \right|
\end{align}
for $P_{k,I}\ne 0$, and $A_k=0$ otherwise. Similarly, in the case of initially large linear damping, we let
\begin{align}\label{eq_projectionlargedamping}
   &A_k(t) = \left| \int_{t-2\pi/\sqrt{K_{Q,Q}}}^t \frac{\sqrt{K_{Q,Q}} u_Q(s)e^{\mathrm{j}\sqrt{K_{Q,Q}}s}}{\pi}\, \mathrm{d} s\right|
\end{align}
for $L_{k,Q}\ne 0$, and $A_k=0$ otherwise. When either of these are coupled to the damping and network dynamics, the variable substitution $s\mapsto t-\sigma$ in either integral yields a system of integro-differential equations with distributed delay of the form
\begin{equation}
    \dot{y}(t)=f\left(y(t),\int_0^\upsilon \kappa(t,\sigma)y(t-\sigma)\,\mathrm{d}\sigma\right)
\end{equation}
for some kernel $\kappa(t,\sigma)$ and $\upsilon=2\pi/\omega_I$ or $2\pi/\sqrt{K_{Q,Q}}$, respectively~\cite{khasawneh2011stability}. In simulations, we approximate the integral term using numerical quadrature, thereby obtaining a system of delay differential equations with finitely many discrete delays~\cite{gallage2017approximation}. Here, we choose to discretize either of the integrals using the trapezoidal rule with four evenly spaced points on the interval $[0,\upsilon]$, since this suffices to compute the amplitude of a harmonic signal of the frequency expected in the limit as $\epsilon\rightarrow 0$. Although the actual response frequency varies with $\epsilon$, we proceed to use this simplified definition of $A_k$ to investigate the possibility of hysteretic variations in response amplitude and damping coefficients also away from the $\epsilon\rightarrow0$ limit. 

\begin{figure}[htbp]
    \centering
    \includegraphics[width=.45\textwidth]{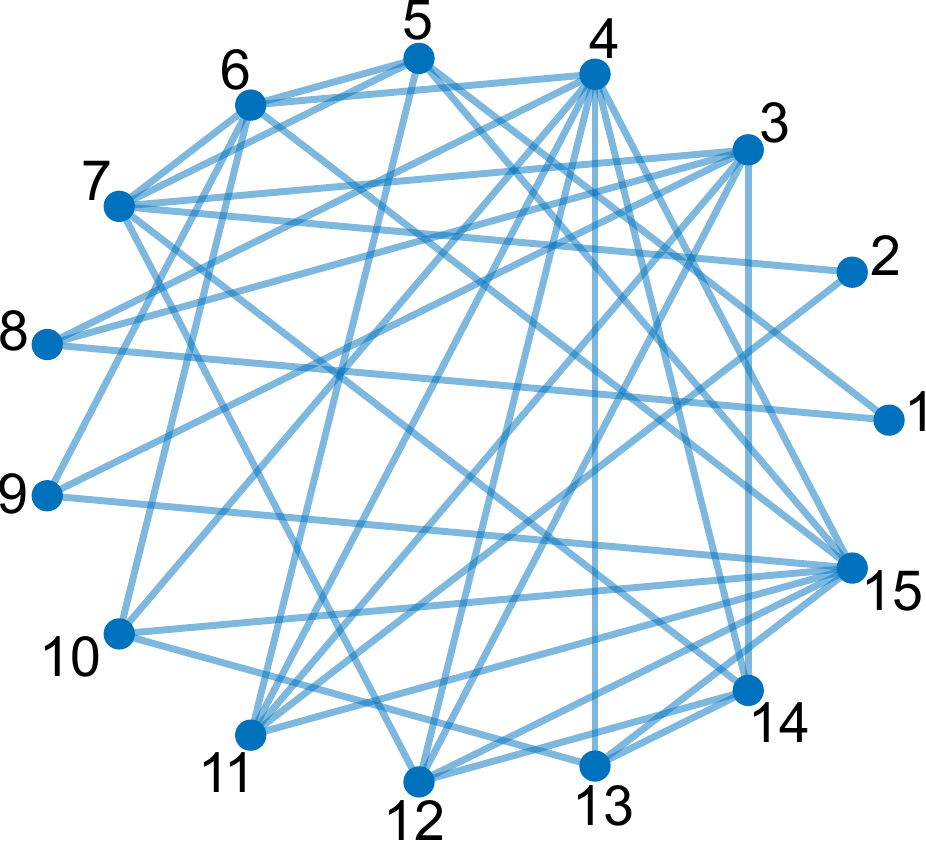}
    \caption{A nondegenerate undirected network topology with 15 nodes corresponding to an adjacency matrix with nonzero entries $A_{1,5}$, $A_{1,8}$, $A_{2,7}$, $A_{2,11}$, $A_{3,7}$, $A_{3,8}$, $A_{3,9}$, $A_{3,11}$, $A_{3,12}$, $A_{3,14}$, $A_{4,6}$, $A_{4,8}$, $A_{4,10}$, $A_{4,11}$, $A_{4,12}$, $A_{4,13}$, $A_{4,14}$, $A_{4,15}$, $A_{5,6}$, $A_{5,7}$, $A_{5,11}$, $A_{5,15}$, $A_{6,7}$, $A_{6,9}$, $A_{6,10}$, $A_{6,15}$, $A_{7,12}$, $A_{7,14}$, $A_{9,15}$, $A_{10,13}$, $A_{10,15}$, $A_{11,15}$, $A_{12,14}$, $A_{12,15}$, $A_{13,14}$, $A_{13,15}$, and so on by symmetry.}
    \label{15node_network_example}
\end{figure}

\begin{figure*}[htbp]
\centering
\subfloat[]{\includegraphics[width=.45\textwidth]{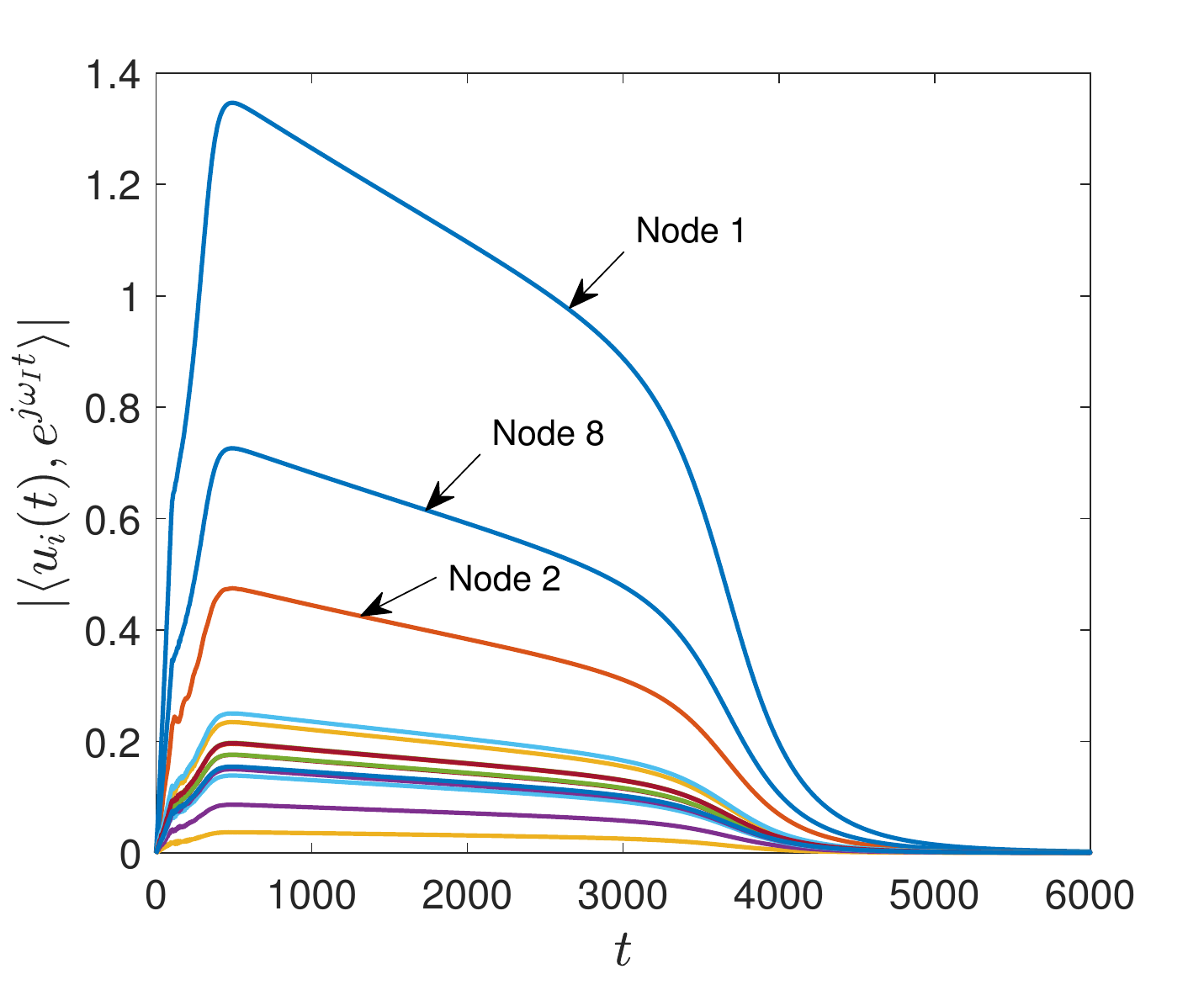}}
\subfloat[]{\includegraphics[width=.45\textwidth]{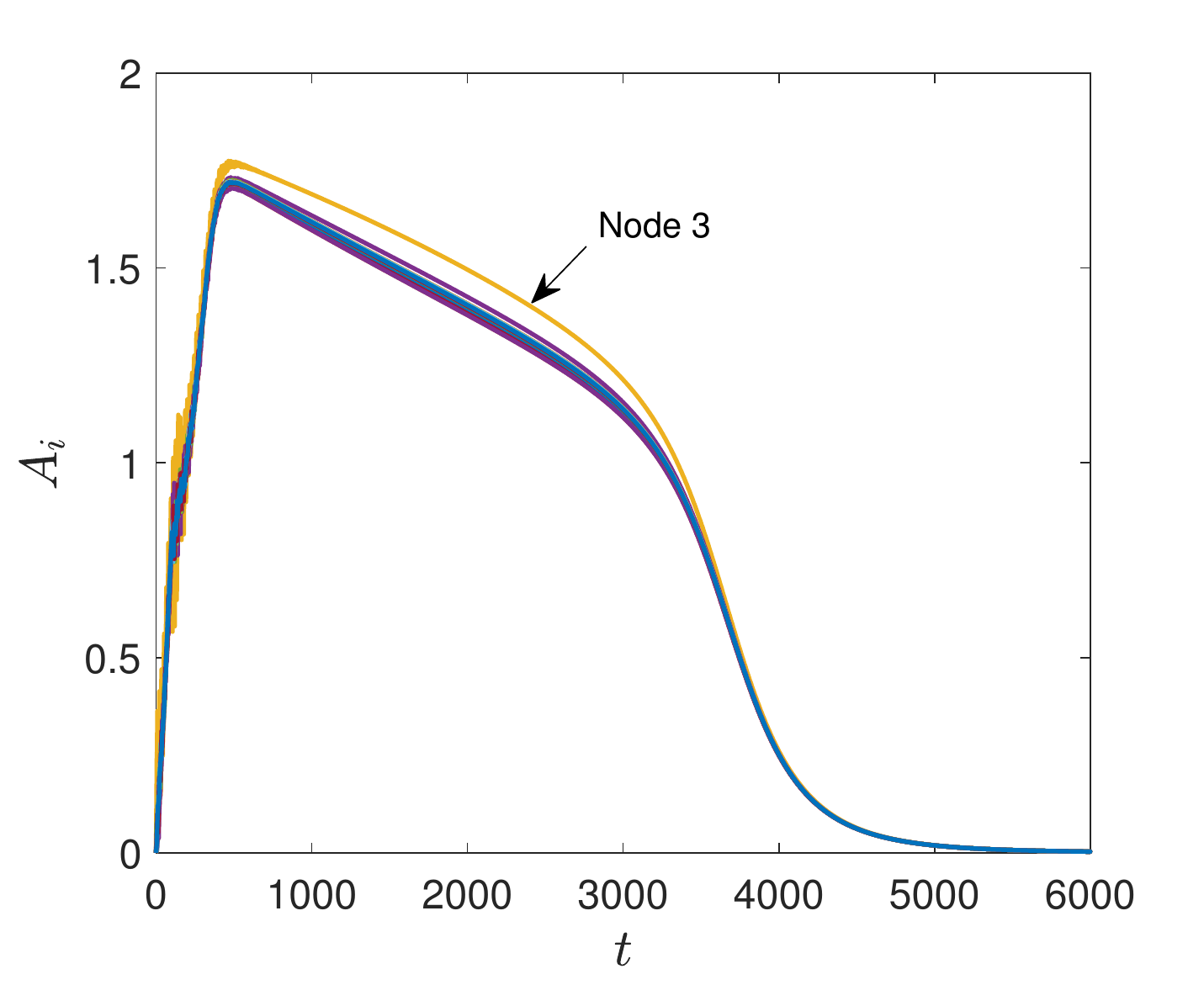}}\\
\subfloat[]{\includegraphics[width=.45\textwidth]{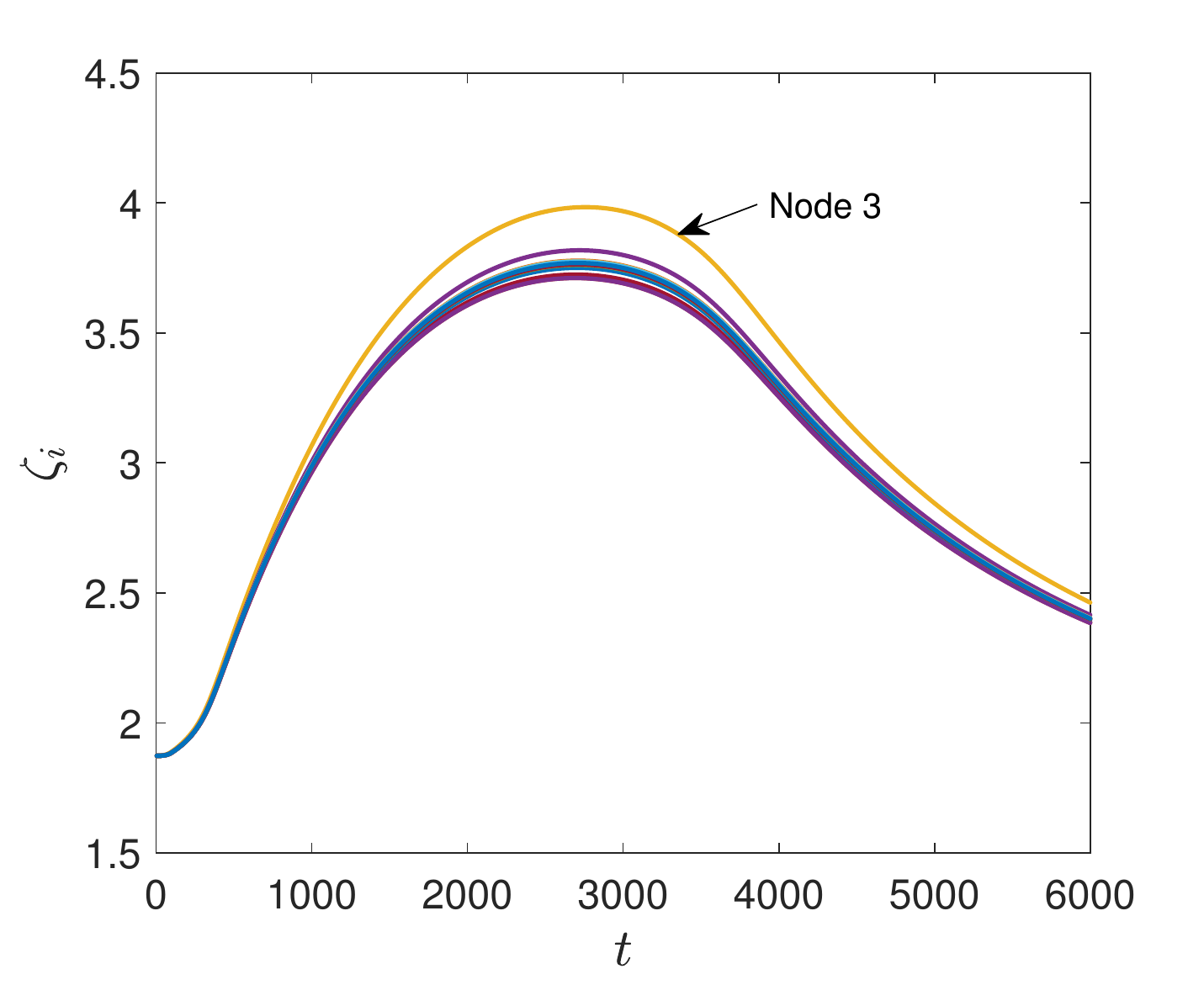}}
\subfloat[]{\includegraphics[width=.45\textwidth]{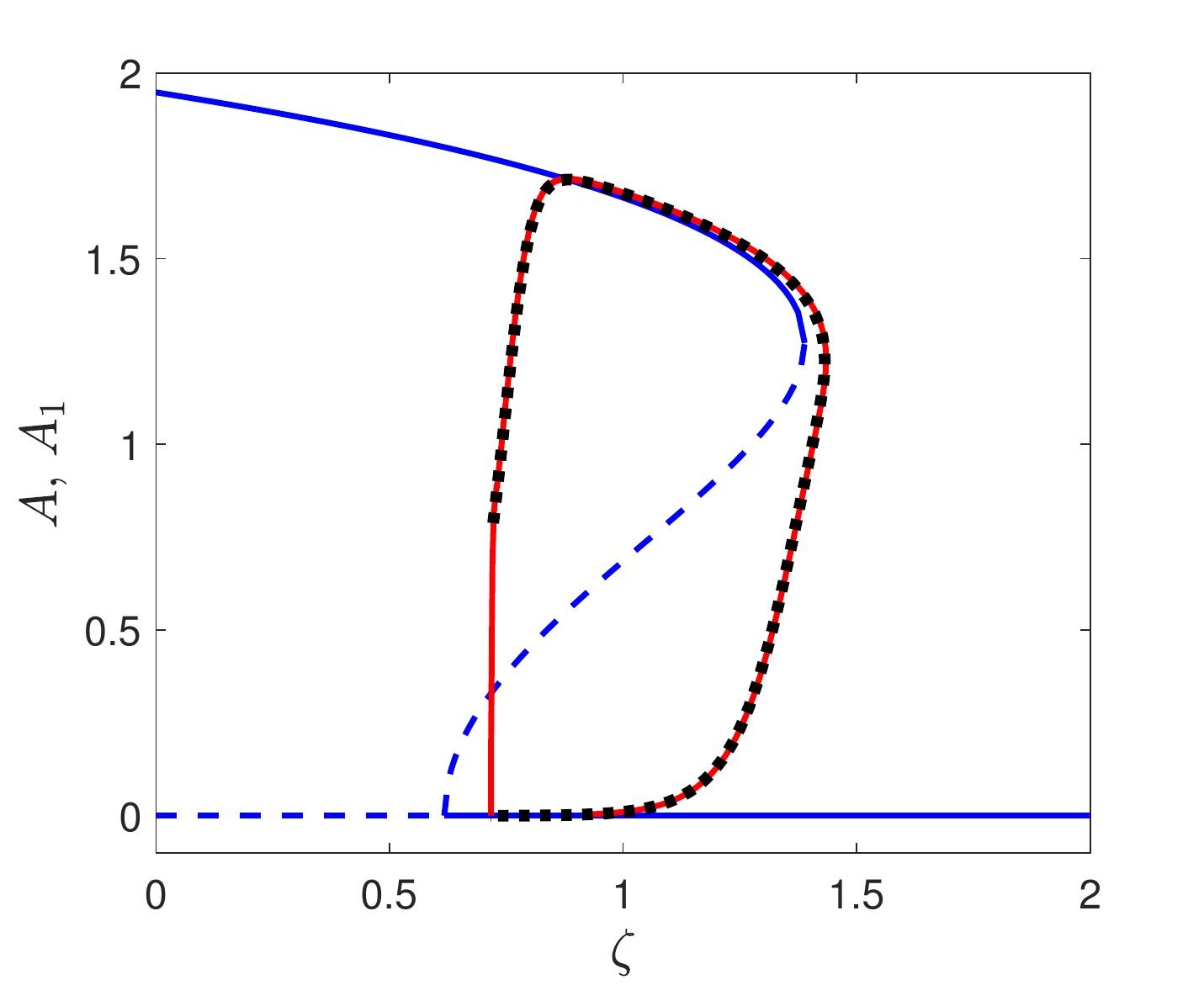}}
\caption{Time histories of (a) nodal amplitudes $|\langle u_i(t),e^{j\omega_I t}\rangle|$, (b) nodal estimates $A_i$, and (c) nodal damping coefficients $\zeta_i$ for the 15-node network with topology shown in Fig.~\ref{15node_network_example} with $Q=1$, $I=2$, and $\omega_I\approx 1.5212$. Panel (d) shows the corresponding hysteretic trajectory (red solid) projected onto $\zeta=\sum_{k\ne Q}P_{k,I}^2\zeta_k$ and $A_1$, the trajectory (black dotted) of the coupled $A$ and $\zeta$ dynamics obtained from Eqs.~\eqref{eq:Adyn} and \eqref{eq_smalldampingdynamics}, and the corresponding $A$ nullcline obtained from Eq.~\eqref{eq:Anull}. The full simulation is initialized with zero initial displacements and velocities, $A_i=0$ for $t\in[-2\pi/\omega_I,0]$ and $\zeta_i=(\delta+P^2_{Q,I}\nu)/(1-P^2_{Q,I})$. The simulation of the coupled $A$ and $\zeta$ dynamics is initialized with the values of $A_1$ and $\zeta$ at the conclusion of the initial period of exogenous excitation. Here, $\epsilon=0.01$, $\nu=1$, $\eta=10$, $F = 3\epsilon \sin (\omega_I t)\mathbf{e}_1$ for $0\leq t\leq 1/\epsilon$, $\delta=0.1$, and $\tau=20$. }\label{15_node_time_example_Q1_ep0_01}
\end{figure*}

\begin{figure*}[htbp]
\centering
\subfloat[]{\includegraphics[width=.45\textwidth]{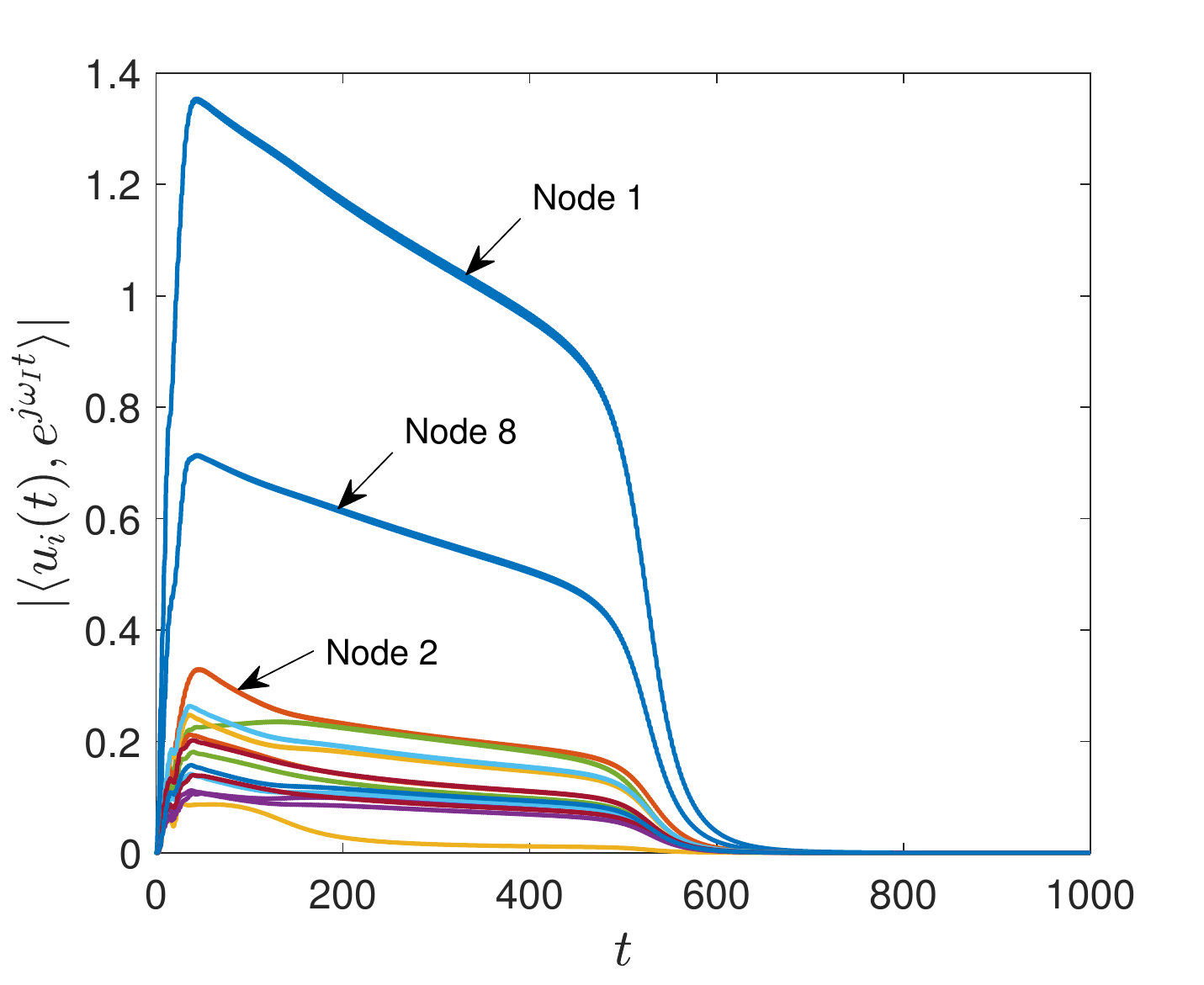}}
\subfloat[]{\includegraphics[width=.45\textwidth]{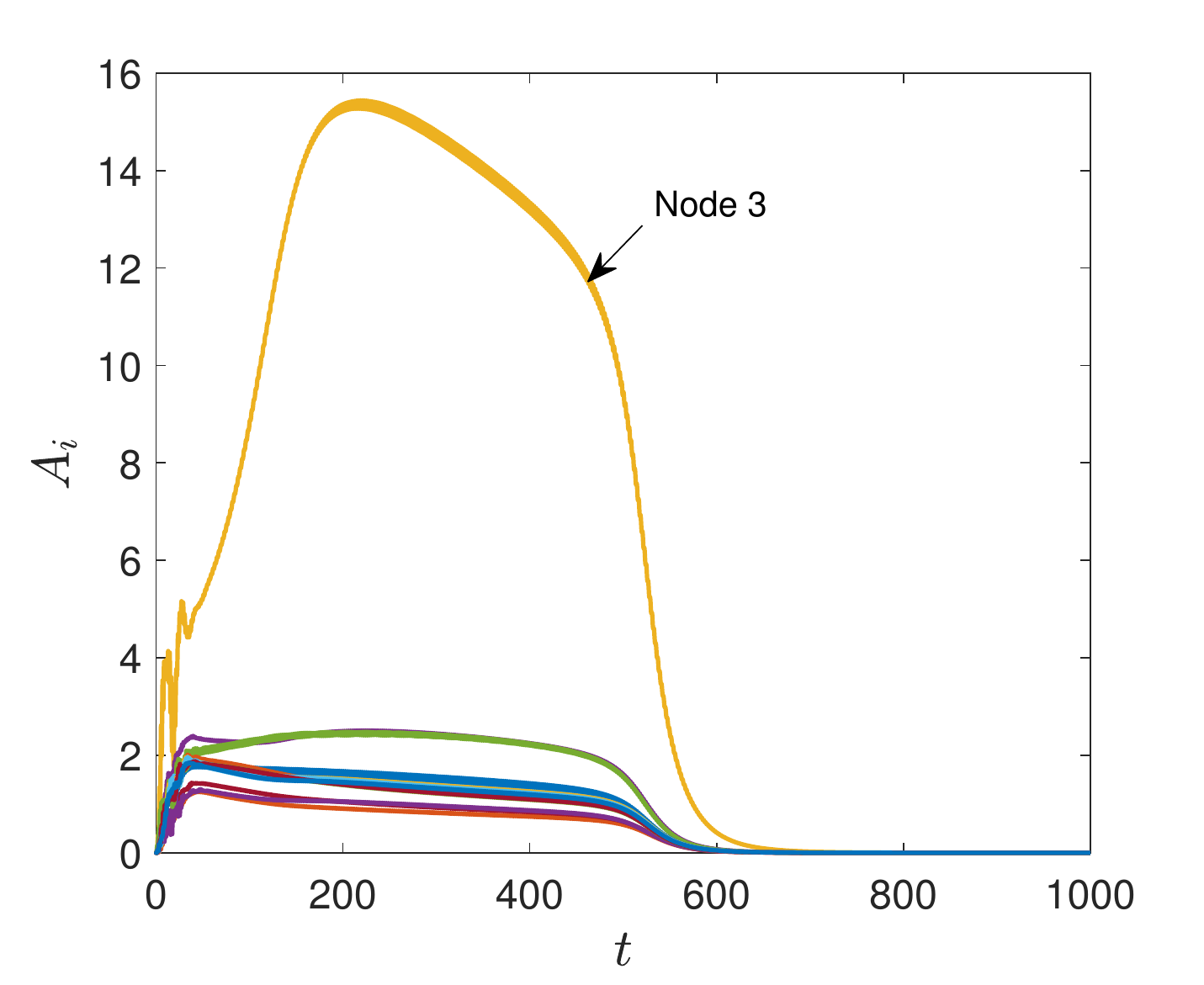}}\\
\subfloat[]{\includegraphics[width=.45\textwidth]{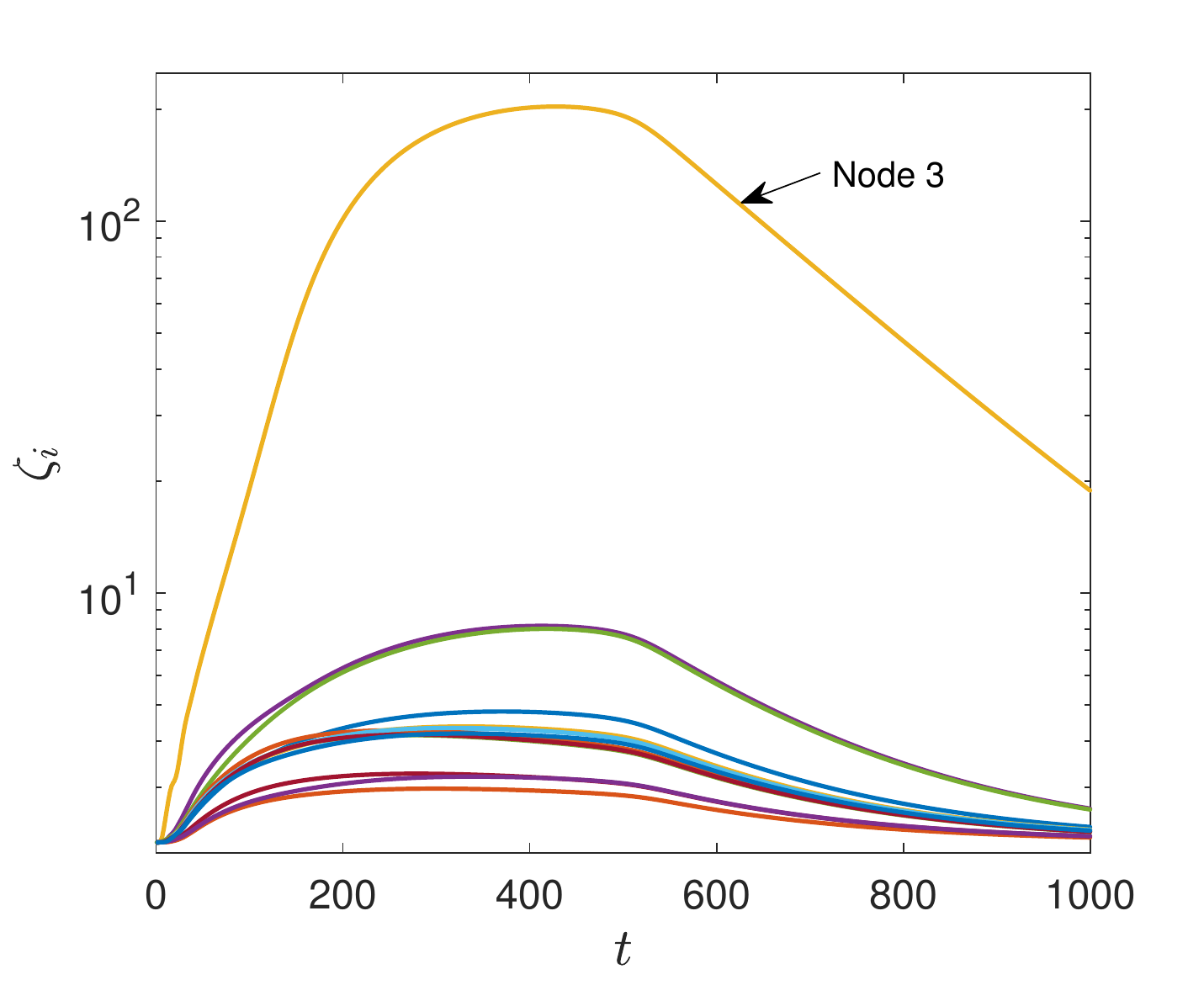}}
\subfloat[]{\includegraphics[width=.45\textwidth]{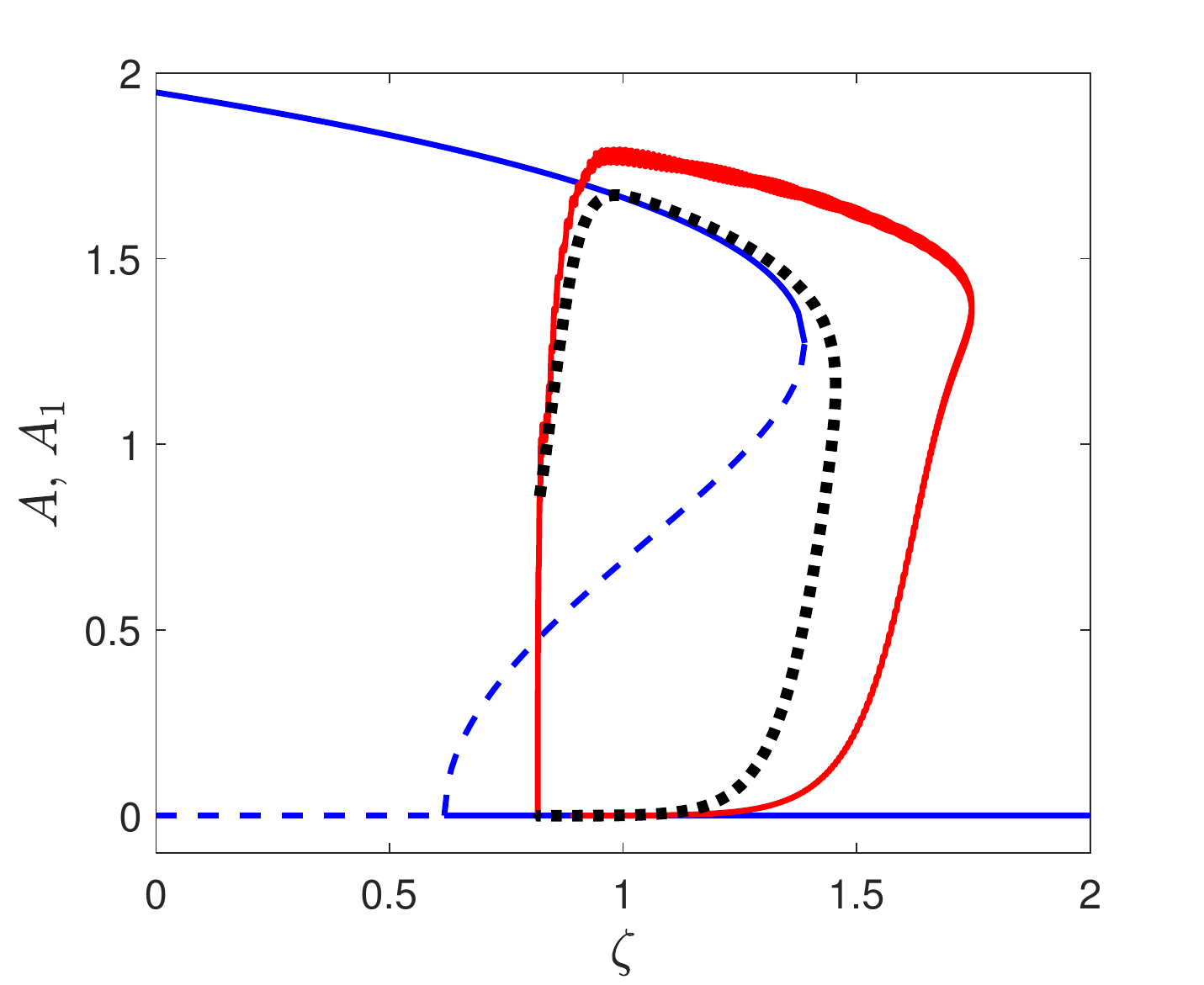}\label{15_node_time_example_Q1_ep0_1_d}}
\caption{Time histories of (a) nodal amplitudes $|\langle u_i(t),e^{j\omega_I t}\rangle|$, (b) nodal estimates $A_i$, and (c) nodal damping coefficients $\zeta_i$ for the 15-node network with topology shown in Fig.~\ref{15node_network_example} with $Q=1$, $I=2$, and $\omega_I\approx 1.5212$. Panel (d) shows the corresponding hysteretic trajectory (red solid) projected onto $\zeta=\sum_{k\ne Q}P_{k,I}^2\zeta_k$ and $A_1$, the trajectory (black dotted) of the coupled $A$ and $\zeta$ dynamics obtained from Eqs.~\eqref{eq:Adyn} and \eqref{eq_smalldampingdynamics}, and the corresponding $A$ nullcline obtained from Eq.~\eqref{eq:Anull}. The full simulation is initialized with zero initial displacements and velocities, $A_i=0$ for $t\in[-2\pi/\omega_I,0]$ and $\zeta_i=(\delta+P^2_{Q,I}\nu)/(1-P^2_{Q,I})$. The simulation of the coupled $A$ and $\zeta$ dynamics is initialized with the values of $A_1$ and $\zeta$ at the conclusion of the initial period of exogenous excitation. Here, $\epsilon=0.1$, $\nu=1$, $\eta=10$, $F = 3\epsilon \sin (\omega_I t)\mathbf{e}_1$ for $0\leq t\leq 1/\epsilon$, $\delta=0.2$, and $\tau=20$. }\label{15_node_time_example_Q1_ep0_1}
\end{figure*}

As an example, consider the 15-node network topology shown in  Fig.~\ref{15node_network_example}. With $Q=1$, the analysis in Sect.~\ref{sec_networkfilters} predicts that $I=2$ and, consequently, that instability in the small damping limit occurs at a Hopf bifurcation with limiting frequency $\omega_I\approx 1.52$ and linear mode shape $P_{\cdot,I} \approx (-0.7856, \allowbreak 0.2785, \allowbreak 0.0210, \allowbreak 0.0500, \allowbreak -0.1150, \allowbreak 0.0807, \allowbreak 0.1030, \allowbreak -0.4239, \allowbreak 0.1138, \allowbreak 0.1364, \allowbreak 0.0881, \allowbreak 0.1027,\quad \allowbreak 0.1456, \allowbreak 0.1146, \allowbreak 0.0901 )^\top$. From the absence of displacement nodes in the corresponding mode shape, our construction generates nontrivial dynamics for all linear damping coefficients. This is verified by the dynamics shown in Figs.~\ref{15_node_time_example_Q1_ep0_01} and \ref{15_node_time_example_Q1_ep0_1} for $\tau=20$, $\nu=1$, $\eta=10$, and $(\epsilon,\delta)=(0.01,0.1)$ and $(0.1,0.2)$, respectively.

In each of these simulations, the network is initialized with zero displacements and velocities, $A_i(t)=0$ for $t\in[-2\pi/\omega_I,0]$, and $\zeta_i=(\delta+P_{Q,I}^2\nu)/(1-P_{Q,I}^2)\approx 1.87$. A burst of exogenous harmonic excitation given by $3\epsilon\sin\omega_I t$ applied to the first node for $\epsilon t\in[0,1]$ provides the impetus for the subsequent dynamics. In Fig.~\ref{15_node_time_example_Q1_ep0_01}, the behavior is as expected from the perturbation analysis. The nodal amplitudes are computed from the projection of $u(t)$ onto $e^{\mathrm{j}\omega_I t}$, and are found to scale according to the mode shape $P_{\cdot,I}$. The nodal estimates $A_i$ and nodal damping coefficients $\zeta_i$, respectively, trace nearly identical time histories for different $i$ with some deviation for node 3 caused by the small value of $|P_{3,I}|$ which exaggerates deviations from the predicted behavior in the $\epsilon\rightarrow 0$ limit. When projected onto the $(\zeta,A_1)$ plane, the transient trajectory subsequent to the termination of the exogenous excitation follows closely that predicted from the perturbation analysis.

As $\epsilon$ increases, deviations from the predictions of the multiple-scale analysis are further amplified by small denominators in Eq.~\eqref{eq_projectionsmalldamping}. This is clearly demonstrated in Fig.~\ref{15_node_time_example_Q1_ep0_1} in the case, again, of significant deviations of $A_3$ and $\zeta_3$ from the remaining nodal estimates and damping coefficients. Since the contribution to $\zeta$ from $\zeta_3$ is scaled by $P_{3,I}^2$, this deviation is not reflected in Fig.~\ref{15_node_time_example_Q1_ep0_1_d} where hysteresis is again observed with $\epsilon=0.1$. In this case, the transient trajectory  only returns to $A_1=0$ after some delay. The apparent thickness of the curve along the upper branch corresponds to rapid oscillations in $A_1$ due to a mismatch between the response frequency and $\omega_I$.

\begin{figure*}[htbp]
\centering
\subfloat[]{\includegraphics[width=.5\textwidth]{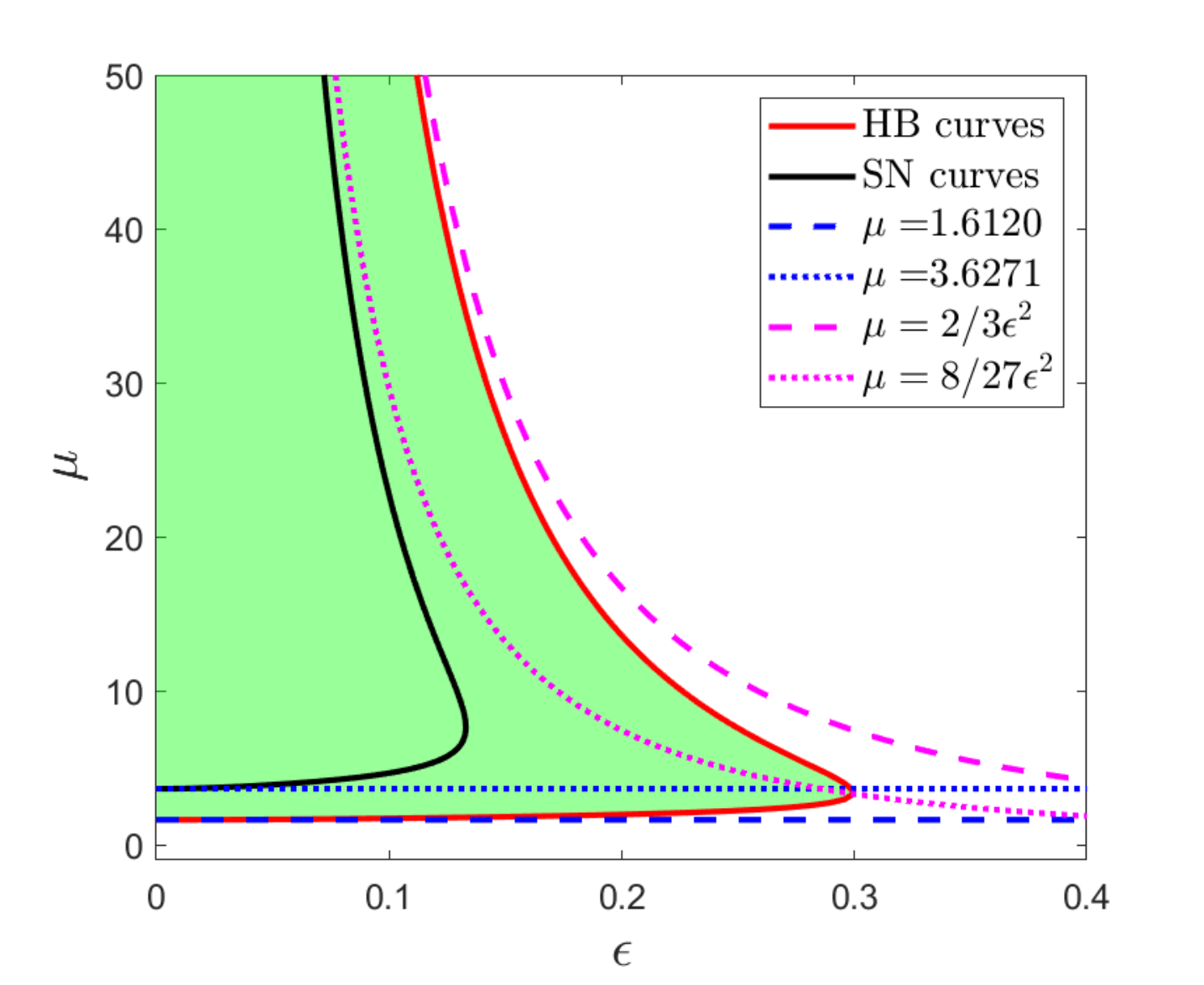}\label{15_node_HB_and_SN_curves_a}}
\subfloat[]{\includegraphics[width=.5\textwidth]{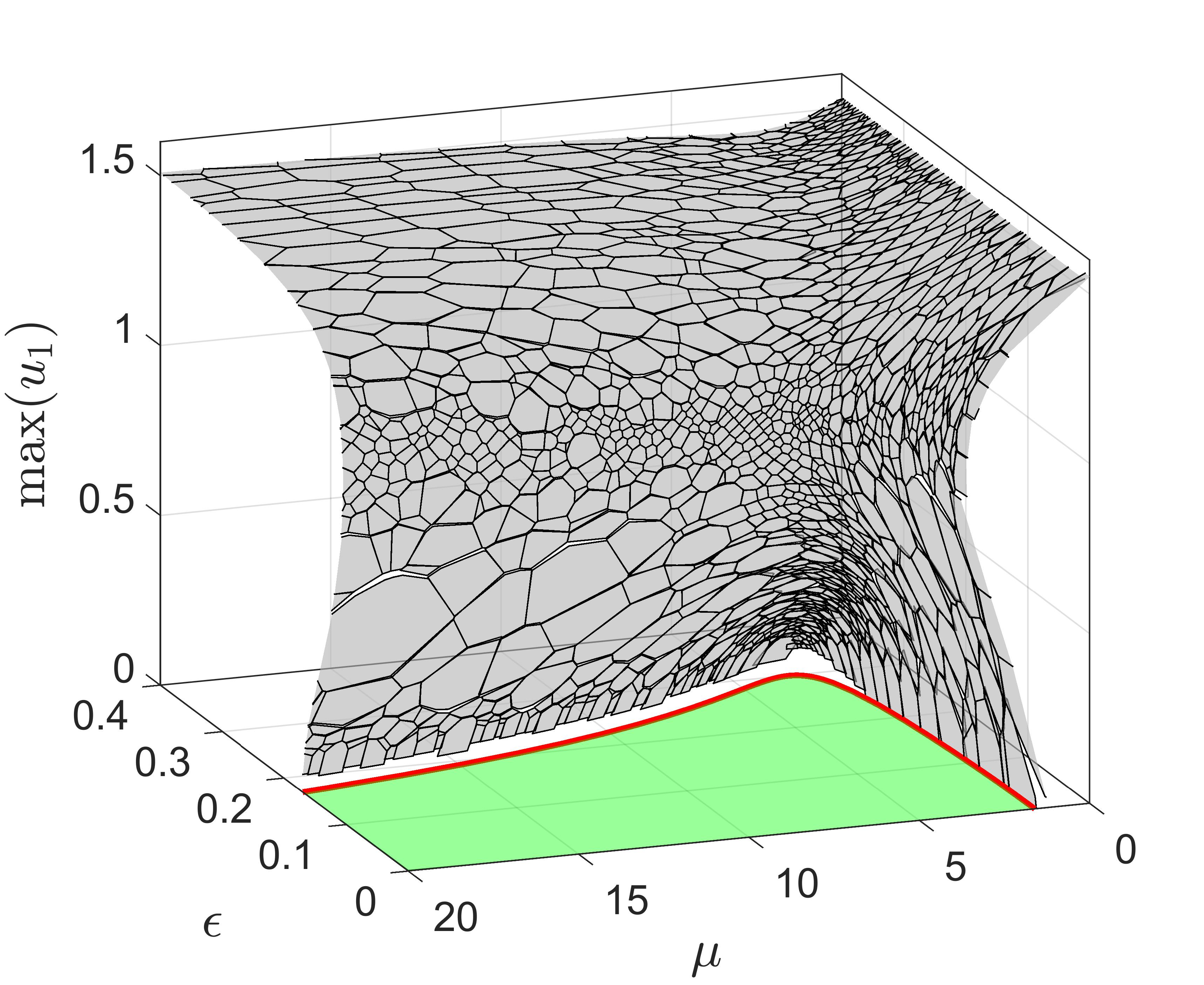}\label{15_node_HB_and_SN_curves_b}}
\caption{Bifurcation analysis for the 15-node network example shown in Fig.~\ref{15node_network_example} with $Q=1$ and $\zeta_i=\mu$ for $i\neq Q$ obtained using the equilibrium (\texttt{eq}) and periodic-orbit (\texttt{po}) toolboxes of the \textsc{coco} software platform~\cite{dankowicz2013recipes}. (a) Hopf (red) and saddle-node (black) bifurcation curves and their asymptotes in the $(\epsilon,\mu)$ parameter plane (b) Families of periodic responses (grey surface) in the $(\epsilon,\mu,\max(u_1))$ projection. The equilibrium $u=0$ is stable in the green shaded region.
}\label{15_node_HB_and_SN_curves}
\end{figure*}

We may interpret these observations by considering deformations to the family of limit cycles born at the Hopf bifurcation due to nonzero $\epsilon$ and for different fixed values of the linear damping coefficients. Figure~\ref{15_node_HB_and_SN_curves_a} shows the corresponding Hopf and saddle-node bifurcation curves for the full network model under variations in $\epsilon$ and $\mu$ with $\zeta_i=\mu$ for all $i$ (consistent with the postulated behavior of the time-dependent damping coefficients in the fully coupled dynamics in the limit as $\epsilon\rightarrow 0$). These curves are tangential to the two horizontal lines corresponding to the predicted bifurcation values for $\mu\sim\mathcal{O}(1)$ in the $\epsilon\rightarrow 0$ limit predicted in Sect.~\ref{sec_networkfilters} and given by

\begin{align}
    \mu_{\text{HB}}&=P^2_{Q,I}\nu/(1-P^2_{Q,I})\approx 1.62,\\
    \mu_{\text{SN}}&=P^2_{Q,I}(\nu+\eta/8)/(1-P^2_{Q,I})\approx 3.62. 
\end{align}
Interestingly, as $\epsilon$ varies, the critical values of $\mu$ first increase monotonically with $\epsilon$ but both curves fold back toward decreasing values of $\epsilon$ at local extrema in $\epsilon$. In fact, the subsequent variations are asymptotic to the predicted relationship between $\mu\sim\mathcal{O}(1/\epsilon^2)$ and $\epsilon$ in the limit as $\epsilon\rightarrow 0$ predicted in Sect.~\ref{sec_networkfilters} and given by
\begin{align}
    \mu_{\text{HB}}&=(K_{Q,Q}-1)/K_{Q,Q}\nu\epsilon^2=2/3\epsilon^2,\\
    \mu_{\text{SN}}&=(K_{Q,Q}-1)/K_{Q,Q}(\nu+\eta/8)\epsilon^2=8/27\epsilon^2.
\end{align}

A fuller picture is provided by the three-dimensional representation in Fig.~\ref{15_node_HB_and_SN_curves_b}. The shape of this surface and the saddle-node and Hopf bifurcation curves explains why a larger value of $\zeta$ is reached before $A_1$ drops to $0$ in Fig.~\ref{15_node_time_example_Q1_ep0_1_d}. As long as the value of $\epsilon$ is to the left of that corresponding to the fold in the saddle-node bifurcation curve, it should be possible to obtain a hysteretic response in the full network dynamics, provided that $\delta$ is chosen large enough. This is no longer the case for values of $\epsilon$ to the right of this fold point, since there the family of stable periodic orbits persists across the entire range of parameter values.

\begin{figure*}[htbp]
\centering
\subfloat[]{\includegraphics[width=.33\textwidth]{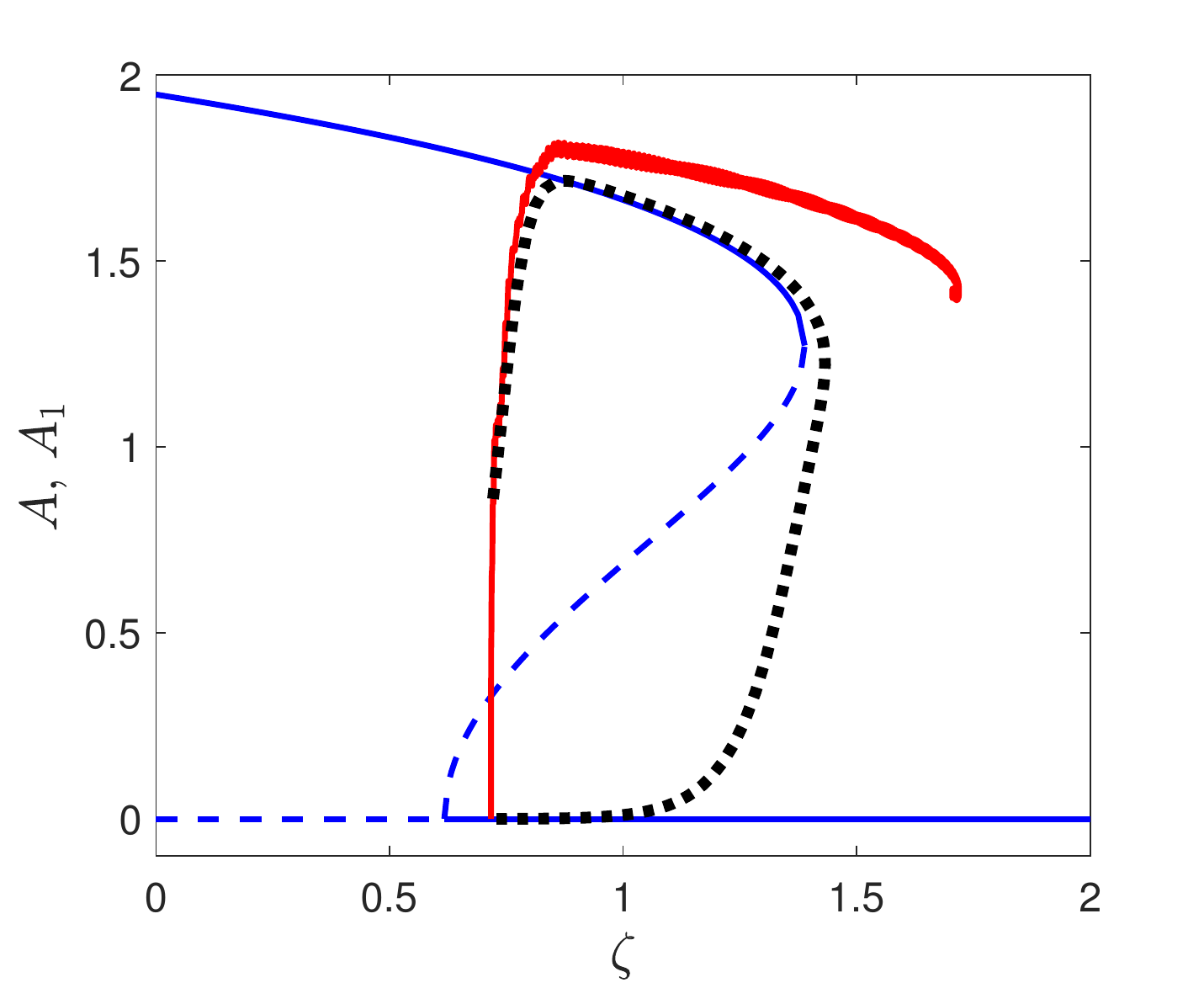}\label{15_node_hysteresis_Q1_a}}
\subfloat[]{\includegraphics[width=.33\textwidth]{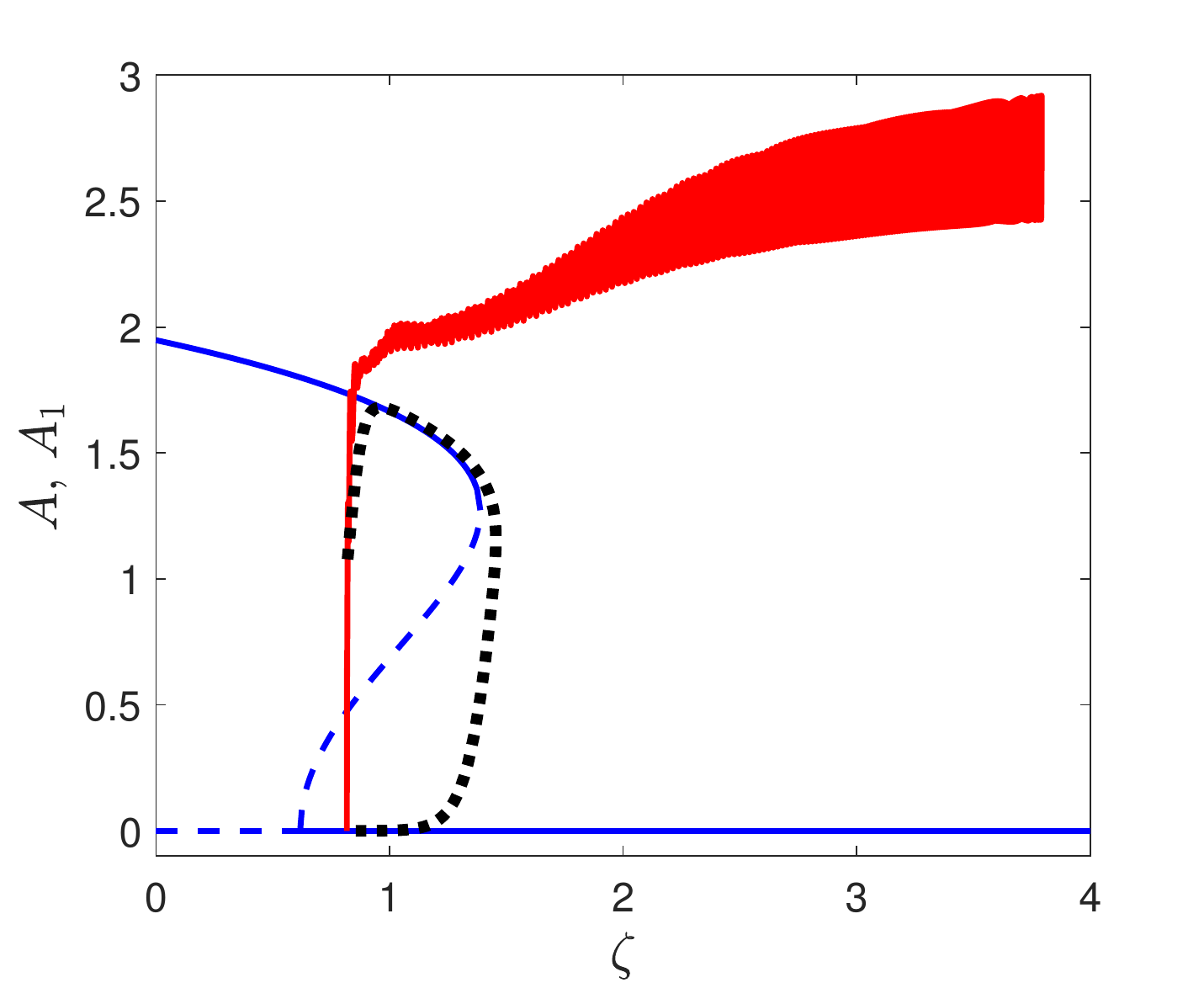}\label{15_node_hysteresis_Q1_b}}
\subfloat[]{\includegraphics[width=.33\textwidth]{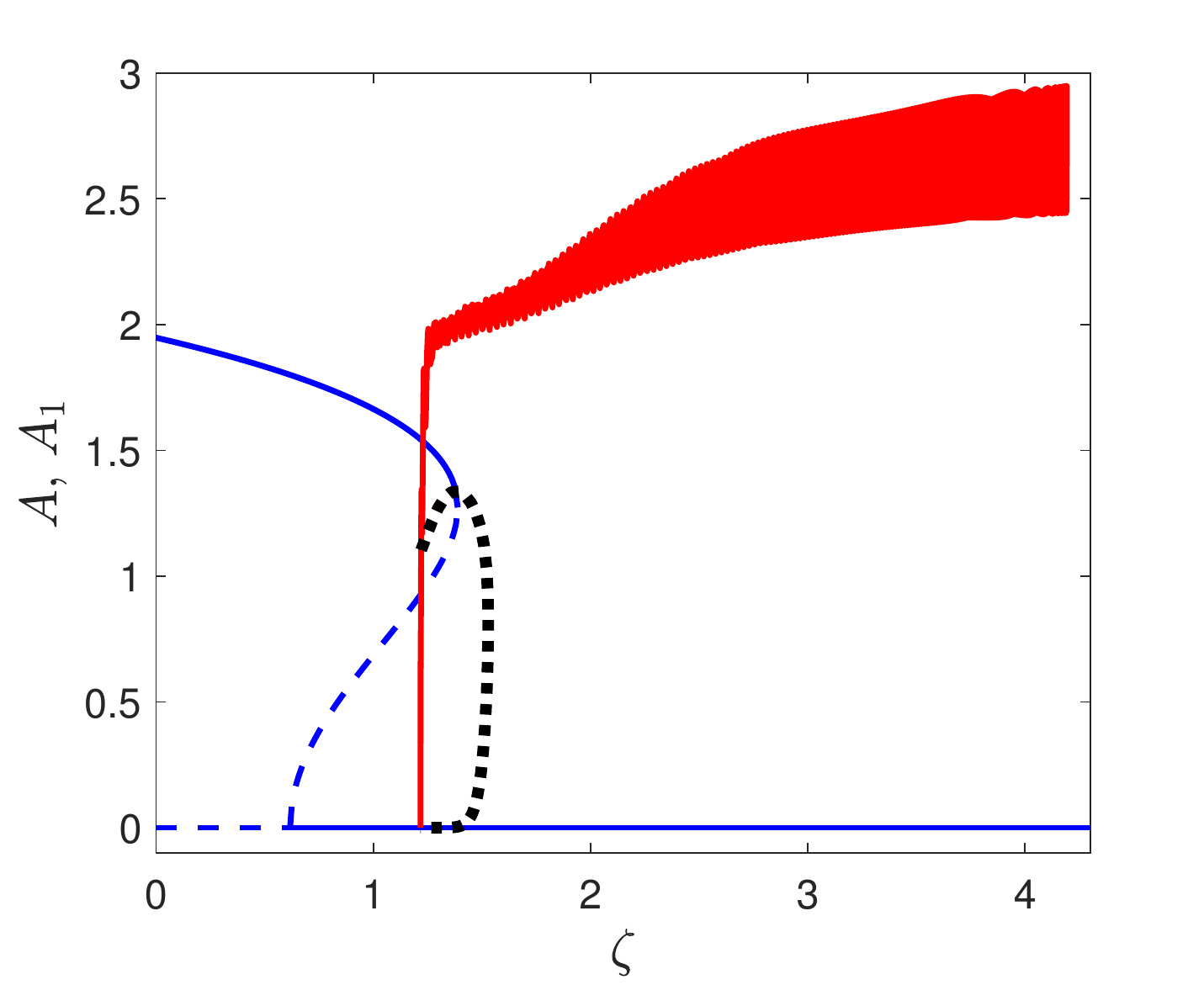}\label{15_node_hysteresis_Q1_c}}
\caption{Non-hysteretic trajectories (red solid) projected onto $\zeta=\sum_{k\ne Q}P_{k,I}^2\zeta_k$ and $A_1$, the trajectory (black dotted) of the coupled $A$ and $\zeta$ dynamics obtained from Eqs.~\eqref{eq:Adyn} and \eqref{eq_smalldampingdynamics}, and the corresponding $A$ nullcline obtained from Eq.~\eqref{eq:Anull}. The full simulation is initialized with zero initial displacements and velocities, $A_i=0$ for $t\in[-2\pi/\omega_I,0]$ and $\zeta_i=(\delta+P^2_{Q,I}\nu)/(1-P^2_{Q,I})$. The simulation of the coupled $A$ and $\zeta$ dynamics is initialized with the values of $A_1$ and $\zeta$ at the conclusion of the initial period of exogenous excitation. Here, $\nu=1$, $\eta=10$, $F = 3\epsilon \sin (\omega_I t)\mathbf{e}_1$ for $0\leq t\leq 1/\epsilon$, and $\tau=20$. (a) $\epsilon=0.1$ and $\delta=0.1$, (b) $\epsilon=0.2$ and $\delta=0.2$, (c) $\epsilon=0.2$ and $\delta=0.6$.}\label{15_node_hysteresis_Q1}
\end{figure*}

These predictions are supported by the results of simulation shown in Fig.~\ref{15_node_hysteresis_Q1}. Hysteresis is not observed with $\epsilon=0.1$ and $\delta=0.1$, and becomes impossible to achieve for $\epsilon=0.2$ for any $\delta$. In each case, the transient trajectory converges to a stable periodic response albeit with a frequency different from $\omega_I$ ($\approx 1.54$ in Fig.~\ref{15_node_hysteresis_Q1_a}, $\approx 1.70$ in Fig.~\ref{15_node_hysteresis_Q1_b}, and $\approx 1.70$ in Fig.~\ref{15_node_hysteresis_Q1_c}) and with non-identical values of $\zeta_i$.

\begin{figure*}[htbp]
\centering
\subfloat[]{\includegraphics[width=.45\textwidth]{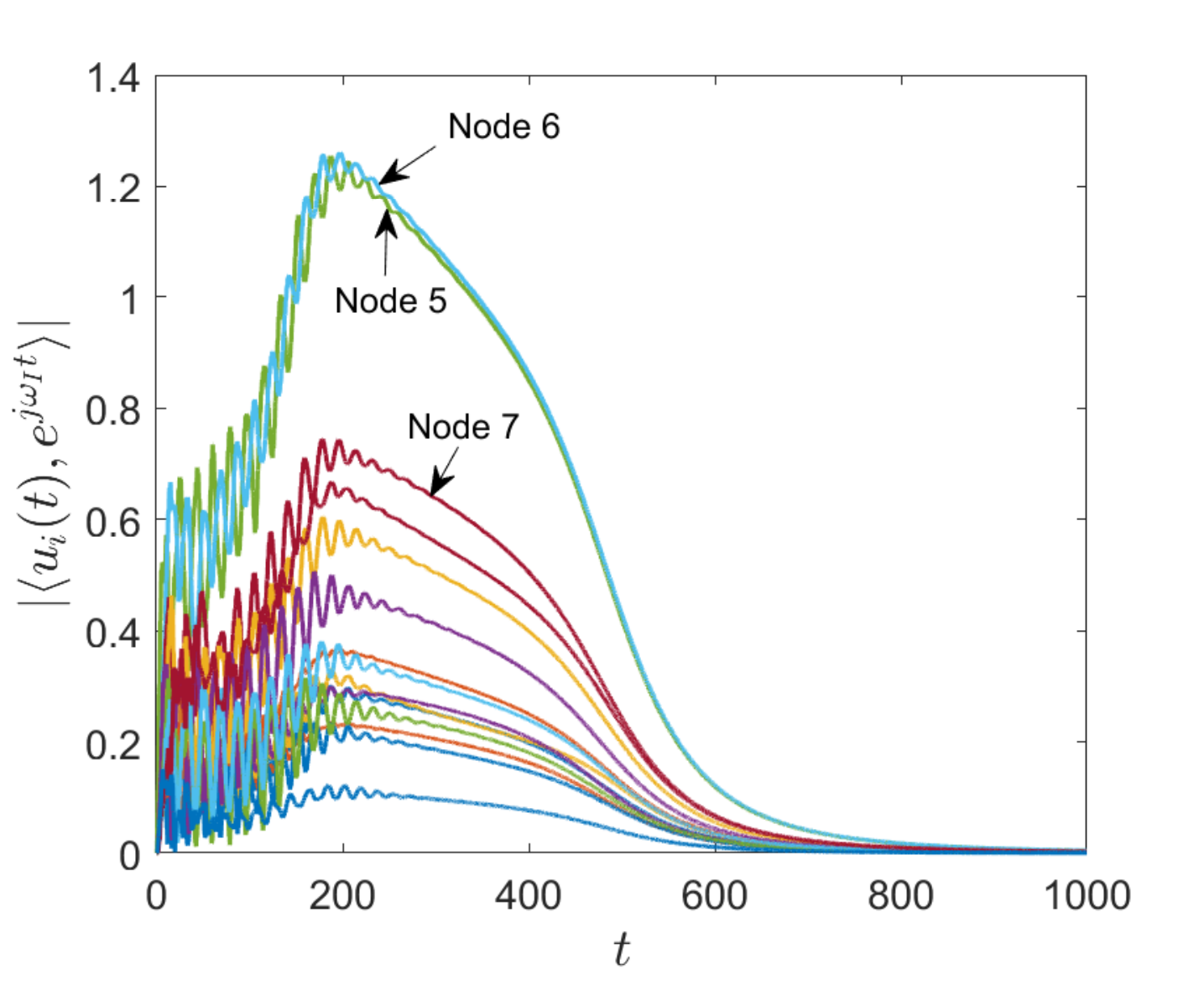}}
\subfloat[]{\includegraphics[width=.45\textwidth]{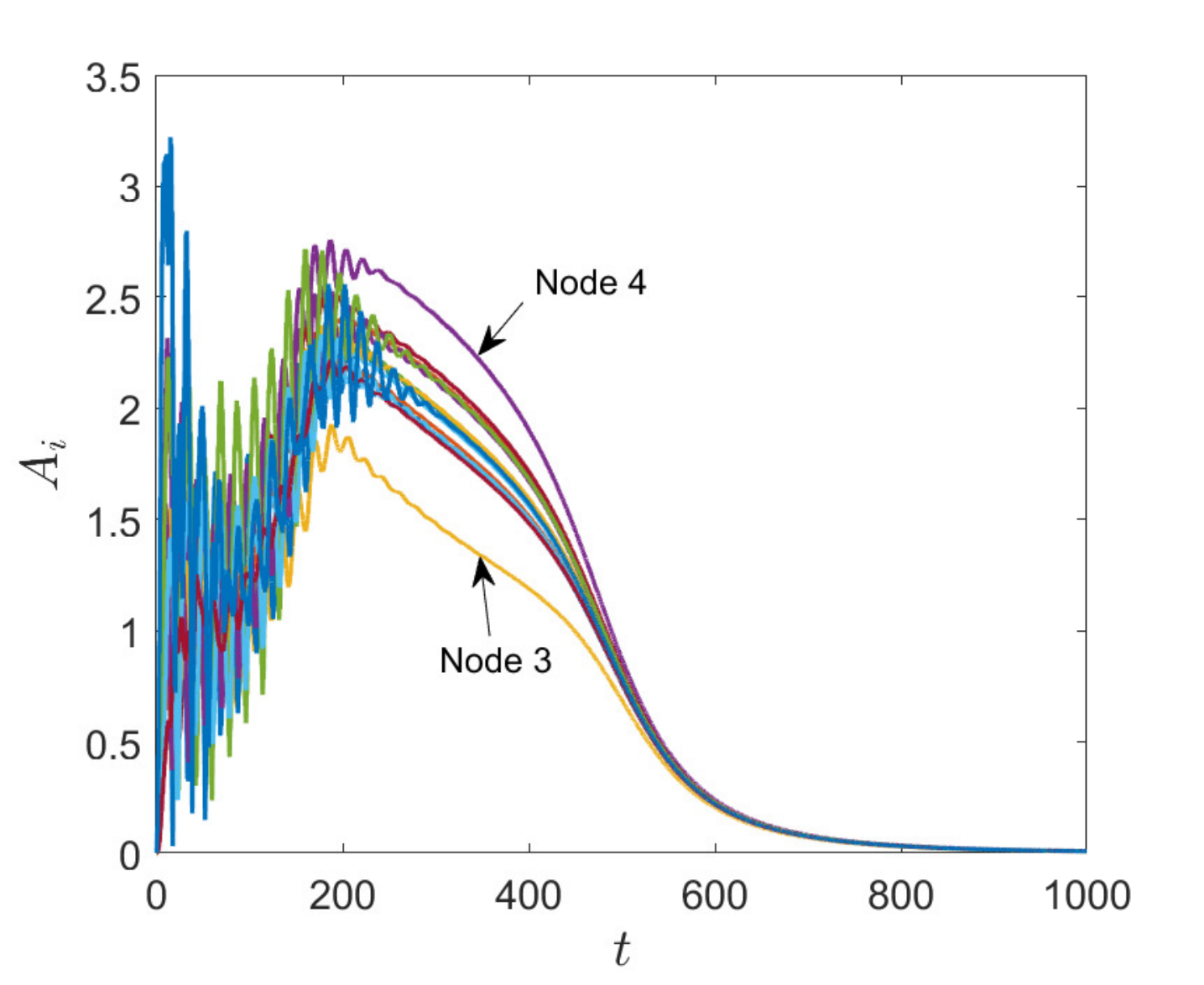}}\\
\subfloat[]{\includegraphics[width=.45\textwidth]{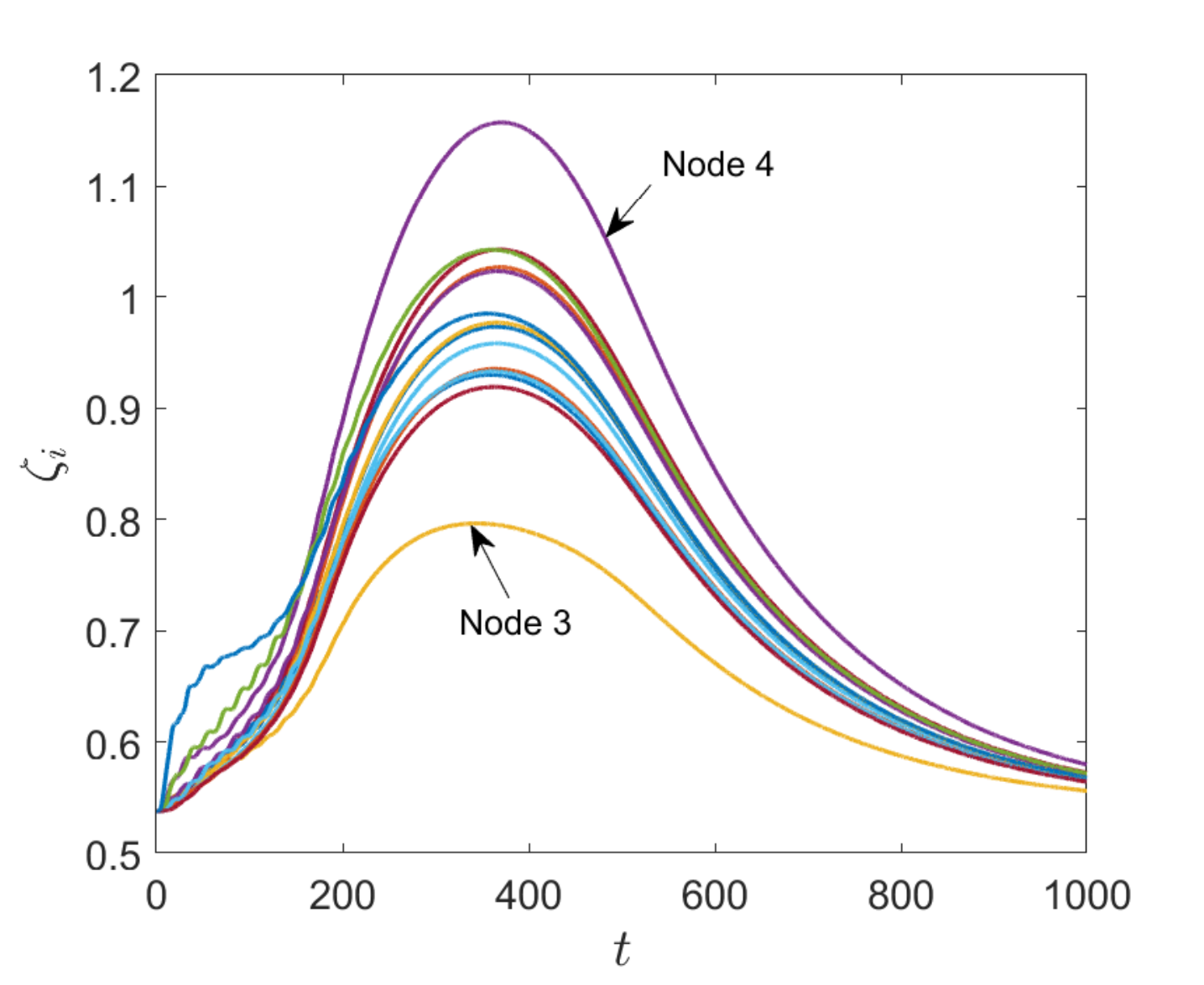}}
\subfloat[]{\includegraphics[width=.45\textwidth]{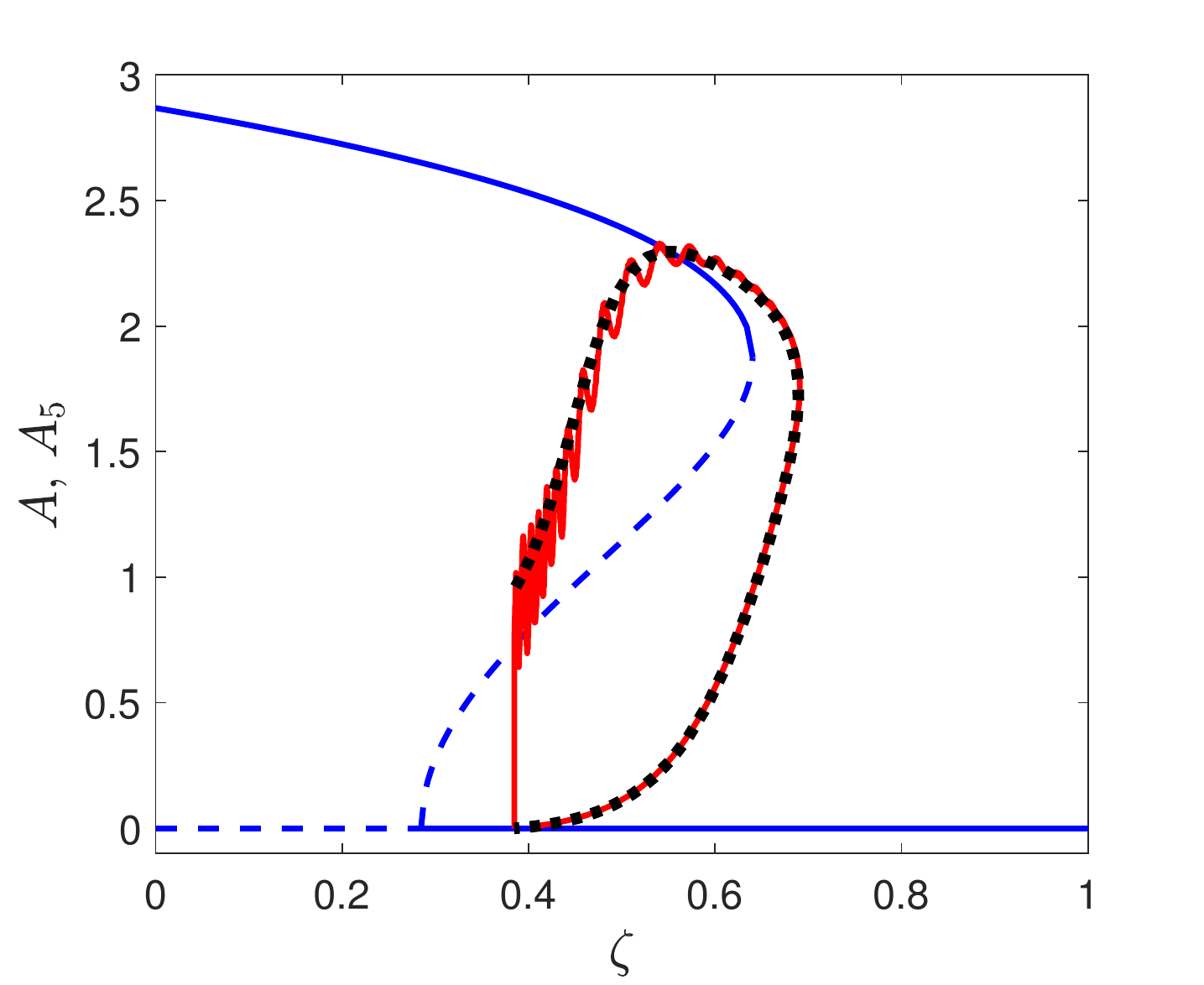}}
\caption{Time histories of (a) nodal amplitudes $|\langle u_i(t),e^{j\omega_I t}\rangle|$, (b) nodal estimates $A_i$, and (c) nodal damping coefficients $\zeta_i$ for the 15-node network with topology shown in Fig.~\ref{15node_network_example} with $Q=5$, $I=12$, and $\omega_I\approx 2.8428$. Panel (d) shows the corresponding hysteretic trajectory (red solid) projected onto $\zeta=\sum_{k\ne Q}P_{k,I}^2\zeta_k$ and $A_5$, the trajectory (black dotted) of the coupled $A$ and $\zeta$ dynamics obtained from Eqs.~\eqref{eq:Adyn} and \eqref{eq_smalldampingdynamics}, and the corresponding $A$ nullcline obtained from Eq.~\eqref{eq:Anull}. The full simulation is initialized with zero initial displacements and velocities, $A_i=0$ for $t\in[-2\pi/\omega_I,0]$ and $\zeta_i=(\delta+P^2_{Q,I}\nu)/(1-P^2_{Q,I})$. The simulation of the coupled $A$ and $\zeta$ dynamics is initialized with the values of $A_5$ and $\zeta$ at the conclusion of the initial period of exogenous excitation. Here, $\epsilon=0.1$, $\nu=1$, $\eta=10$, $F = 9\epsilon \sin (\omega_I t)\mathbf{e}_5$ for $0\leq t\leq 1/\epsilon$, $\delta=0.1$, and $\tau=20$. }\label{15_node_time_example_Q5_ep0_1}
\end{figure*}

\begin{figure*}[htbp]
\centering
\subfloat[]{\includegraphics[width=.45\textwidth]{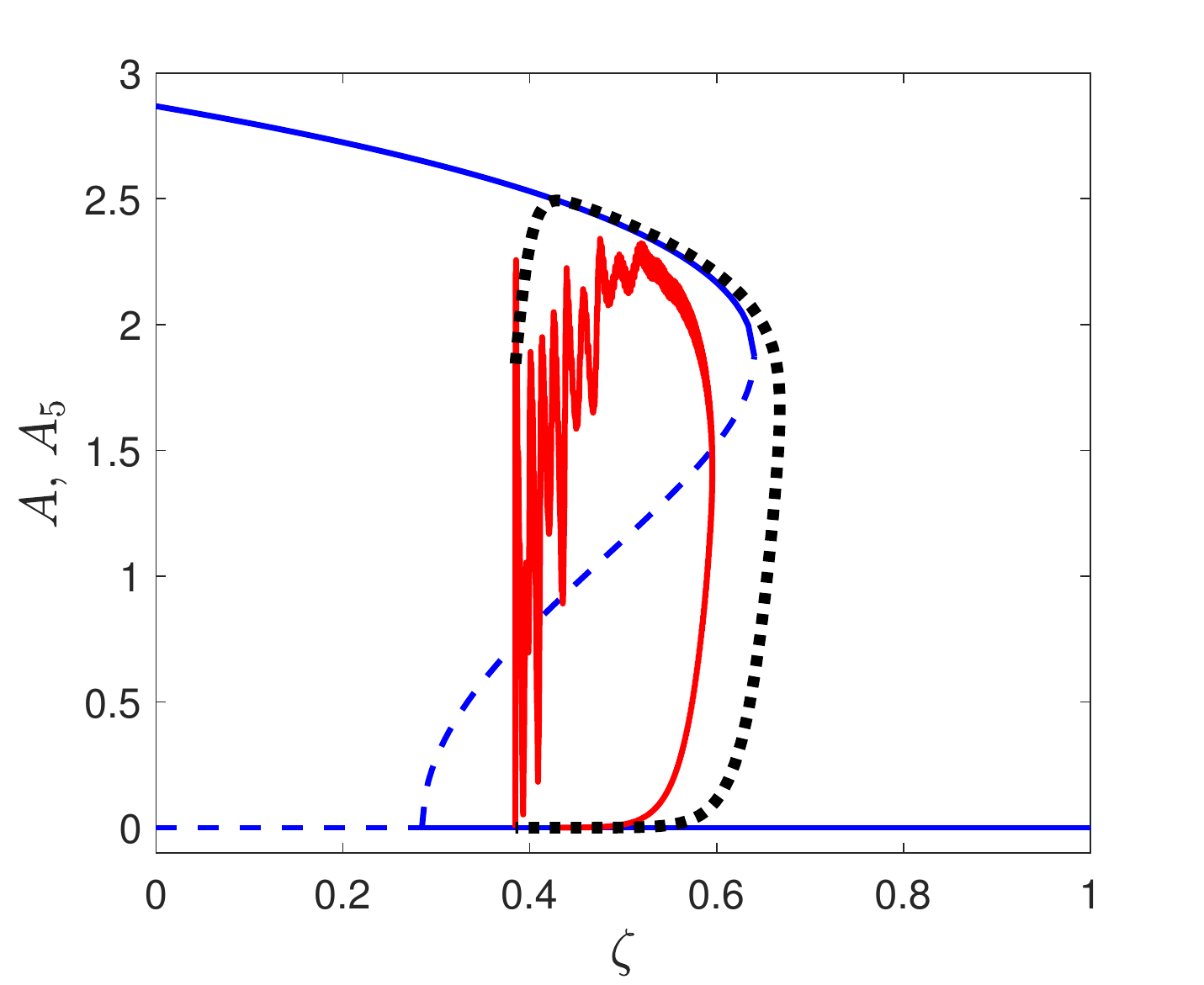}\label{15_node_hysteresis_Q5_a}}
\subfloat[]{\includegraphics[width=.45\textwidth]{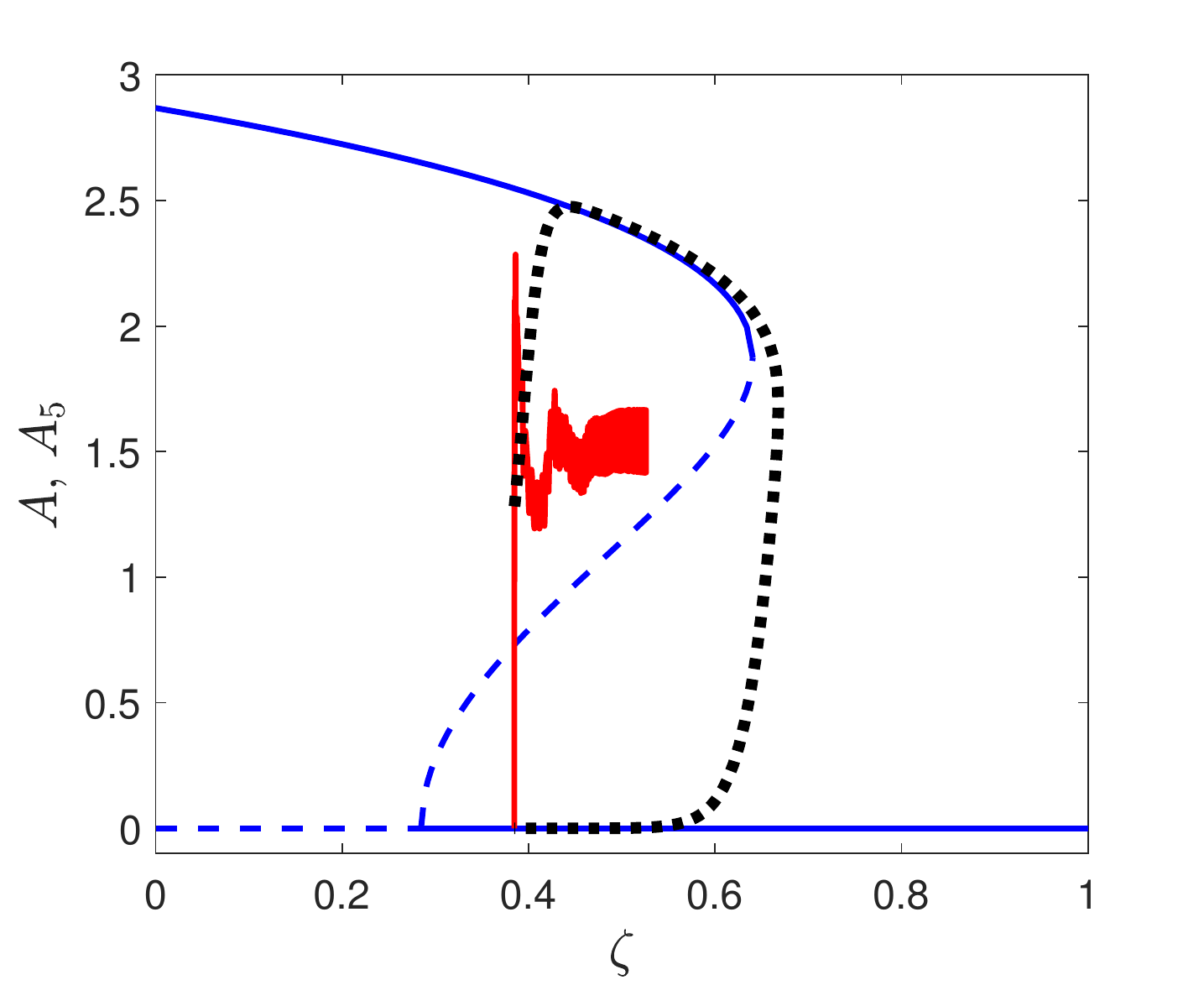}\label{15_node_hysteresis_Q5_b}}
\caption{Hysteretic and non-hysteretic trajectories (red solid) for the 15-node network with topology shown in Fig.~\ref{15node_network_example} with $Q=5$, $I=12$, and $\omega_I\approx 2.8428$, projected onto $\zeta=\sum_{k\ne Q}P_{k,I}^2\zeta_k$ and $A_5$, the trajectory (black dotted) of the coupled $A$ and $\zeta$ dynamics obtained from Eqs.~\eqref{eq:Adyn} and \eqref{eq_smalldampingdynamics}, and the corresponding $A$ nullcline obtained from Eq.~\eqref{eq:Anull}. The full simulation is initialized with zero initial displacements and velocities, $A_i=0$ for $t\in[-2\pi/\omega_I,0]$ and $\zeta_i=(\delta+P^2_{Q,I}\nu)/(1-P^2_{Q,I})$. The simulation of the coupled $A$ and $\zeta$ dynamics is initialized with the values of $A_5$ and $\zeta$ at the conclusion of the initial period of exogenous excitation. Here, $\nu=1$, $\eta=10$, $F = 9\epsilon \sin (\omega_I t)\mathbf{e}_1$ for $0\leq t\leq 1/\epsilon$, and $\tau=60$. (a) $\epsilon=0.3$ and $\delta=0.1$, (b) $\epsilon=0.4$ and $\delta=0.1$.}\label{15_node_hysteresis_Q5}
\end{figure*}

The largest value of $\epsilon$ that allows for hysteretic behavior varies as the nonlinear oscillator is put at different locations. When $Q=5$, the system can exhibit hysteretic behavior and a small spread for the time histories of $\zeta_i$ even for $\epsilon=0.1$ and $\delta=0.1$, as shown in Fig.~\ref{15_node_time_example_Q5_ep0_1}. In this case, $I=12$ and periodic orbits with limiting frequency $\omega_I\approx 2.84$ are born from the Hopf bifurcation for $\zeta_i=\mathcal{O}(1)$ in the limit as $\epsilon\rightarrow0$. The hysteresis loop in the $(\zeta,A)$ plane tracks closely the trajectory from the perturbation analysis apart from transient oscillations due to excitation of another nearby natural frequency ($\omega_9\approx 2.49$) that become negligible as $\epsilon$ decreases. The system is still able to exhibit hysteretic behavior when $\epsilon=0.3$, as shown in Fig.~\ref{15_node_hysteresis_Q5_a}, albeit with a more significant deviation from the prediction of the multiple-scale analysis. Hysteresis is not observed with $\epsilon=0.4$ and the transient trajectory converges to a stable periodic response with a frequency different from $\omega_I$ ($\approx 2.34$ in Fig.~\ref{15_node_hysteresis_Q5_b}). 

\begin{figure*}[htbp]
\centering
\subfloat[]{\includegraphics[width=.49\textwidth]{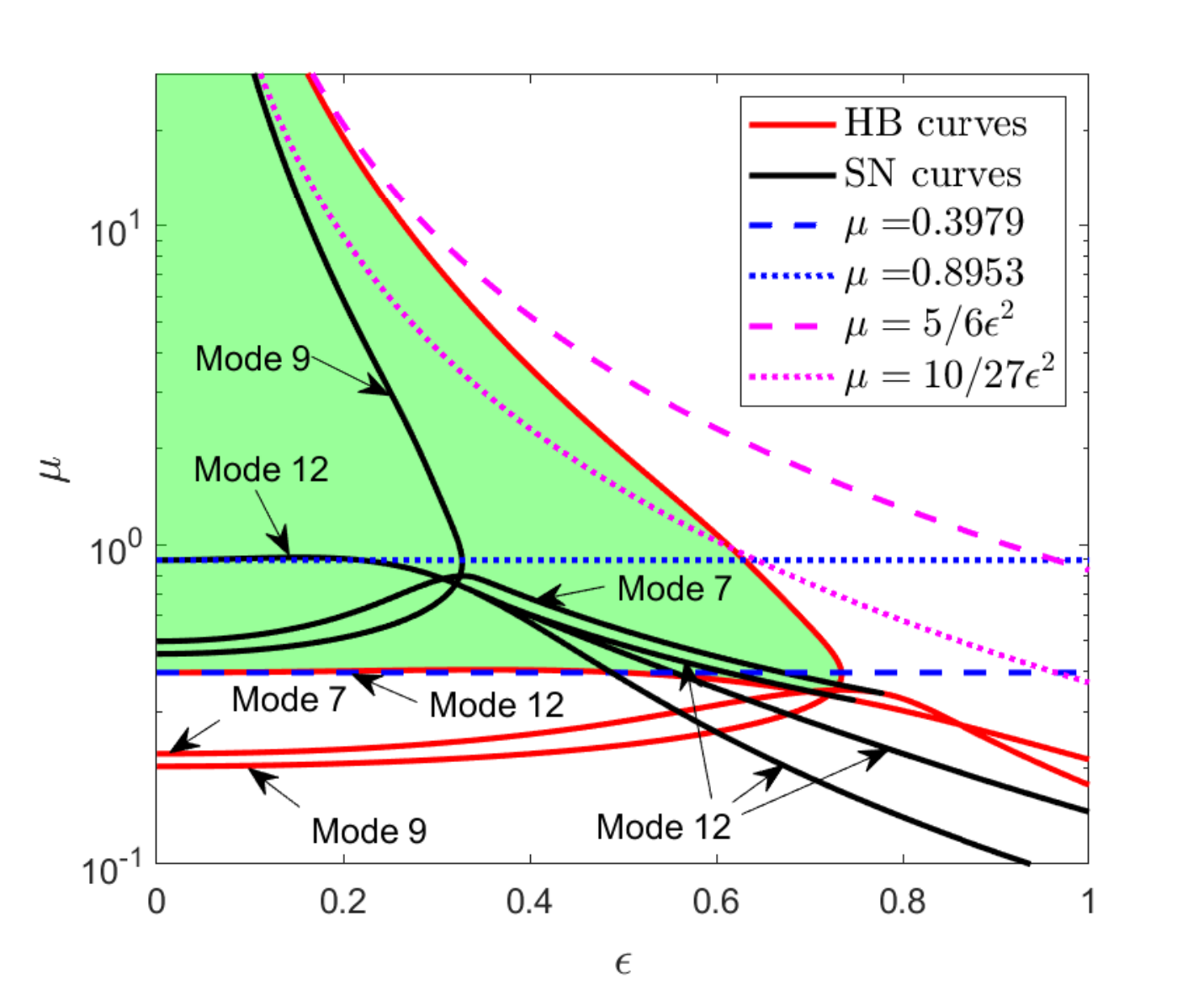}\label{15_node_Q5_HB_and_SN_curves_a}}
\subfloat[]{\includegraphics[width=.49\textwidth]{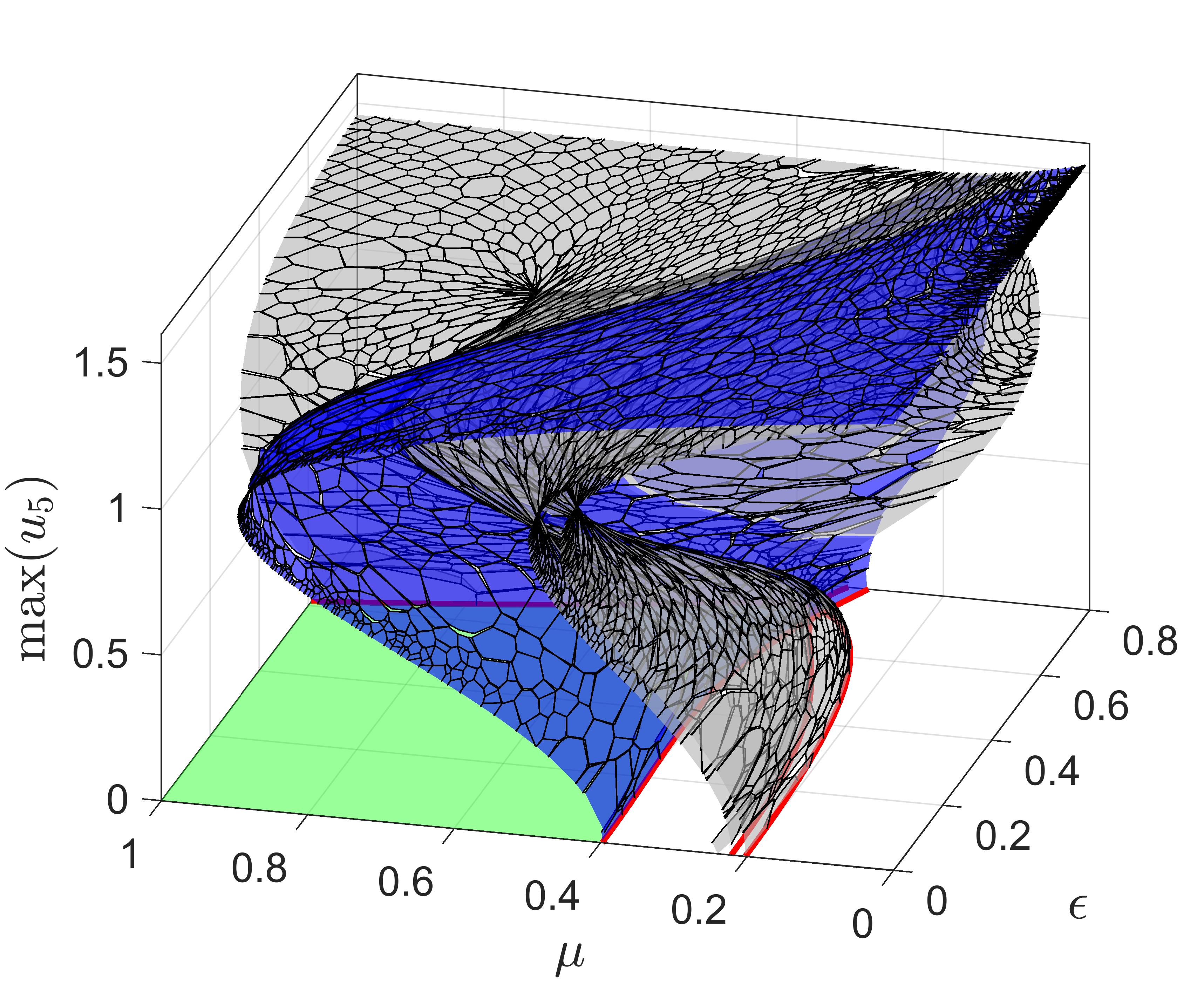}\label{15_node_Q5_HB_and_SN_curves_b}}
\caption{Bifurcation analysis for the 15-node network example shown in Fig.~\ref{15node_network_example} with $Q=5$ and $\zeta_i=\mu$ for $i\neq Q$ obtained using the equilibrium (\texttt{eq}) and periodic-orbit (\texttt{po}) toolboxes of the \textsc{coco} software platform~\cite{dankowicz2013recipes}. (a) Hopf (red) and saddle-node (black) bifurcations and select asymptotes in the $(\epsilon,\mu)$ parameter plane. Two of the SN curves associated with mode 12 appear in a cusp bifurcation, and one of them disappears when the amplitude becomes $0$. (b) Families of periodic responses (grey surface is associated with modes 7 and 9, blue surface is associated with mode 12) in the $(\epsilon,\mu,\max(u_5))$ projection. The equilibrium $u=0$ is stable in the green shaded region. }\label{15_node_Q5_HB_and_SN_curves}
\end{figure*}

\begin{figure*}[htbp]
\centering
\subfloat[]{\includegraphics[width=.49\textwidth]{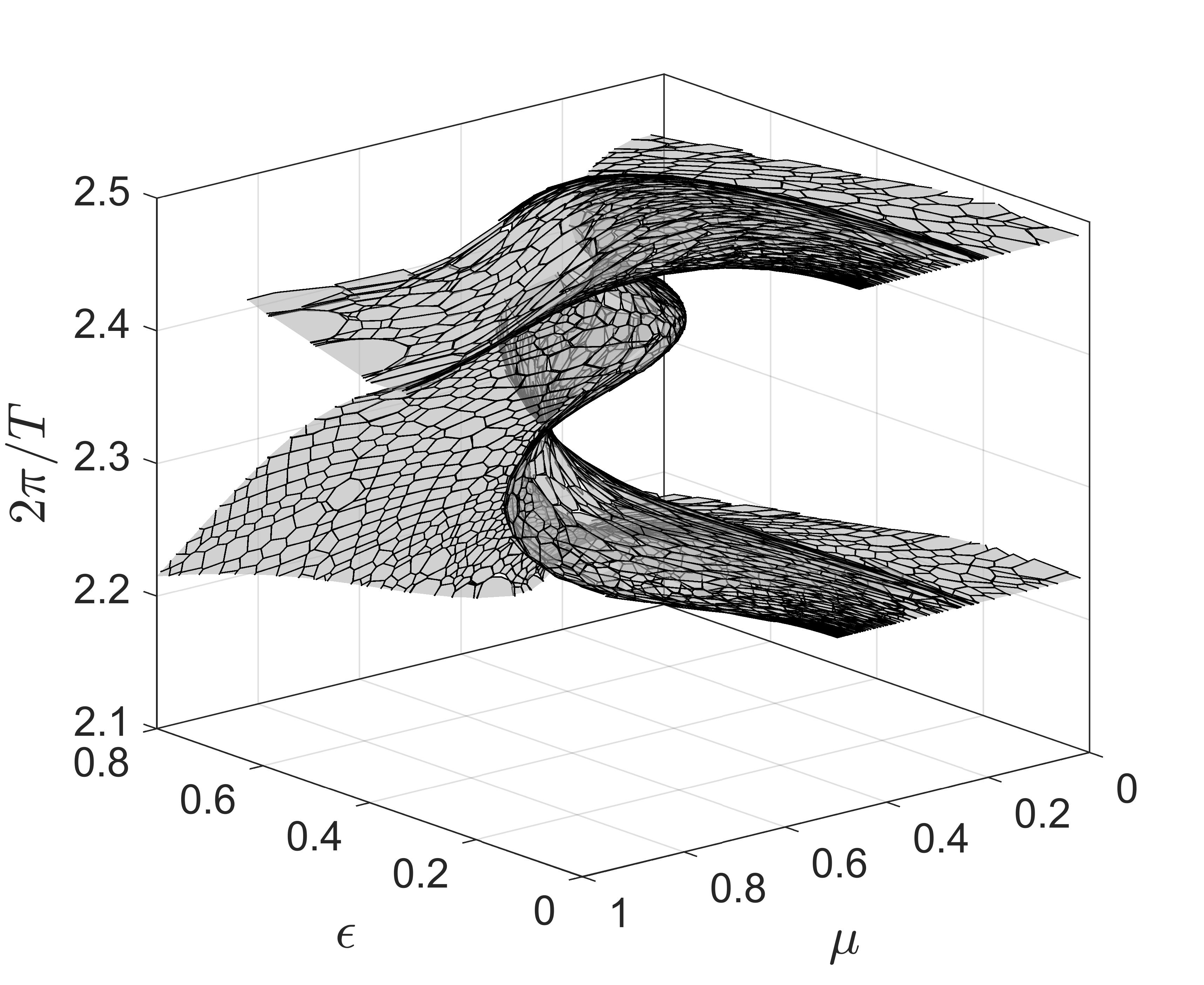}}
\subfloat[]{\includegraphics[width=.49\textwidth]{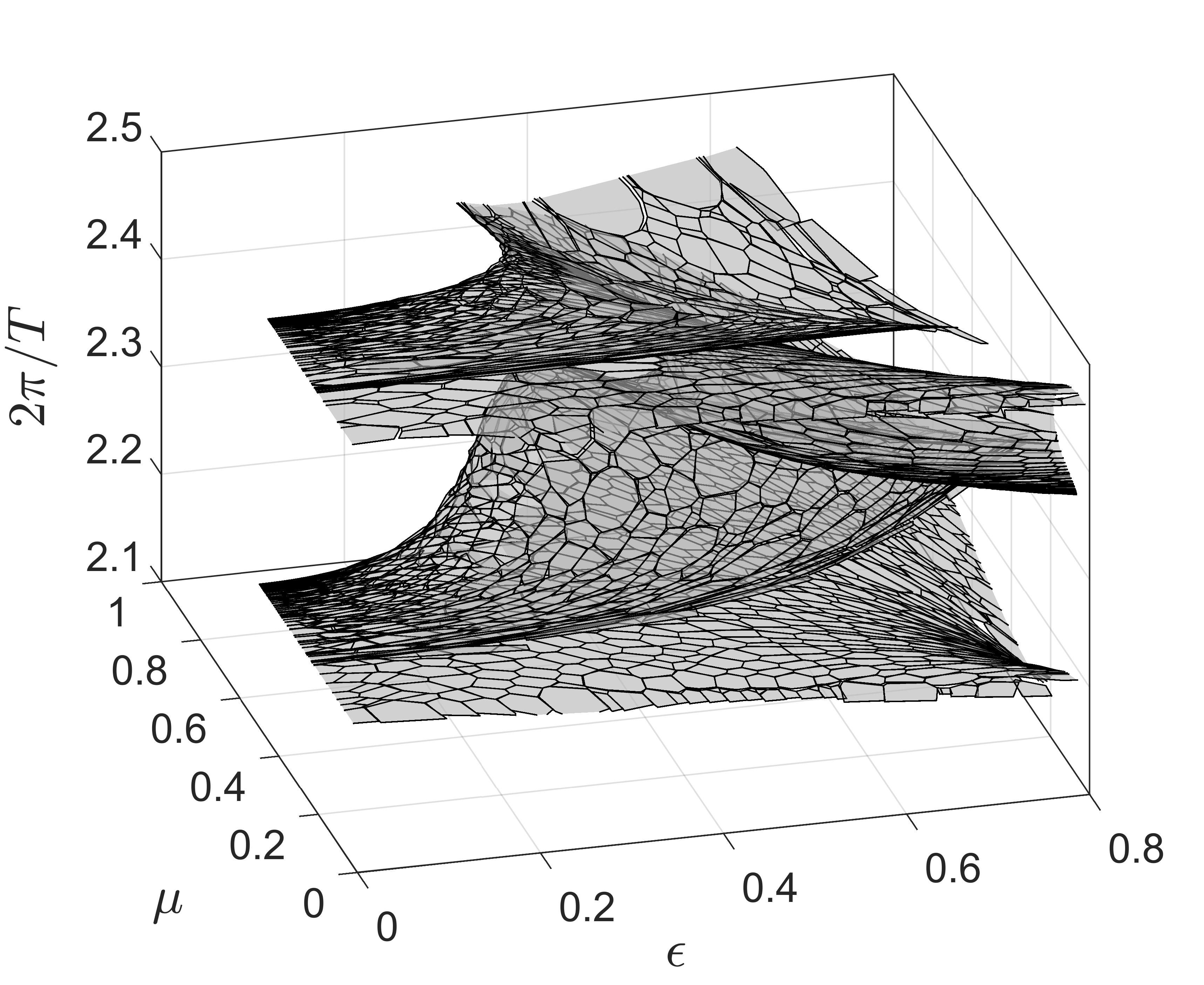}}
\caption{Surface of periodic responses of period $T$ in the $(\epsilon,\mu,2\pi/T)$ projection for the 15-node network example shown in Fig.~\ref{15node_network_example} with $Q=5$ and $\zeta_i=\mu$ for $i\neq Q$ obtained using the periodic-orbit (\texttt{po}) toolbox of the \textsc{coco} software platform~\cite{dankowicz2013recipes}. In the limit as $\epsilon\rightarrow 0$ the surface asymptotes to families of limit cycles associated with modes 7 and 9.}\label{15_node_surf_Q5_large_freq}
\end{figure*}

Figure~\ref{15_node_Q5_HB_and_SN_curves_a} shows several Hopf and saddle-node bifurcation curves under variations in $\epsilon$ and $\mu$ (again assuming $\zeta_i=\mu$ for all $i$) together with select asymptotes predicted by the multiple-scale analysis for the $\mathcal{O}(1)$ and $\mathcal{O}(1/\epsilon^2)$ limits. We have labeled the curves by the corresponding linear mode in the limit as $\epsilon\rightarrow 0$ and $\mu=\mathcal{O}(1)$.
In contrast to the case with $Q=1$, the region of stability for the trivial equilibrium is here bounded by a piecewise-defined contour obtained from two distinct Hopf bifurcation curves associated with modes 12 and 9, respectively. We note that for $\epsilon\approx0.3$, the saddle-node bifurcation associated with mode 12 is closer to the corresponding Hopf bifurcation curve (which is almost horizontal over a large range of values of $\epsilon$) than in the $\epsilon\rightarrow 0$ limit, explaining why $A_5$ drops to $0$ before the predicted fold in Fig.~\ref{15_node_hysteresis_Q5_a}. A three-dimensional representation of the corresponding families of periodic responses projected onto $(\epsilon,\mu,\max(u_5))$ is shown in Fig.~\ref{15_node_Q5_HB_and_SN_curves_b}. Interestingly, the surfaces of periodic orbits associated with modes 7 and 9 in the $\epsilon\rightarrow 0$ limit connect smoothly at large $\mu$, as shown by the projection in Fig.~\ref{15_node_surf_Q5_large_freq} onto $(\epsilon,\mu,\omega)$ in terms of the angular frequency $\omega$. The saddle-node and Hopf bifurcation curves associated with mode 12 intersect near $\epsilon=0.56$, corresponding to a loss of bistability and making the desired hysteretic behavior impossible beyond this point. In fact, interference from stable periodic orbits on the surfaces associated with modes 7 and 9 makes it difficult to achieve hysteresis already past the intersection near $\epsilon=0.32$ of the saddle-node bifurcation curves associated with modes 12 and 9. Indeed, beyond this point, stable periodic orbits are found across the entire range of damping values.

\begin{figure*}[htbp]
\centering
\subfloat[]{\includegraphics[width=.45\textwidth]{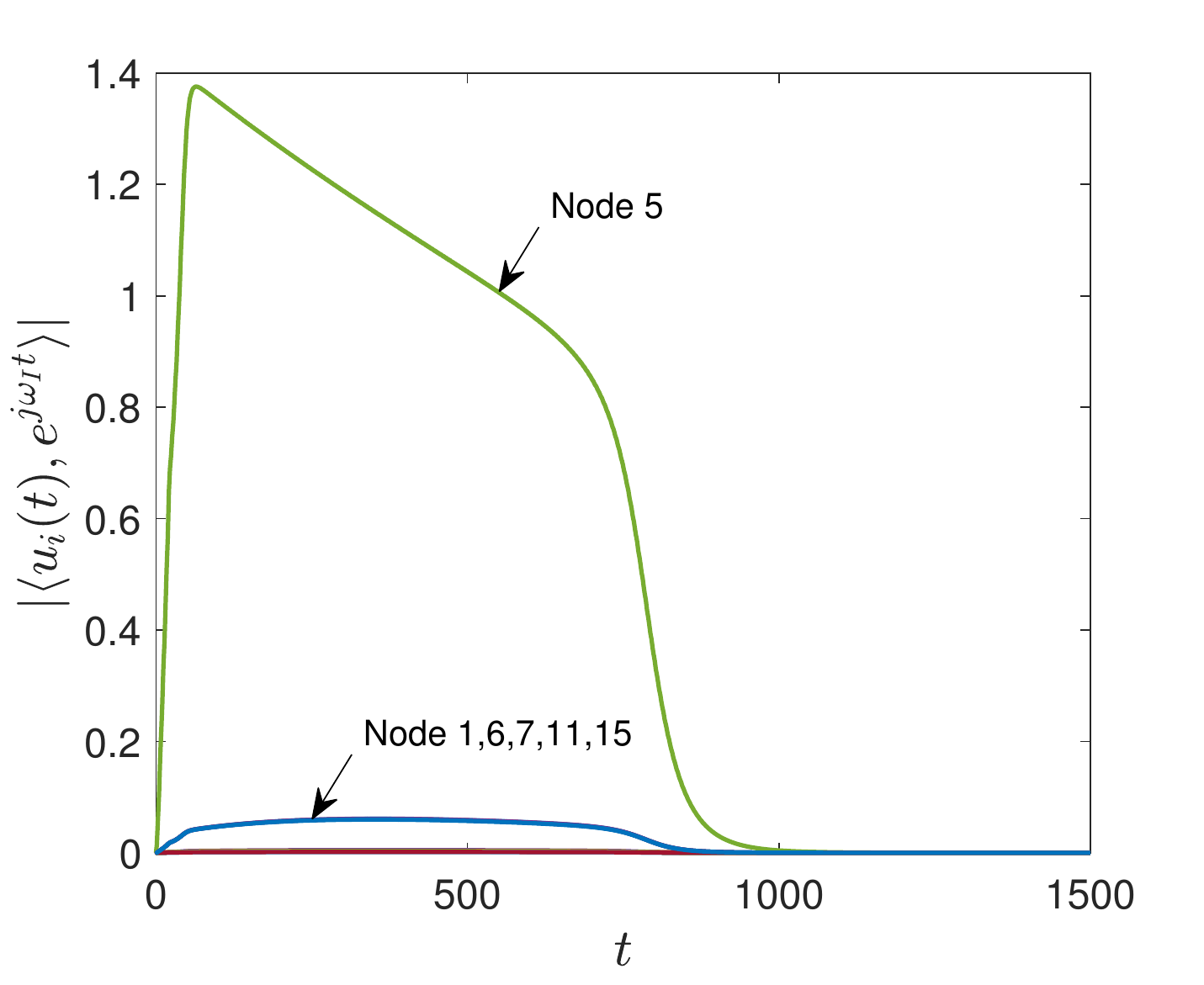}}
\subfloat[]{\includegraphics[width=.45\textwidth]{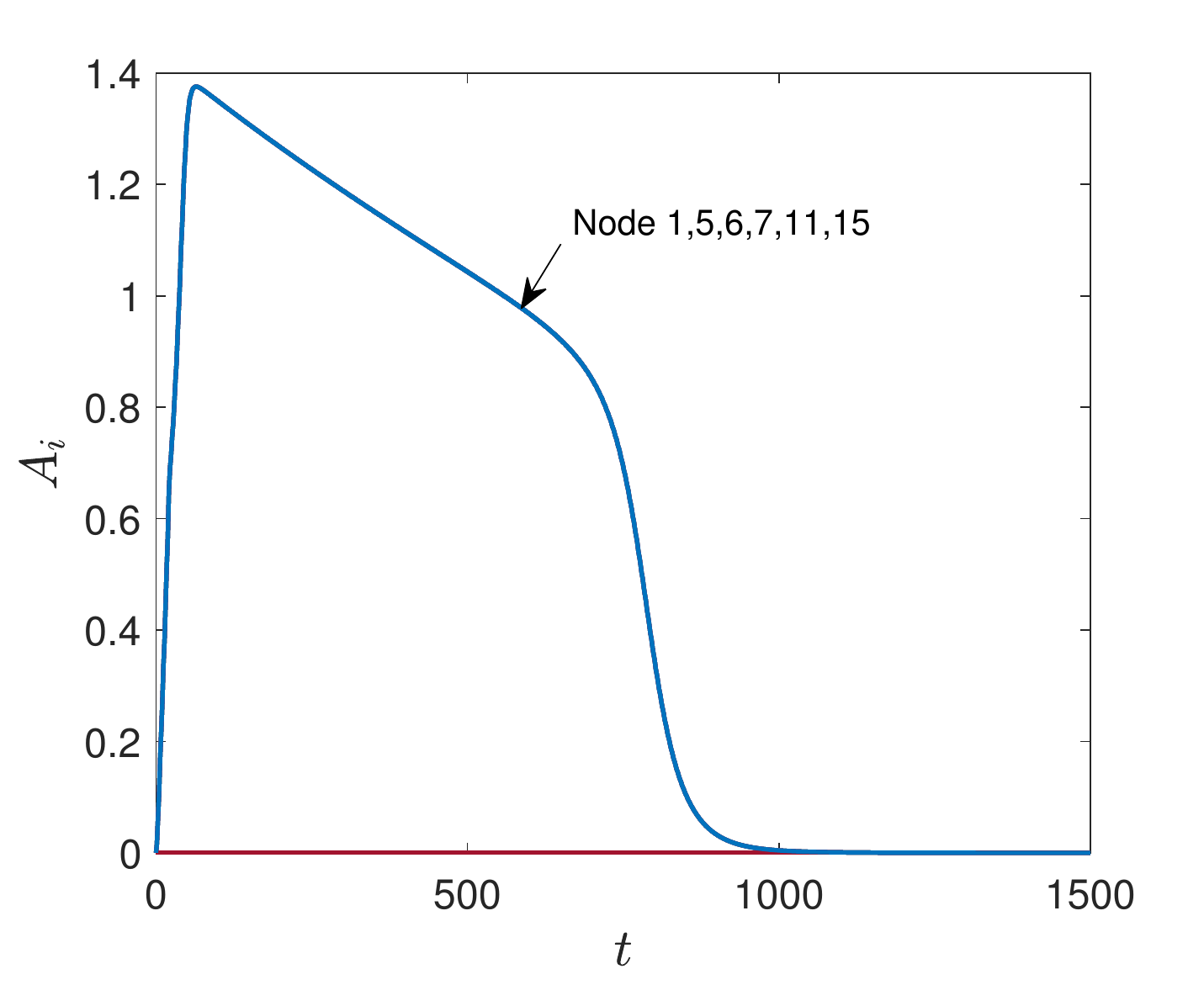}}\\
\subfloat[]{\includegraphics[width=.45\textwidth]{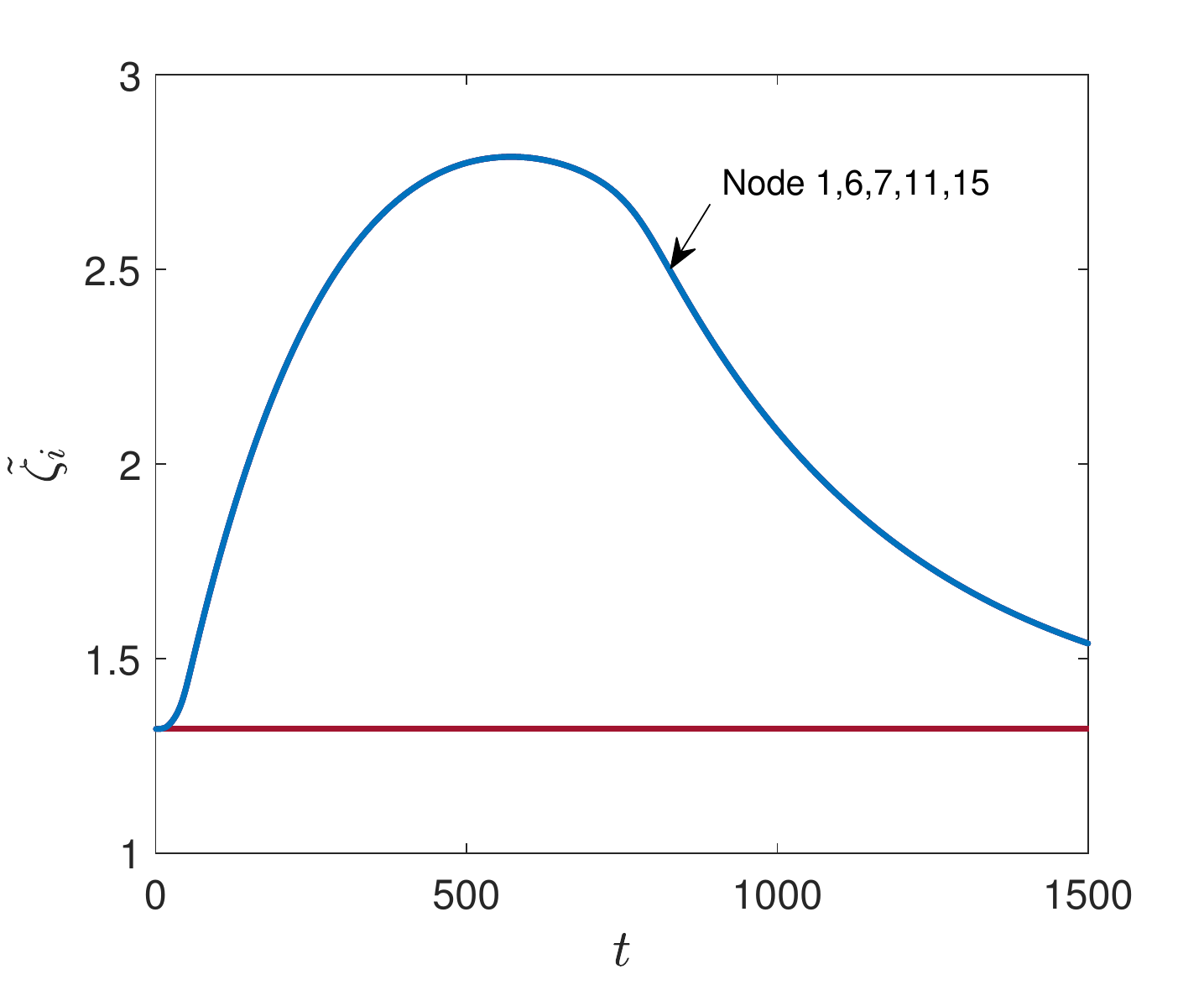}}
\subfloat[]{\includegraphics[width=.45\textwidth]{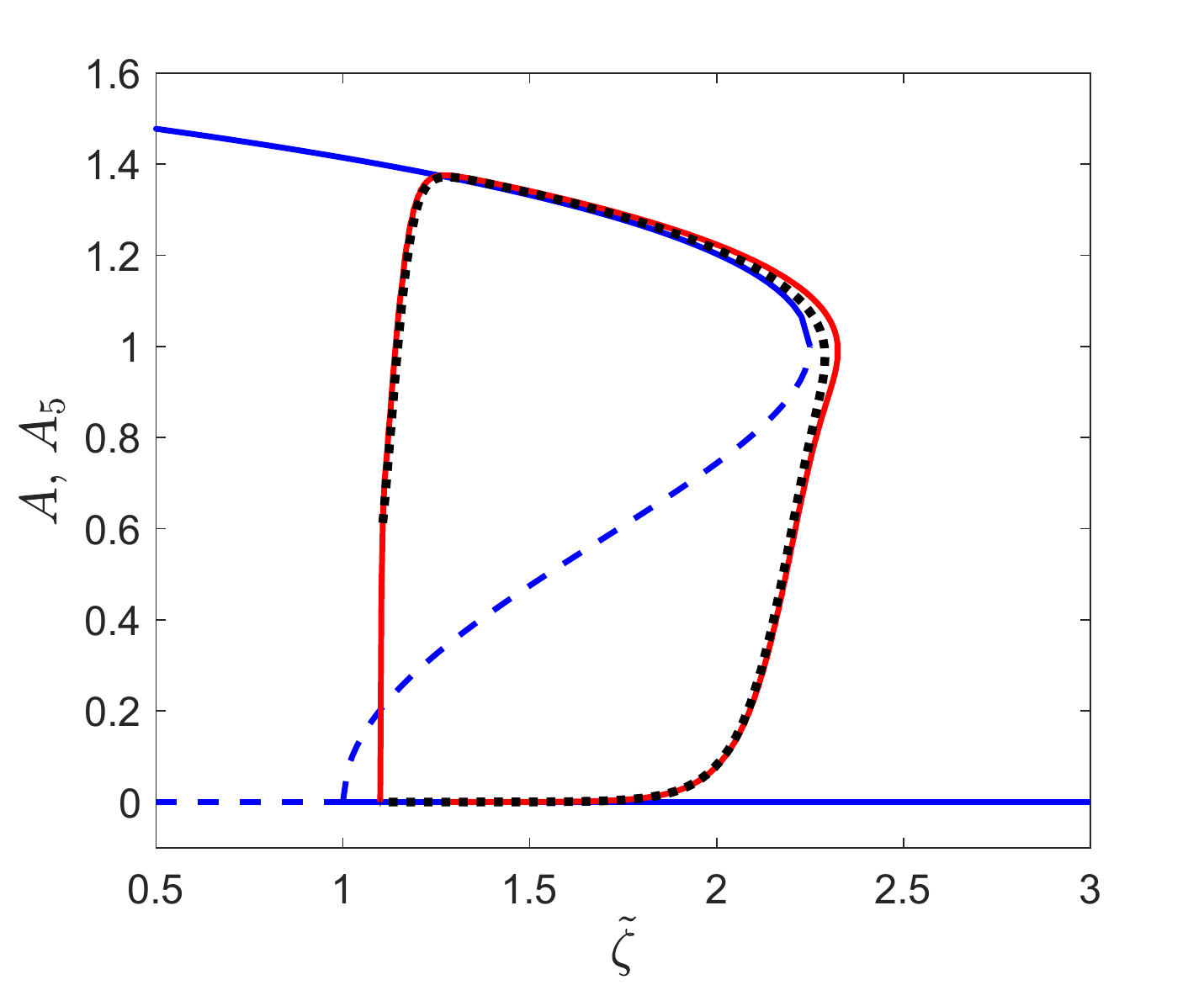}}
\caption{Time histories of (a) nodal amplitudes $|\langle u_i(t),e^{j\omega_I t}\rangle|$, (b) nodal estimates $A_i$, and (c) rescaled nodal damping coefficients $\tilde{\zeta}_i$ for the 15-node network with topology shown in Fig.~\ref{15node_network_example} with $Q=5$ and $\sqrt{K_{Q,Q}}=\sqrt{6}$. Panel (d) shows the corresponding hysteretic trajectory (red solid) projected onto $\tilde{\zeta}=\sum_{k\ne Q}K_{Q,k}K_{k,Q}\tilde{\zeta}_k/K_{Q,Q}$ and $A_5$, the trajectory (black dotted) of the coupled $A$ and $\tilde{\zeta}$ dynamics obtained from Eqs.~\eqref{eq:ampeqhighdamp1} and \eqref{eq:ampeqhighdamp3}, and the corresponding $A$ nullcline obtained from Eq.~\eqref{eq:ampeqhighdamp2}. The full simulation is initialized with zero initial displacements and velocities, $A_i=0$ for $t\in[-2\pi/\omega_I,0]$ and $\tilde{\zeta_i}=K_{Q,Q}(\delta+\nu)/(K_{Q,Q}-1)$. The simulation of the coupled $A$ and $\tilde{\zeta}$ dynamics is initialized with the values of $A_5$ and $\tilde{\zeta}$ at the conclusion of the initial period of exogenous excitation. Here, $\epsilon=0.05$, $\nu=1$, $\eta=10$, $F = 3\epsilon \sin \sqrt{6}\mathbf{e}_5$ for $0\leq t\leq 1/\epsilon$, $\delta=0.1$, and $\tau=20$. }\label{15_node_time_example_Q5_ep0_05_large}
\end{figure*}

Finally, as predicted in Sect.~\ref{sec_networkfilters}, for $\zeta_i=\mathcal{O}(1/\epsilon^2)$, the trivial equilibrium for the 15-node network example with $Q=5$ loses stability at a Hopf bifurcation out of which emanates a family of periodic orbits with limiting frequency $\omega=\sqrt{6}$ as $\epsilon\rightarrow 0$. Figure~\ref{15_node_time_example_Q5_ep0_05_large} shows a simulated hysteretic trajectory in this asymptotic limit with $\epsilon=0.05$. We observe close agreement with the predictions of the perturbation analysis. Per the network topology, nonzero nodal estimates occur only for nodes 1, 6, 7, 11 and 15 and these same nodes are associated with slowly-varying, rescaled damping coefficients. 

%%%%%%%%%%%%%%%%%%%%%%%%%%%%%%%%%%%%%%%%%%%%%%%%%%%%%%%%%%%%%%%%%%%%%%%%%%

\section{Parameter robustness}\label{sec_arameterrobustness}

Based on the results in Sect.~\ref{sec:examples}, we hypothesize that the desired hysteretic response is possible for values of $\epsilon$ smaller than a maximal value along a piecewise-defined contour composed from saddle-node bifurcation curves under simultaneous variations of $\epsilon$ and a suitably defined damping parameter $\mu$. In Fig.~\ref{15_node_HB_and_SN_curves}, we assigned $\zeta_i=\mu$ for all $i\ne Q$ and obtained a smooth contour given by a single saddle-node bifurcation curve associated with limit cycles with angular frequencies $\approx 1.52$ and $\sqrt{3}$ in the asymptotic limits $\mu=\mathcal{O}(1)$ and $\mu=\mathcal{O}(1/\epsilon^2)$, respectively. In contrast, in Fig.~\ref{15_node_Q5_HB_and_SN_curves}, with the same definition for $\zeta_i$, the contour was composed from two saddle-node bifurcation curves associated with limit cycles with angular frequencies $\approx 2.84$ and $\sqrt{6}$ in the asymptotic limits $\mu=\mathcal{O}(1)$ and $\mu=\mathcal{O}(1/\epsilon^2)$, respectively.

\begin{figure}[htbp]
\centering
\includegraphics[width=.45\textwidth]{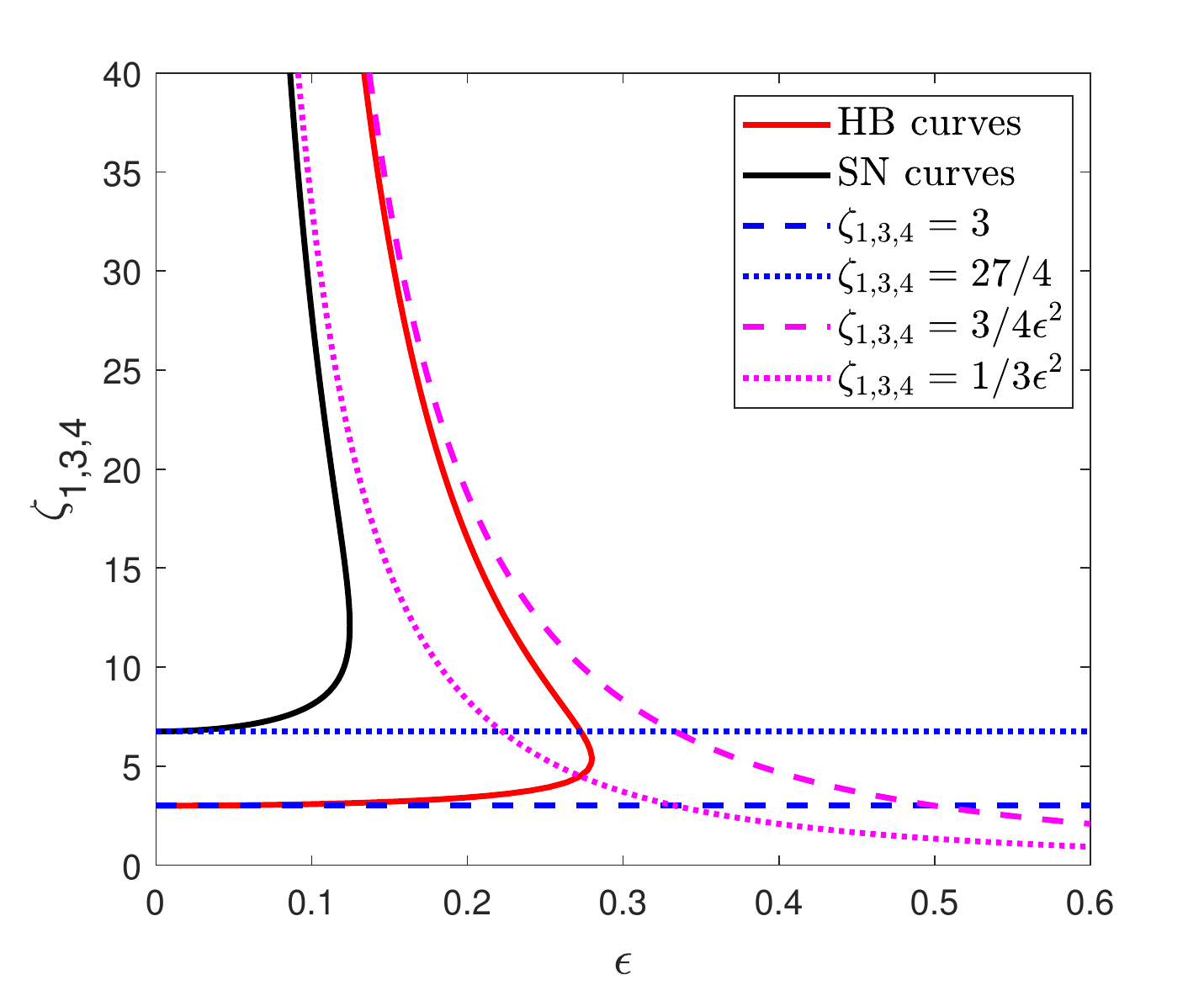}
\caption{Hopf and saddle-node bifurcations in the $(\epsilon,\zeta)$ parameter plane for the 4-node network example shown in Fig.~\ref{4node_network_example} with nonlinear oscillator at node 2 ($Q=2$). $\zeta_1=\zeta_3=\zeta_4=\mu$.}\label{4_node_HB_and_SN_curves_Q2}
\end{figure}

We may estimate the maximal value $\epsilon_{\text{max}}$ along this contour by considering the intersection between the corresponding asymptotes, scaled by an empirical factor of correction. From the data reported in Fig.~\ref{15_node_HB_and_SN_curves}, we obtain the intersection at $\epsilon\approx 0.29$, while $\epsilon_{\text{max}}\approx 0.13$. Similarly, from the data reported in Fig.~\ref{15_node_Q5_HB_and_SN_curves}, these values equal $\approx 0.64$ and $\approx 0.32$, respectively. In general, for a topology with all nonzero $P_{k,I}$ (and, therefore, $\zeta_i=\mu$ for $i\ne Q$), we obtain the equality
\begin{equation}\label{eq_upperlimit}
    (\nu+\eta/8)\epsilon =\sqrt{ \frac{(K_{Q,Q}-1)(1-P_{Q,I}^2)}{K_{Q,Q}P_{Q,I}^2}}
\end{equation}
at the intersection of the two asymptotes and divide the predicted value of $\epsilon$ by an empirical factor of $2$ to obtain an estimated upper bound for the possibility of hysteresis. These conditions apply to the 4-node network example with $Q=2$, for which we estimate $\epsilon_{\text{max}}\approx 1/9$. In this case, Fig.~\ref{4_node_HB_and_SN_curves_Q2} shows a single smooth contour composed from the saddle-node bifurcation curve associated with limit cycles with angular frequencies $\sqrt{5}$ and $2$ in the asymptotic limits $\mu=\mathcal{O}(1)$ and $\mu=\mathcal{O}(1/\epsilon^2)$, respectively. Our estimate for $\epsilon_{\text{max}}$ is clearly roughly accurate.

\begin{figure*}[htbp]
\centering
\subfloat[]{\includegraphics[width=.45\textwidth]{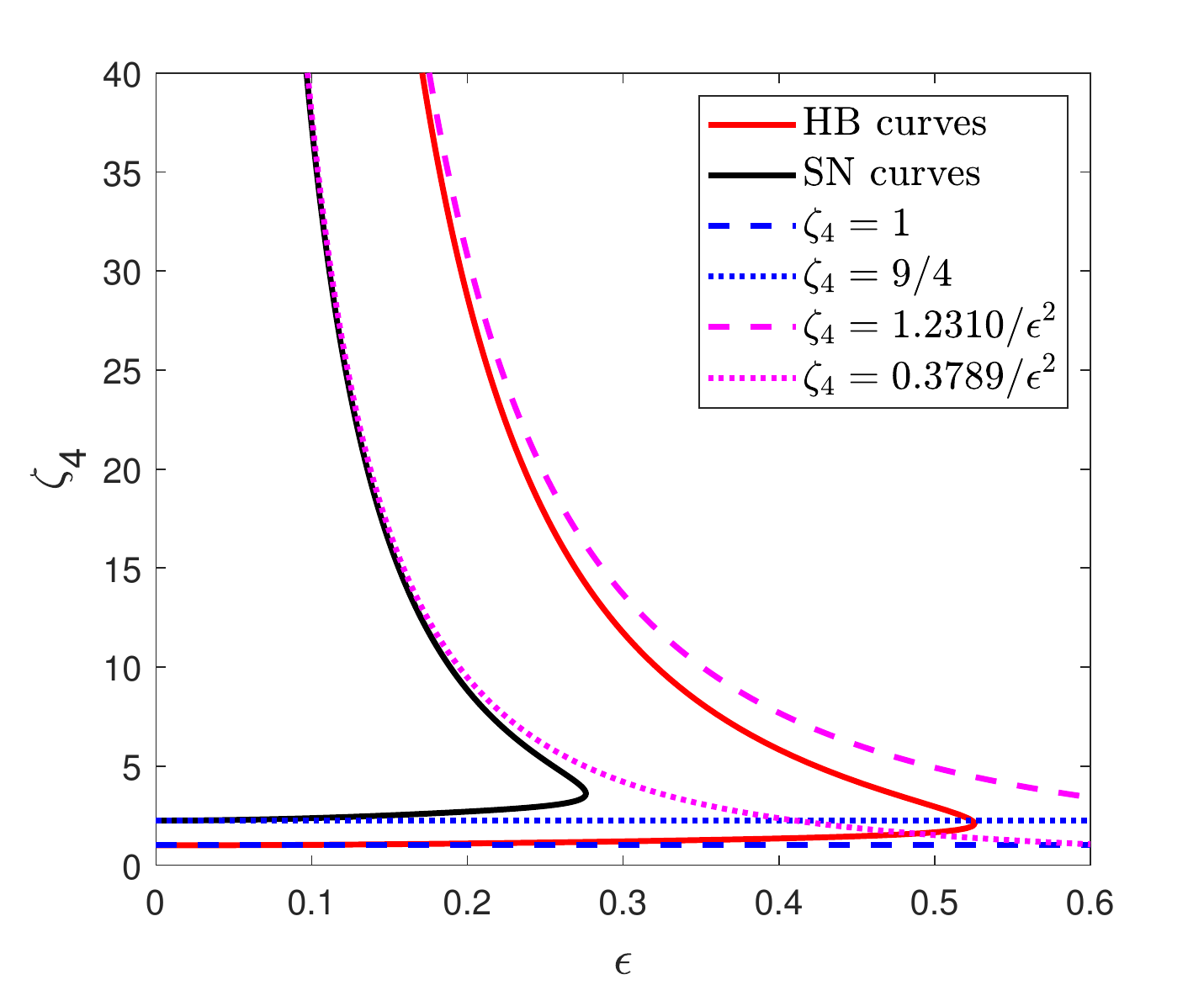}\label{4_node_HB_and_SN_curves_Q1_a}}
\subfloat[]{\includegraphics[width=.45\textwidth]{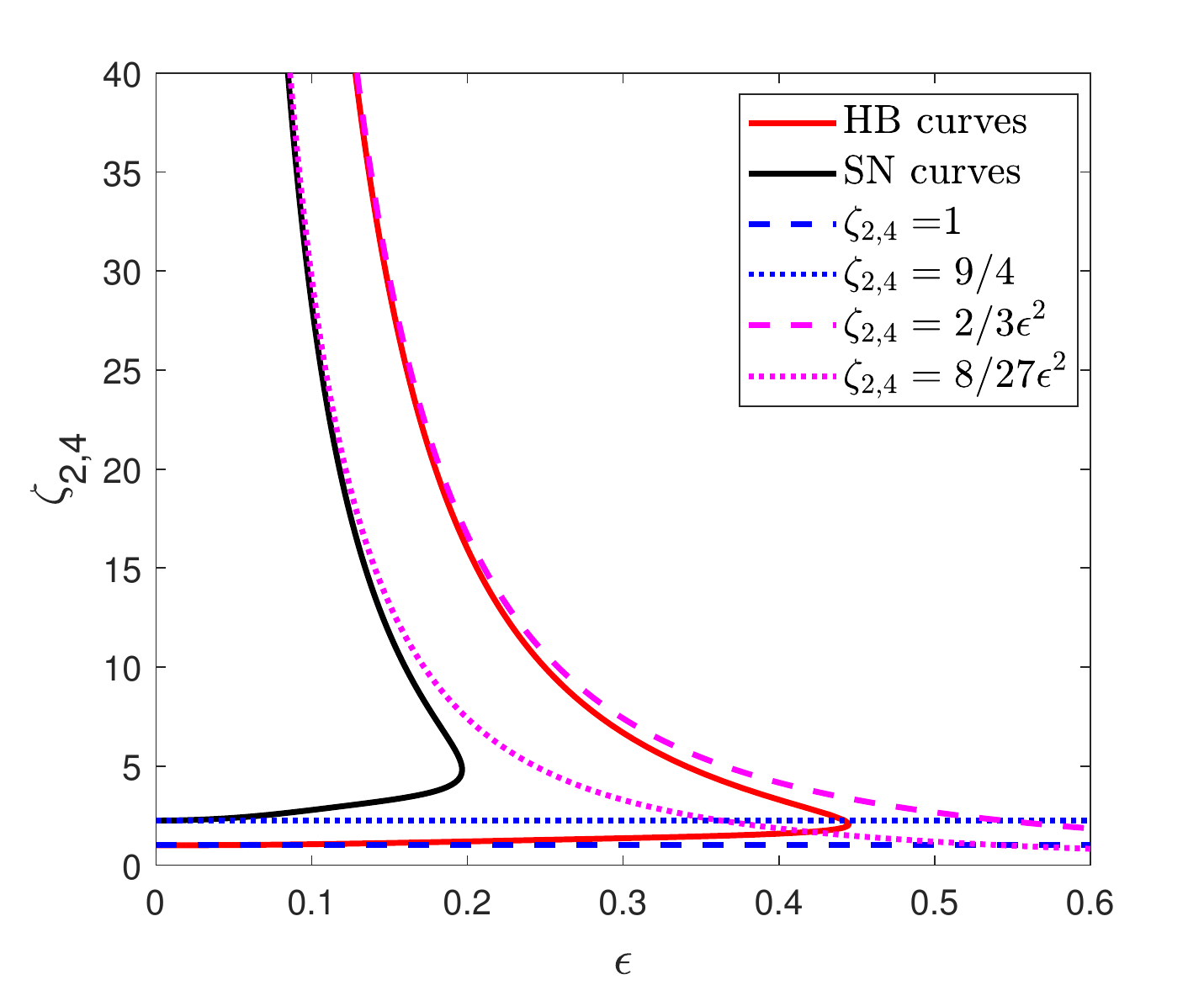}\label{4_node_HB_and_SN_curves_Q1_b}}
\caption{Hopf and saddle-node bifurcations in the $(\epsilon,\zeta)$ parameter plane for the 4-node network example shown in Fig.~\ref{4node_network_example} with nonlinear oscillator at node 1 ($Q=1$). Solid lines denote the continuation results, dashed lines denote the perturbation results for Hopf bifurcation, and dotted lines denote the perturbation results for saddle-node bifurcation. (a) $\zeta_2=\zeta_3=1+\delta$, (b)$\zeta_2=\zeta_4$ and $1/\epsilon^2\zeta_3=(3+\delta)/2$, $\delta=0.1$.}\label{4_node_HB_and_SN_curves_Q1}
\end{figure*}

\begin{figure*}[htbp]
\centering
\subfloat[]{\includegraphics[width=.45\textwidth]{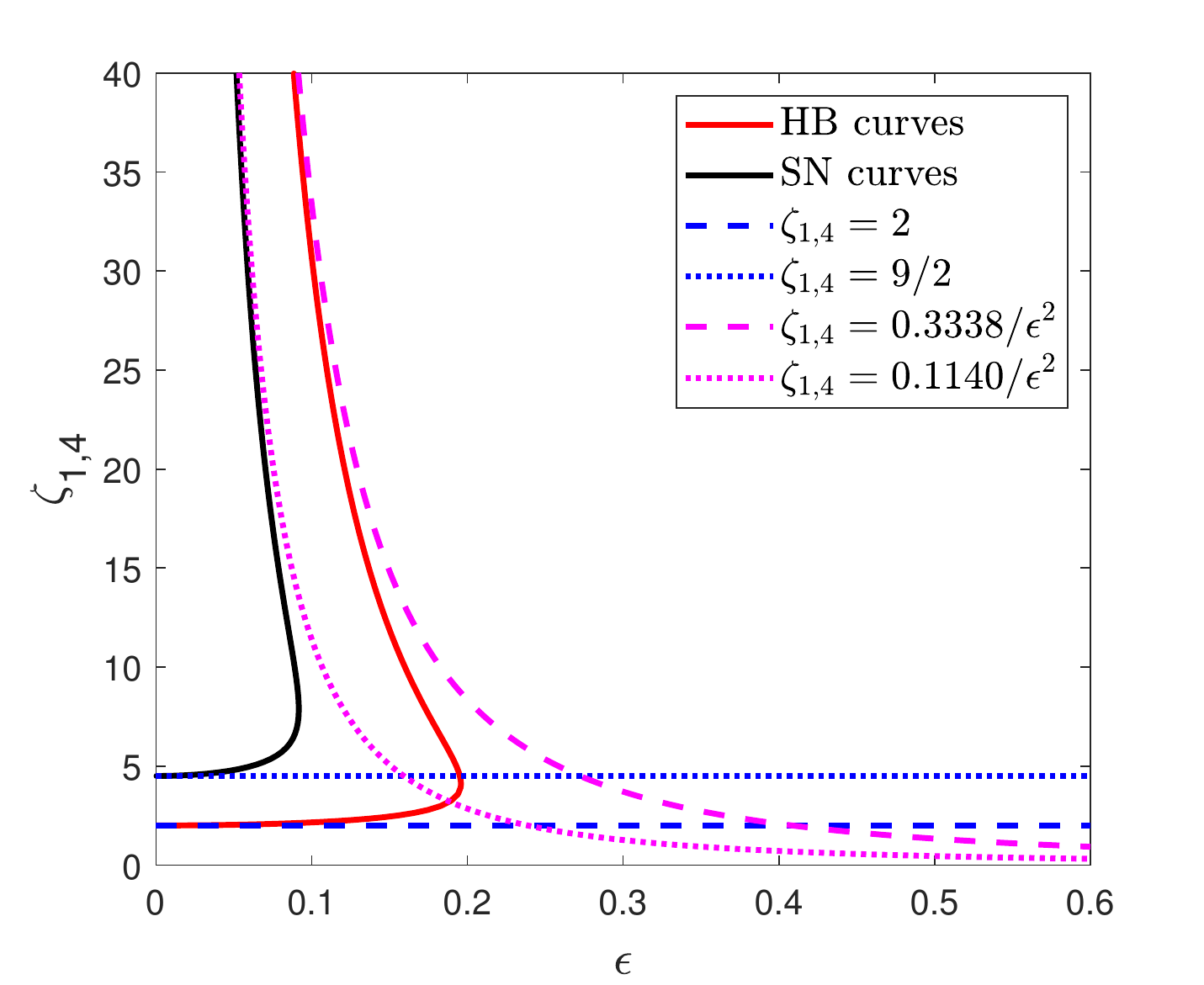}\label{4_node_HB_and_SN_curves_Q3_a}}
\subfloat[]{\includegraphics[width=.45\textwidth]{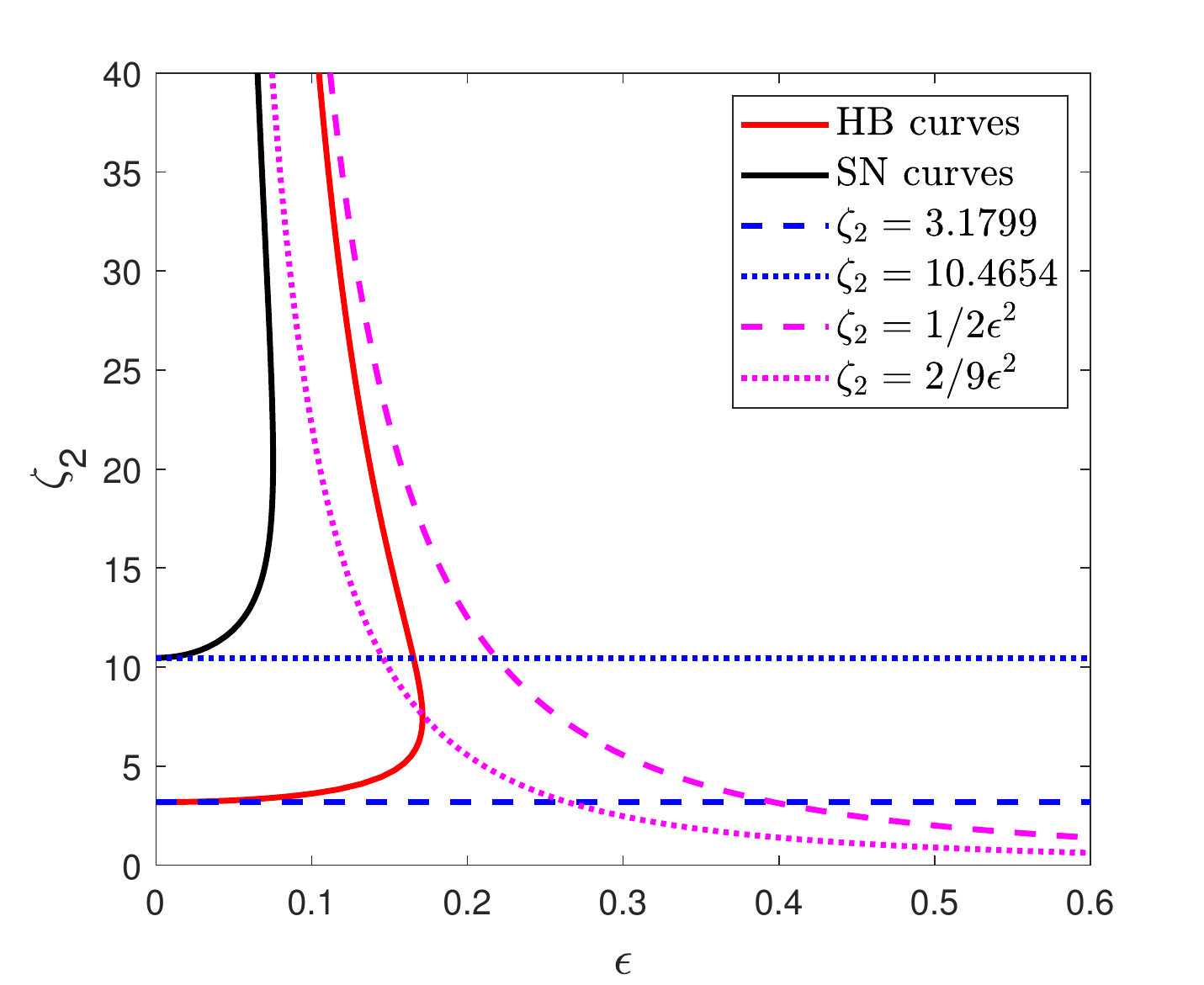}\label{4_node_HB_and_SN_curves_Q3_b}}
\caption{Hopf and saddle-node bifurcations in the $(\epsilon,\zeta)$ parameter plane for the 4-node network example shown in Fig.~\ref{4node_network_example} with nonlinear oscillator at node 3 ($Q=3$). (a) $\zeta_1=\zeta_4$ and $\zeta_2=(4+\delta)/2$, (b)$1/\epsilon^2\zeta_1=1/\epsilon^2\zeta_4=2+\delta$, $\delta=0.1$.}\label{4_node_HB_and_SN_curves_Q3}
\end{figure*}

For the cases with $Q=1$ and $Q=3$, the definition of $\mu$ depends on whether we are concerned with hysteretic behavior with $\zeta_i=\mathcal{O}(1)$ or $\zeta_i=\mathcal{O}(1/\epsilon^2)$ as $\epsilon\rightarrow 0$. For example, as discussed in Sect.~\ref{sec:a model problem}, for $Q=1$ and $\zeta_i=\mathcal{O}(1)$, we expect variations only in $\zeta_4$ along a hysteretic trajectory, while $\zeta_2=\zeta_3=1+\delta$ throughout. With $\zeta_4=\mu$ and $\zeta_2=\zeta_3=1.1$, we obtain the Hopf and saddle-node bifurcation curves shown in Fig.~\ref{4_node_HB_and_SN_curves_Q1_a} under simultaneous variations in $\epsilon$ and $\mu$. The horizontal asymptotes are identical to those obtained in Sect.~\ref{sec:a model problem}. In contrast, the hyperbolic asymptotes are here obtained by assuming that $\mu=\mathcal{O}(1/\epsilon^2)$ rather than $\zeta_i=\mathcal{O}(1/\epsilon^2)$ as in Sect.~\ref{sec:a model problem}. In this case, we estimate $\epsilon_{\text{max}}\approx 0.21$. If, instead, we consider the limit $\zeta_i=\mathcal{O}(1/\epsilon^2)$, then the analysis in Sect.~\ref{sec:a model problem} suggests identical variations in $\tilde{\zeta}_2$ and $\tilde{\zeta}_4$ along a hysteretic trajectory, while $\tilde{\zeta}_3=(3+\delta)/2$ throughout. Indeed, with $\zeta_2=\zeta_4=\mu$ and $\tilde{\zeta}_3=1.55$, we obtain the Hopf and saddle-node bifurcation curves shown in Fig.~\ref{4_node_HB_and_SN_curves_Q1_b} under simultaneous variations in $\epsilon$ and $\mu$. Our empirical estimate for $\epsilon_{\text{max}}$ yields $ 0.18$. Finally, for $Q=3$, we either let $\zeta_1=\zeta_4=\mu$ and $\zeta_2=(4+\delta)/2$ or $\zeta_2=\mu$ and $\tilde{\zeta}_1=\tilde{\zeta}_4=2+\delta$ to estimate an upper bound for the value of $\epsilon$ that supports hysteretic behavior for $\zeta_i=\mathcal{O}(1)$ or $\zeta_i=\mathcal{O}(1/\epsilon^2)$, respectively. We show the results of numerical continuation and the corresponding asymptotes predicted by the multiple-scale analysis in Fig.~\ref{4_node_HB_and_SN_curves_Q3}. In both cases, the empirical estimate for $\epsilon_{\text{max}}$ are roughly accurate.

It is clear from these results that the upper bound $\epsilon_{\text{max}}$ depends on the topology of the network. For example, since $K_{Q,Q}-1$ equals the degree of the $Q$-th node, the value of the estimated upper limit given in Eq.~\eqref{eq_upperlimit} depends on the network topology both directly through its connectivity and indirectly through the matrix $P$, whose columns span the eigenspace of the corresponding Laplacian. For a given topology, we may search among the possible values of $Q$ for that which yields the greatest estimated value of $\epsilon_{\text{max}}$, as this is likely to correspond to a greater range of bistability and possible hysteretic behavior under the coupled dynamics proposed in Sect.~\ref{sec:examples}. In the 15-node network example, the largest estimated value for $\epsilon_{\text{max}}$ is obtained for $Q=10$.

%%%%%%%%%%%%%%%%%%%%%%%%%%%%%%%%%%%%%%%%%%%%%%%%%%%%%%%%%%

\section{The triggering impulse}\label{sec_hysteresisresonance}

The examples in Sect.~\ref{sec:examples} used a brief pulse of harmonic excitation at the natural frequency associated with the critical Hopf bifurcation to excite a transition to self-sustained oscillations and the subsequent hysteretic return to quiescence. In this section, we consider briefly the relationship between the amplitude of excitation and the required duration of the impulse that will ensure such a transition.

To this end, consider again the case when $\zeta_i=\mathcal{O}(1)$ for all $i\ne Q$ and assume that $F(t)=\epsilon f\cos\omega_I t$ for some amplitude vector $f$. From the analogous perturbation analysis as in Sect.~\ref{sec_networkfilters}, we obtain secular terms in the dynamics of $\sum_i P_{i,I}v_i$ with coefficients
\begin{equation}
    -2 \omega_I A \phi' - \sum_i P_{i,I}f_i \cos\phi
\end{equation}
and
\begin{equation}
\begin{split}
  &-2\omega_I A'+\left( \nu  P_{Q,I}^2 - \zeta \right)\omega_I A +\frac{1}{4} P_{Q,I}^4 \eta \omega_I A^3 \\
  & - \frac{1}{8} P_{Q,I}^6 \eta \omega_I A^5 - \sum_i P_{i,I}f_i\sin\phi
\end{split}
\end{equation}
in front of $\cos(\omega_I t+\phi)$ and $\sin(\omega_I t+\phi)$, respectively. With $\zeta = \zeta_{\text{HB}}+\delta = P_{Q,I}^2\nu + \delta$, the predicted slow dynamics of $A$ and $\phi$ are then governed by the differential equations
\begin{equation}\label{eq_Nforcedslowdynamics}
\begin{split}
A' =&  -\frac{A}{2} \delta  +\frac{1}{8} P_{Q,I}^4 \eta A^3 - \frac{1}{16} P_{Q,I}^6 \eta  A^5  \\
& - \frac{\sin\phi}{2\omega_I} \sum_i P_{i,I}f_i \\
\phi' =&  - \frac{\cos\phi}{2 A \omega_I} \sum_i P_{i,I}f_i
\end{split}
\end{equation}

From Eq.~\eqref{eq:Anull}, we see that a transition to self-excited oscillations follows the termination of the exogenous excitation provided that $A$ exceeds the threshold
\begin{equation}\label{eq_Nunstable}
    \bar{A}=\sqrt{\frac{1}{P^2_{Q,I}}-\sqrt{\frac{P^2_{Q,I}\eta-8\delta}{P^6_{Q,I}\eta}}}
\end{equation}
at that time. Since $\bar{A}=\mathcal{O}(\sqrt{\delta})$ for $\delta\ll 1$, we may assume that $A=\mathcal{O}(\sqrt{\delta})$ for the duration of the triggering impulse in this limit. In this case, we expect from Eq.~\eqref{eq_Nforcedslowdynamics} that $\phi$ will quickly become equal to $\pi$ or $-\pi$ independently of $\phi(0)$. Moreover, $A(0)=0$ implies that
\begin{equation}
    |A(t)| = \frac{\epsilon t}{2\omega_I}  \left|\sum_i P_{i,I}f_i\right| .
\end{equation}
The minimum time $t_{\text{req}}$ required to achieve a transition to self-excited oscillations is then estimated to equal
\begin{equation}\label{eq_treq_approximate}               \frac{2\omega_I\bar{A}}{\epsilon\left| \sum_i P_{i,I}f_i \right| } .
\end{equation}

\begin{figure}[htbp]
\centering
\includegraphics[width=.45\textwidth]{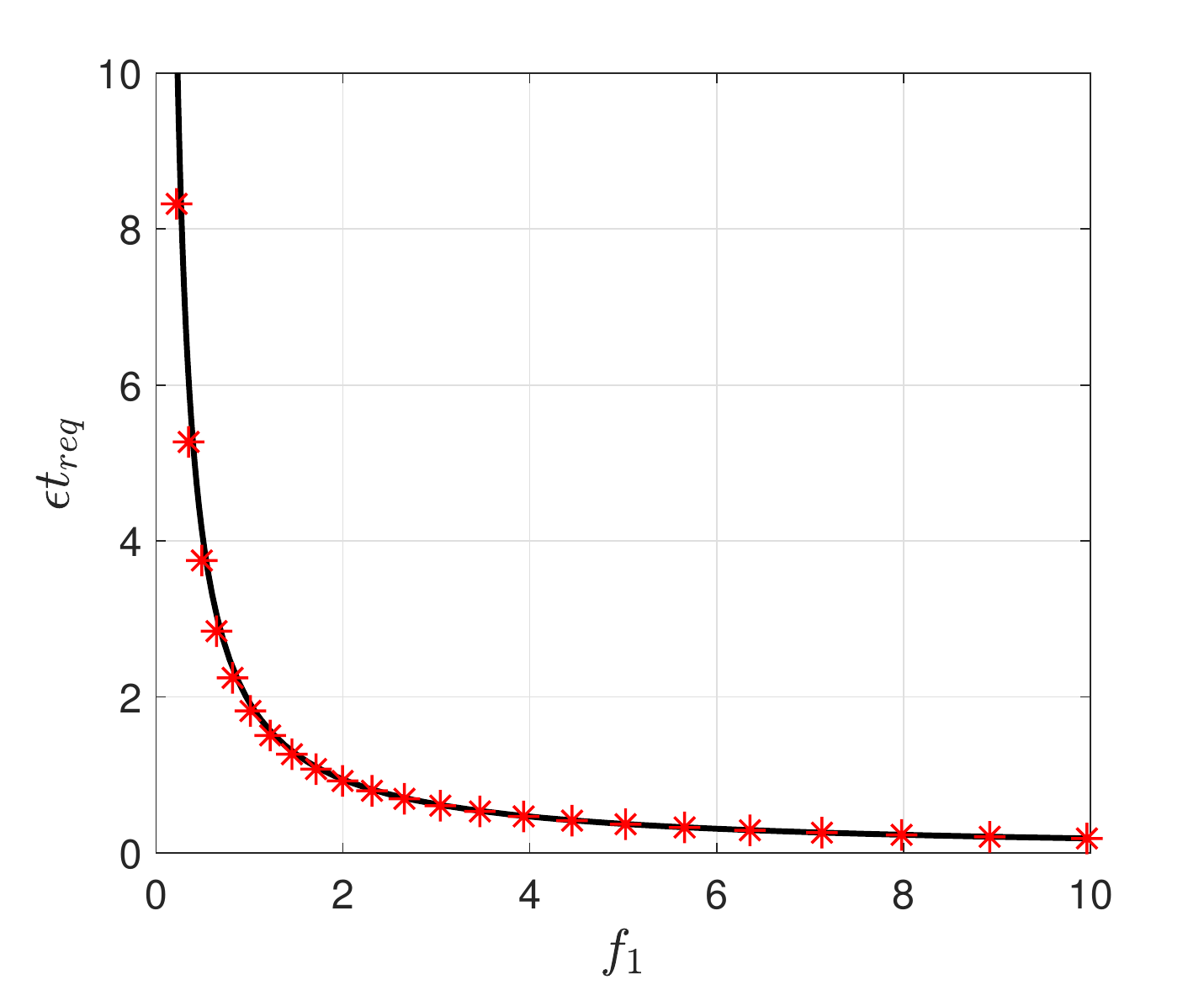}
\caption{Minimum required time for harmonic excitation on node 1 in order to trigger hysteresis in the 15-node network example shown in Fig.~\ref{15node_network_example}. Black line denotes the results obtained via numerical continuation, and red stars denote the results obtained from Eq.~\eqref{eq_treq_approximate}. $f_i=0$ for $i\neq 1$, $\delta=0.2$, $\mu_1=-1$, $\eta_1=10$, $\mu_i=\zeta_{\text{HB}}+\delta$ for $i\neq 1$. }\label{15_node_example_treq}
\end{figure}

Figure~\ref{15_node_example_treq} graphs the estimate Eq.~\eqref{eq_treq_approximate}  for $t_{\text{req}}$ against $f_1$ for the 15-node network example with $Q=1$ and $f_i=0$ for $i\ne Q$. We include additionally the predictions of numerical continuation obtained by solving a two-point boundary-value problem for the differential equation Eq.~\eqref{eq_Nforcedslowdynamics} with $A(0)=0$ and $A(t_{\text{req}})=\bar{A}$ for $\delta=0.2$. We observe close agreement between these predictions. We thus predict that, given an arbitrary network topology, the required duration of the triggering impulse may be minimized by choosing $Q$ so as to maximize $P_{Q,I}^2=\max_k P_{Q,k}^2$.  

%%%%%%%%%%%%%%%%%%%%%%%%%%%%%%%%%%%%

\section{Concluding discussion} \label{sec_conclusion}

A central goal of this work has been to design an oscillator architecture that exhibits a transient, hysteretic response to an external trigger. For this purpose, we have relied on the coupling of linear and nonlinear oscillators with a proposed definition of distributed internal dynamics. In the adiabatic limit, these guide a bifurcation parameter past a cyclic fold bifurcation and back to its initial value near a subcritical Hopf bifurcation. The analysis has shown the robustness of this behavior over some range of model parameter values away from the adiabatic limit. It has also uncovered a complicated bifurcation structure of stable and unstable periodic orbits associated with the network topology.

Our interest in the general phenomenology is inspired by the possibility that active processes in individual agents in a complex system may create conditions for collective, sustained action also in response to short, localized excitation. In a colony of social insects, localized disturbances may be amplified through a network of interactions to ensure a colony-wide change, also after the disturbance has been removed. In a social community, triggering information of particular resonance may result in sustained activity, modulated by internal variables governing attention spans. We imagine that qualitatively similar dynamic models may be proposed for such contexts, albeit with different underlying structure. 

A consequence of operation near a subcritical Hopf bifurcation is the possibility of noise-induced transitions to self-excited oscillations and subsequent return to quiescence. Were such unexpected pulses of spontaneous activity to occur in a physical, biological, or social system, it is reasonable to assume a similar mechanism of operating close to an explosive transition. Such operation has been used in the literature to explain enhanced frequency selectivity~\cite{wang2017explosive} in cochlea. Provided that the system has a way to return to quiescence, such spontaneous pulses may simply serve to keep the wheels greased, so to speak.

There are several opportunities for further work on the class of models discussed in this paper. These include alternative definitions of the internal dynamics of the damping coefficients (or other model parameters available for slow modulation), a rigorous derivation of topology-dependent upper bounds on $\epsilon$ that support hysteresis, and analysis of networks with additional nonlinear nodes.

As an example, our preliminary analysis of the persistence of hysteresis for finite values of $\epsilon$ has uncovered two distinct possibilities differentiated by the relative configurations of contours of cyclic fold bifurcations. We do not have a systematic theory for when either will occur in a particular network and for a particular location of a nonlinear node, nor can we exclude more complicated scenarios associated with closely spaced natural frequencies. In certain cases, we have been able to estimate an upper bound for $\epsilon$ using properties of the modal matrix $P$. This analysis was used to propose optimal locations for the nonlinear node that would maximize this upper bound. While our network model only considers unit coupling between nodes, networks with nonidentical coupling strengths expand the design possibilities and are worth further study. 

Furthermore, although it may be attractive to consider distributions of nonlinearity across multiple nodes of a network, we anticipate that this will be accompanied by transitions from periodic to quasiperiodic self-excited oscillations in the asymptotic limit of large damping for the linear oscillators (cf.~\cite{storti1982dynamics}). Whether anything would be gained by such complexity remains to be determined.

Finally, it remains of great interest to explore integration of the proposed phenomenology in a physical sensor network, for example, for mass sensing by microcantilevers. As part of ongoing work and inspired by \cite{KaYaYaMa2020}, we envision a hybrid realization that combines physical and in-silico components and that realizes internodal coupling through suitable actuator interfaces. For any such implementation, additional considerations of energy consumption and complexity would invariably drive further developments.

%% For one-column wide figures use
%\begin{figure}
%% Use the relevant command to insert your figure file.
%% For example, with the graphicx package use
%  \includegraphics{example.eps}
%\end{figure}
%%
%% For two-column wide figures use
%\begin{figure*}
%% Use the relevant command to insert your figure file.
%% For example, with the graphicx package use
%  \includegraphics[width=0.75\textwidth]{example.eps}
%\end{figure*}
%\begin{table}
%\begin{tabular}{lll}
%\hline\noalign{\smallskip}
%\end{tabular}
%\end{table}

\begin{acknowledgements}
This material is based upon work supported by the National Science Foundation under Grant No. BCS-1246920. 
\end{acknowledgements}

\section*{Conflict of interest}
The authors declare that they have no conflict of interest.

\section*{Data availability}
Data sharing not applicable to this article as no datasets were generated or analysed during the current study.

\begin{appendices}

\section{Nomenclature}

\setlist[itemize,1]{label=,labelwidth=0.9in,align=parleft,itemsep=0.1\baselineskip,leftmargin=!}

\begin{itemize}
\item[$\delta$] A small offset for damping dynamics 
\item[$\epsilon$] Scaling parameter of linear and nonlinear damping in the network model
\item[$\eta$] Nonlinear damping parameter of the nonlinear oscillator in the network model
\item[$\mu$] Damping parameter for linear oscillator at each node when $\zeta_k=\mu$ for all $k$
\item[$\nu$] Linear damping parameter of the nonlinear oscillator in the network model
\item[$\omega_k$] Natural frequency of $k$-th mode in the network model 
\item[$\phi(\epsilon t)$] Phase of $\mathcal{O}(1)$ solution in multiple-scale analysis
\item[$\tau$] A large time scale for damping dynamics
\item[$\tilde{\zeta}$] A damping parameter in bifurcation analysis, $\tilde{\zeta}=\sum_{k\ne Q}\frac{K_{Q,k}K_{k,Q}}{\epsilon^2\zeta_k K_{Q,Q}}$
\item[$\tilde{C}(x)$] Modal damping matrix
\item[$\tilde{K}$] Modal stiffness matrix
\item[$\zeta$] A damping parameter in bifurcation analysis, $\zeta=\sum_{k\ne Q}P^2_{k,I}\zeta_k$\\
\item[$\zeta_k$] Damping parameter for linear oscillator at node $k$ in the network model
\item[$\zeta_{\text{HB}}$, $\tilde{\zeta}_{\text{HB}}$, $\mu_{\text{HB}}$] Damping parameter at Hopf bifurcation
\item[$\zeta_{\text{SN}}$, $\tilde{\zeta}_{\text{SN}}$, $\mu_{\text{SN}}$] Damping parameter at saddle-node bifurcation
\item[$A(\epsilon t)$] Amplitude of $\mathcal{O}(1)$ solution in multiple-scale analysis\\
\item[$A_k(t)$] Amplitude of node $k$ in the network
\item[$C(u)$] Damping matrix of the network model
\item[$f$] Amplitude vector for exogenous excitation
\item[$F(t)$] Time-dependent exogenous excitation vector
\item[$I_N$] $N\times N$ identity matrix
\item[$K$] Stiffness matrix of the network model, determined by network topology
\item[$L$] The Laplacian matrix of network
\item[$P$] An orthogonal matrix whose columns span the eigenspace of the network Laplacian
\item[$Q$] The index of the nonlinear oscillator node in the network model
\item[$u(t)$] Nodal displacement vector in the network model
\item[$x$] Modal coordinate vector
\end{itemize}

\section{Multiple-scales perturbation analysis for the 4-node network with small linear damping}\label{appendix_4node_multiplescale}

In this appendix, we present a detailed derivation of the slow-flow dynamics in Eqs.~\eqref{eq_4nodeslowflow1} and \eqref{eq_4nodeslowflow2} in the limit of $\epsilon\ll 1$ using the method of multiple scales \cite{nayfeh1995nonlinear}. For the case of constant $\zeta_2,\zeta_3,\zeta_4\sim\mathcal{O}(1)$, we seek to arrive at a zeroth-order (in $\epsilon$) description of the displacement vector $u$ that mirrors the form of the corresponding linear normal-mode, albeit with slowly varying amplitude and phase. To this end, we assume two time scales $t_0=t$ and $t_1=\epsilon t$ and, per Eq.~\eqref{eq_4nodemultscale}, a response of the form
\begin{align}
    u&=\frac{A(t_1)}{\sqrt{2}}\begin{pmatrix}1\\0\\0\\-1\end{pmatrix}\cos\left(2t_0+\phi(t_1)\right)+\epsilon v(t_0,t_1)
\end{align}
and compute derivatives using the relationship
\begin{equation}
    \frac{d}{dt}=\frac{\partial}{\partial t_0}+\epsilon\frac{\partial}{\partial t_1}.
\end{equation}
Using this notation, we obtain
\begin{align}
    \dot{u}_1&=-\sqrt{2}A(t_1)\sin\left(2t_0+\phi(t_1)\right)\nonumber\\
    &\qquad+\epsilon\frac{A'(t_1)}{\sqrt{2}}\cos\left(2t_0+\phi(t_1)\right)\nonumber\\
    &\qquad-\epsilon\phi'(t_1)\frac{A(t_1)}{\sqrt{2}}\sin\left(2t_0+\phi(t_1)\right)\nonumber\\
    &\qquad+\left(\frac{\partial}{\partial t_0}+\epsilon\frac{\partial}{\partial t_1}\right)\epsilon v_1(t_0,t_1)
\end{align}
and so on. Substitution into the governing equation \eqref{eq_4nodegenform} with $F=0$ results in the perfect cancellation of all $\mathcal{O}(1)$ terms. A similar cancellation of all terms of $\mathcal{O}(\epsilon)$ requires that
\begin{equation}\label{eq_oeps}
    \frac{\partial^2 v}{\partial t_0^2} + K v = q(t_0,t_1)
\end{equation}
where
\begin{align}
    q_1 &= \frac{5 A(t_1)^5 }{8 \sqrt{2}}s_1+\frac{15 A(t_1)^5 }{16 \sqrt{2}}s_3 +\frac{5 A(t_1)^5 }{16 \sqrt{2}}s_5\nonumber\\ 
    &\qquad-\frac{5 A(t_1)^3 }{2 \sqrt{2}}s_1  -\frac{5 A(t_1)^3 }{2 \sqrt{2}}s_3+2 \sqrt{2} A(t_1) \phi' c_1\nonumber\\
    &\qquad-\sqrt{2} A(t_1) s_1+2 \sqrt{2} A'(t_1) s_1 \\
    q_2&= q_3=0 \\
    q_4&= -\sqrt{2} A(t_1) \zeta_4 s_1-2 \sqrt{2} A(t_1) \phi' c_1\nonumber\\
    &\qquad-2 \sqrt{2} A'(t_1) s_1,
\end{align}
where we use the shorthand $c_k=\cos\left(2kt_0+k\phi(t_1)\right)$ and $s_k=\sin\left(2kt_0+k\phi(t_1)\right)$. It follows that
\begin{align}
    \frac{\partial^2 w}{\partial t_0^2} &+ 4w = \frac{5}{16} A^5(t_1) s_1+\frac{15}{32} A^5(t_1) s_3+\frac{5}{32} A^5(t_1) s_5\nonumber\\&\qquad-\frac{5}{4} A^3(t_1) s_1-\frac{5}{4} A^3(t_1) s_3+A(t_1) \zeta_4 s_1\nonumber\\&\qquad
    +4 A(t_1) \phi'(t_1) c_1-A(t_1) s_1+4 A'(t_1) s_1,
\end{align}
where $w=(v_1-v_4)/\sqrt{2}$. For a consistent approximation, the \emph{secular terms} on the right-hand side proportional to $\cos(2t_0+\phi(t_1))$ and $\sin(2t_0+\phi(t_1))$, respectively, must cancel, resulting in the conditions
\begin{equation}
    4 A(t_1)\phi'(t_1)=0
\end{equation}
and
\begin{equation}
    4 A'(t_1)-(1-\zeta_4)A(t_1)-\frac{5}{4}A^3(t_1)+\frac{5}{16}A^5(t_1)=0 .
\end{equation}

\section{4-node network with small linear damping and $Q=2$, $3$, and $4$, respectively}\label{appendix_4node_smalldamping}

The existence of a subcritical Hopf bifurcation and an associated branch of limit cycles found in the case of $Q=1$ carries over also to other configurations of the network with the nonlinear oscillator located at node 2, 3, or 4. Indeed, when the nonlinear oscillator is located at node 4 with $\zeta_1,\zeta_2,\zeta_3\sim\mathcal{O}(1)$, the network symmetry results in an identical set of results with $\zeta_4$ replaced by $\zeta_1$. In contrast, if the nonlinear oscillator is located at node 2 with $\zeta_1,\zeta_3,\zeta_4\sim\mathcal{O}(1)$, the complex exponential rates in Eq.~\eqref{eq_4nodemultscale} become
\begin{align}
\begin{split}
    \lambda_{1,2}&=\pm\mathrm{j}+\frac{1-\zeta_1-\zeta_3-\zeta_4}{8}\epsilon+\mathcal{O}(\epsilon^2),\\
    \lambda_{3,4}&=\pm\sqrt{2}\mathrm{j}-\frac{\zeta_1+4\zeta_3+\zeta_4}{12}\epsilon+\mathcal{O}(\epsilon^2),\\
    \lambda_{5,6}&=\pm2\mathrm{j}-\frac{\zeta_1+\zeta_4}{4}\epsilon+\mathcal{O}(\epsilon^2),\\
    \lambda_{7,8}&=\pm\sqrt{5}\mathrm{j}+\frac{9-\zeta_1-\zeta_3-\zeta_4}{24}\epsilon+\mathcal{O}(\epsilon^2).
\end{split}
\end{align}
In this case, with $\zeta:=\zeta_1+\zeta_3+\zeta_4$, $u=0$ is asymptotically stable for $\zeta>9$ and unstable for $\zeta<9$. Substitution of the ansatz
\begin{equation}
\label{eq_4nodemultscale_2}
\begin{split}
    &u_1(t)=\frac{1}{2\sqrt{3}}A(\epsilon t)\cos\left( \sqrt{5}t+\phi(\epsilon t )\right)+\epsilon v_1(t),\\
    &u_2(t)=-\frac{\sqrt{3}}{2}A(\epsilon t)\cos\left( \sqrt{5}t+\phi(\epsilon t )\right)+\epsilon v_2(t),\\
    &u_3(t)=\frac{1}{2\sqrt{3}}A(\epsilon t)\cos\left( \sqrt{5}t+\phi(\epsilon t )\right)+\epsilon v_3(t),\\
    &u_4(t)=\frac{1}{2\sqrt{3}}A(\epsilon t)\cos\left( \sqrt{5}t+\phi(\epsilon t )\right)+\epsilon v_4(t)
\end{split}
\end{equation}
into the fully nonlinear governing equations then yields the differential equation
\begin{equation}
    A'=\frac{1}{24}(9-\zeta)A+\frac{45}{64}A^3-\frac{135}{512}A^5
\end{equation}
governing the slow dynamics of the amplitude $A$ with nonzero equilibria obtained from
\begin{equation}
    A^2=\frac{4}{3}\pm \frac{4}{9}\sqrt{\frac{81-4\zeta}{5}},
\end{equation}
i.e., for $\zeta\in[0,81/4]$ with two co-existing solutions on the interval $\zeta\in[9,81/4)$. We obtain the equivalent of Eq.~\eqref{eq_4node_damplaw1}, for example, by letting
\begin{equation}
    \zeta_{1,3,4}'=\tau^{-1}\left(-\zeta_{1,3,4}+\frac{1}{3}\delta+3+\frac{45}{16}A^2\right),
\end{equation}
since these imply that
\begin{equation}
\zeta'=\tau^{-1}\left(-\zeta+\delta+9+\frac{135}{16}A^2\right).
\end{equation}

Finally, with the nonlinear oscillator located at node 3 with $\zeta_1,\zeta_2,\zeta_4\sim\mathcal{O}(1)$, the complex exponential rates in Eq.~\eqref{eq_4nodemultscale} become
\begin{equation}
\begin{split}
    \lambda_{1,2}&=\pm\mathrm{j}+\frac{1-\zeta_1-\zeta_2-\zeta_4}{8}\epsilon+\mathcal{O}(\epsilon^2),\\
    \lambda_{3,4}&=\pm\sqrt{2}\mathrm{j}+\frac{4-\zeta_1-\zeta_4}{12}\epsilon+\mathcal{O}(\epsilon^2),\\
    \lambda_{5,6}&=\pm2\mathrm{j}-\frac{\zeta_1+\zeta_4}{4}\epsilon+\mathcal{O}(\epsilon^2),\\
    \lambda_{7,8}&=\pm\sqrt{5}\mathrm{j}+\frac{1-\zeta_1-9\zeta_2-\zeta_4}{24}\epsilon+\mathcal{O}(\epsilon^2).
\end{split}
\end{equation}
In this case, with $\zeta:=\zeta_1+\zeta_4$, $u=0$ is asymptotically stable for $\zeta>4$ and unstable for $\zeta<4$. Substitution of the ansatz
\begin{equation}
\label{eq_4nodemultscale_3}
\begin{split}
    &u_1(t)=\frac{1}{\sqrt{6}}A(\epsilon t)\cos\left( \sqrt{2}t+\phi(\epsilon t )\right)+\epsilon v_1(t),\\
    &u_2(t)=\epsilon v_2(t),\\
    &u_3(t)=-\sqrt{\frac{2}{3}}A(\epsilon t)\cos\left( \sqrt{2}t+\phi(\epsilon t )\right)+\epsilon v_3(t),\\
    &u_4(t)=\frac{1}{\sqrt{6}}A(\epsilon t)\cos\left( \sqrt{2}t+\phi(\epsilon t )\right)+\epsilon v_4(t)
\end{split}
\end{equation}
into the fully nonlinear governing equations then yields the differential equation
\begin{equation}
    A'=\frac{1}{12}(4-\zeta)A+\frac{5}{9}A^3-\frac{5}{27}A^5
\end{equation}
governing the slow dynamics of the amplitude $A$ with nonzero equilibria obtained from
\begin{equation}
    A^2=\frac{3}{2}\pm \frac{3}{2}\sqrt{\frac{9-\zeta}{5}},
\end{equation}
i.e., for $\zeta\in[0,9]$ with two co-existing solutions on the interval $\zeta\in[4,9)$. We obtain the equivalent to Eq.~\eqref{eq_4node_damplaw1}, for example, by letting
\begin{equation}
\label{eq:damp_node3}
    \zeta_{1,4}'=\tau^{-1}\left(-\zeta_{1,4}+\frac{1}{2}\delta+2+\frac{5}{3}A^2\right),
\end{equation}
since these imply that
\begin{equation}
    \zeta'=\tau^{-1}\left(-\zeta+\delta+4+\frac{10}{3}A^2\right).
\end{equation}

\section{4-node network with large linear damping and $Q=2$, $3$, and $4$, respectively}\label{appendix_4node_largedamping}

When the nonlinear oscillator is located at node 4 with $\zeta_1,\zeta_2,\zeta_3\sim\mathcal{O}(1/\epsilon^2)$, the analysis from Sect.~\ref{sec:large linear damping} still applies, albeit with $\zeta_4$ and $\tilde{\zeta}_4$ replaced by $\zeta_1$ and $\tilde{\zeta}_1$, respectively. If, instead, the nonlinear oscillator is located at node 2 with $\zeta_1,\zeta_3,\zeta_4\sim\mathcal{O}(1/\epsilon^2)$, we obtain the complex exponential rates
\begin{equation}
\begin{split}
    &\lambda_{1,2,3}\approx -\epsilon\zeta_{1,3,4}+\mathcal{O}(1),\\
    &\lambda_{4,5,6}\lesssim 0, \\
    &\lambda_{7,8}=\pm2\mathrm{j}+\frac{1}{8}\left(4-\tilde{\zeta}_1-\tilde{\zeta}_3-\tilde{\zeta}_4\right)\epsilon+\mathcal{O}(\epsilon^2),
\end{split}
\end{equation}
where $\tilde{\zeta}_1:=1/\epsilon^2\zeta_1$, $\tilde{\zeta}_3:=1/\epsilon^2\zeta_3$, and $\tilde{\zeta}_4:=1/\epsilon^2\zeta_4$. With $\tilde{\zeta}:=\tilde{\zeta}_1+\tilde{\zeta}_3+\tilde{\zeta}_4$, it follows that the trivial equilibrium is asymptotically stable for $\tilde{\zeta}>4$ and unstable for $\tilde{\zeta}<4$, with $\tilde{\zeta}=4$ corresponding to a Hopf bifurcation out of which emanates a branch of periodic orbits approximated in the small-amplitude limit by the normal-mode oscillations $u_2=A\cos(2t+\phi)$, $u_1(t)=u_3(t)=u_4(t)=0$ for constant amplitude $A$ and phase $\phi$. A consistent multiple-scale ansatz now yields the amplitude equation
\begin{equation}
    A'=\frac{1}{8}\left(4-\tilde{\zeta}\right)A+\frac{5}{4}A^3-\frac{5}{8}A^5
\end{equation}
with nontrivial equilibria at
\begin{equation}
    A^2=1\pm\sqrt{\frac{9-\tilde{\zeta}}{5}},
\end{equation}
i.e., for $\tilde{\zeta}\in(0,9]$ with two co-existing solutions on the interval  $\tilde{\zeta}\in[4,9)$. We obtain the equivalent of Eq.~\eqref{eq_4node_damplaw1}, for example, by letting
\begin{equation}
    \tilde{\zeta}'_{1,3,4}=\tau^{-1}\left(-\tilde{\zeta}_{1,3,4}+\frac{1}{3}\delta+\frac{4}{3}+\frac{5}{3}A^2\right),
\end{equation}
since these imply that
\begin{equation}
\tilde{\zeta}'=\tau^{-1}\left(-\tilde{\zeta} + \delta + 4 + 5 A^2\right).
\end{equation}

Finally, when the nonlinear oscillator is located at node 3 with $\zeta_1,\zeta_2,\zeta_4\sim\mathcal{O}(1/\epsilon^2)$, we obtain the complex exponential rates
\begin{equation}
\begin{split}
    &\lambda_{1,2,3}\approx -\epsilon\zeta_{1,2,4}+\mathcal{O}(1),\\
    &\lambda_{4,5,6}\lesssim 0, \\
    &\lambda_{7,8}=\pm\mathrm{j}\sqrt{2}+\frac{1}{4}\left(2-\tilde{\zeta}_2\right)\epsilon+\mathcal{O}(\epsilon^2),
\end{split}
\end{equation}
where $\tilde{\zeta}_2:=1/\epsilon^2\zeta_2$. It follows that the trivial equilibrium is asymptotically stable for $\tilde{\zeta}_2>2$ and unstable for $\tilde{\zeta}_2<2$, with $\tilde{\zeta}_2=2$ corresponding to a Hopf bifurcation out of which emanates a branch of periodic orbits approximated in the small-amplitude limit by the normal-mode oscillations $u_3=A\cos(\sqrt{2}t+\phi)$, $u_1(t)=u_2(t)=u_4(t)=0$ for constant amplitude $A$ and phase $\phi$. A consistent multiple-scale ansatz now yields the amplitude equation
\begin{equation}
    A'=\frac{1}{4}\left(2-\tilde{\zeta}_2\right)A+\frac{5}{4}A^3-\frac{5}{8}A^5
\end{equation}
with nontrivial equilibria at
\begin{equation}
    A^2=1\pm\sqrt{\frac{9-2\tilde{\zeta}_2}{5}},
\end{equation}
i.e., for $\tilde{\zeta}_2\in(0,9/2]$ with two co-existing solutions on the interval  $\tilde{\zeta}_2\in[2,9/2)$. We obtain the equivalent of Eq.~\eqref{eq_4node_damplaw1}, for example, by letting
\begin{equation}
    \tilde{\zeta}'_{2}=\tau^{-1}\left(-\tilde{\zeta}_{2}+\delta+2+\frac{5}{2}A^2\right).
\end{equation}

\end{appendices}

\bibliographystyle{spmpsci}      % mathematics and physical sciences
\bibliography{references}   % name your BibTeX data base

%% Non-BibTeX users please use
%\begin{thebibliography}{}
%%
%% and use \bibitem to create references. Consult the Instructions
%% for authors for reference list style.
%\end{thebibliography}

\end{document}